\documentclass{amsart}

\usepackage[a4paper]{geometry}

\usepackage{float}

\usepackage[colorlinks=true]{hyperref}
\hypersetup{urlcolor=blue, citecolor=blue}
\usepackage{enumerate}
\usepackage{graphicx}
\usepackage{subcaption}

\usepackage{tikz}
\usepackage{pgfplots}
\usepackage{pgfplotstable}
\usetikzlibrary{calc}
\usepackage{tkz-euclide}
\usepackage{tkz-base}
\usetkzobj{all}

\newtheorem{thm}{Theorem}[section]
\newtheorem{lem}[thm]{Lemma}
\newtheorem{prop}[thm]{Proposition}

\newtheorem{rem}[thm]{Remark}

\def\R{\mathbb{R}}

\def\be{\begin{equation}}
\def\ee{\end{equation}}
\def\g{{\bf g}}
\def\n{{\boldsymbol n}}
\def\p{\partial}
\def\grad{\boldsymbol{\nabla}}
\def\div{\grad\cdot}

\def\O{\Omega}

\def\sig{\sigma}

\def\k{\kappa}
\def\a{\alpha}
\def\eps{\epsilon}
\def\x{{\boldsymbol x}}

\def\d{{\rm d}}
\def\ov#1{\overline{#1}}

\def\wh#1{\widehat{#1}}

\def\bphi{\boldsymbol{\phi}}

\def\bW{{\boldsymbol  W}}

\def\bs{{\boldsymbol  s}}
\def\bv{{\boldsymbol v}}
\def\bm{{\boldsymbol m}}
\def\by{{\boldsymbol y}}

\def\0{{\bf 0}}
\def\1{{\bf 1}}

\def\bsigma{{\boldsymbol  \sigma}}

\def\bbK{{\mathbb K}}

\def\Aa{\mathcal{A}}
\def\bAa{\boldsymbol{\Aa}}

\def\Ee{\mathcal{E}}

\def\Hh{\mathcal{H}}

\def\Mm{\mathcal{M}}

\def\Xx{\mathcal{X}}
\def\bXx{{\boldsymbol{\Xx}}}

\def\bDelta{{\boldsymbol \Delta}}
\def\bPsi{{\boldsymbol \Psi}}
\def\bpi{{\boldsymbol{\pi}}}

\def\Eee{\mathfrak{E}}

\def\Nnn{\mathfrak{N}}

\def\Ttt{\mathfrak{T}}

\DeclareMathOperator*{\argmin}{argmin}

\def\bp{{\boldsymbol{p}}}

\def\clem#1{{#1}}

\begin{document}

\title[Simulation of multiphase porous media flows]{Simulation of multiphase porous media flows with minimizing movement and finite volume schemes}
\author{C. Canc\`es}\address{Cl\'ement Canc\`es:  Inria, Univ. Lille, CNRS, UMR 8524 - Laboratoire Paul Painlev\'e, F-59000 Lille (\href{mailto:clement.cances@inria.fr}{\tt clement.cances@inria.fr})}
\author{T. O. Gallouet}\address{Thomas Gallou\"et: Inria, Project team  MOKAPLAN and Mathematics department, Universit\'e de Li\`ege, Belgium, (\href{mailto:thomas.gallouet@inria.fr}{\tt thomas.gallouet@inria.fr})}   
\author{M. Laborde}\address{Maxime Laborde: Department of Mathematics and Statistics, McGill University, Montreal, CANADA (\href{mailto:maxime.laborde@mcgill.ca}{\tt maxime.laborde@mcgill.ca})}
\author{L. Monsaingeon}\address{L\'eonard Monsaingeon: IECL Universit\'e de Lorraine, Nancy, FRANCE \& GFM Universidade de Lisboa, Lisbon, PORTUGAL (\href{mailto:leonard.monsaingeon@univ-lorraine.fr}{\tt leonard.monsaingeon@univ-lorraine.fr})}

\begin{abstract}
The Wasserstein gradient flow structure of the PDE system governing 
multiphase flows in porous media was recently highlighted in [C. Canc\`es, T. O. Gallou\"et, and L. Monsaingeon, {\it Anal. PDE} 10(8):1845--1876, 2017].
The model can thus be approximated by means of the minimizing movement (or JKO) scheme, 
that we solve thanks to the ALG2-JKO scheme proposed in [J.-D. Benamou, G. Carlier, and M. Laborde, {\it ESAIM Proc. Surveys}, 57:1--17, 2016]. 
The numerical results are compared to a classical upstream mobility Finite Volume scheme, for which 
strong stability properties can be established. 
\end{abstract}

\maketitle

\noindent
{\small {\bf Keywords.}
Multiphase porous media flows; 
Wasserstein gradient flow; minimizing movement scheme; Augmented Lagrangian method; Finite Volumes
\vspace{5pt}

\noindent
{\bf AMS subjects classification. }
35K65, 35A15, 49M29, 65M08, 76S05
}

\tableofcontents

\section{Multiphase porous media flows as Wasserstein gradient flow}\label{sec:intro}

Because of their wide range of interest in the applications, multiphase flows in porous media 
have been the object of countless scientific studies.
In particular, there has been an 
extensive effort in order to develop reliable and efficient tools for \clem{the} simulation of such flows. 
\clem{In many practical situations, the characteristic size of the pores (typically of the order the 
$\mu$m for regular sandstones) is much smaller than the characteristic size of the domain 
of interest. The direct numerical simulation of fluid flows at the pore scale is therefore not tractable. 
The use of homogenized models of Darcy type is therefore commonly used to simulate porous media flows. 
The derivation of such models is the purpose of a very extended literature. We refer for instance 
to~\cite{BB90} for an extended introduction to the modeling of porous media flows. But let us stress that, as far as we know, 
there is no rigorous mathematical derivation of homogenized models for multiphase porous media flows.}

\clem{
Because of the very large friction of the fluid with the porous matrix, the energy is dissipated and 
inertia is often naturally neglected in the Darcy type models. The resulting models therefore have 
a formal gradient flow structure, as highlighted in~\cite{CGM15} for immiscible incompressible multiphase 
porous media flows. This was then rigorously established in~\cite{CGM17} that the equations governing 
such flows can be reinterpreted as a gradient flow in some appropriate Wasserstein space.}
%
The goal of this paper is to explore how this new point of view can be used to simulate multiphase
flows in porous media.

\subsection{Incompressible immiscible multiphase flows}\label{ssec:model}

As a first step, let us recall the equations governing multiphase porous media flows. 
We remain synthetic here and refer to the monograph~\cite{BB90} for a rather complete 
presentation of the models. 
The porous medium is represented by a convex bounded open subset $\O$ of $\R^d$ ($d\leq 3$).
Within this porous medium, $N+1$ phases are flowing. Denoting by 
$\bs = \left(s_0,\dots,s_N\right)$ the saturations, i.e., the volume 
ratios of the various phases in the fluid, the following total saturation relations has to be fulfilled:
\begin{subequations}\label{eq:EL}
\be\label{eq:cons}
s_0 + s_1 + \dots + s_N =1.
\ee
In what follows, we denote by 
\[
\bDelta = \left\{\bs \in \R_+^N\; \middle| \;    s_0 + s_1 + \dots + s_N =1 \;\right\},
\]
and by
\[
\bXx = \left\{\bs: \O \to \R^N \; \middle| \; \bs(\x) \in \bDelta \;\text{for a.e.}\; \x \in \O\;\right\}.
\]
As a consequence of~\eqref{eq:cons}, the composition of the fluid is fully characterized by the knowledge of 
\[
\bs^\ast = (s_1, \dots, s_N) \in \bDelta^* = \left\{(s_1, \dots, s_N) \in \R^N_+ \middle|\; \sum_{i=1}^N s_i \leq 1 \right\}.  \] 
Concerning the evolution, each phase is convected with its own speed 
\be\label{eq:conv}
\omega \p_t s_i + \div(s_i \bv_i) = 0, 
\ee
where $\omega$ stands for the porosity of the medium $\O$ and is assumed to 
be constant in the sequel for simplicity.
Then a straightforward rescaling in time allows 
to choose $\omega = 1$.
We further assume a no flux condition across the boundary $\p\O$ for each phase, 
hence the mass is conserved along time.
This motivates the introduction of the set 
\[
\bAa = \left\{ \bs \in L^1_+(\O)^N \; \middle| \; \int_\O s_i \d\x = \int_\O s_i^0 \d\x \right\}, 
\]
where $\bs^0=\left(s_i^0\right): \O \to \bDelta$ is a prescribed initial data.

The phase speeds $\bv_i$ are prescribed by the Darcy law~\cite{Darcy1856}
\be\label{eq:Darcy}
\bv_i = - \frac{\k}{\mu_i} \left(\grad p_i - \rho_i \g\right), \qquad i \in \{0,\dots, N\}. 
\ee
In~\eqref{eq:Darcy}, $\k$ denotes the permeability of the porous medium. For simplicity, 
it is assumed to be constant and positive.
We refer to~\cite{CGM17} for the case of space-dependent anisotropic permeability tensors.
The fluid viscosity and density are denoted by $\mu_i>0$ and $\rho_i\geq0$, respectively, whereas $\g= - g \boldsymbol{e}_z$ 
denotes the gravity.
The unknown phase pressures $\bp = \left(p_i\right)_{0 \leq i \leq N}$ are
related to the saturations by $N$ capillary pressure relations 
\be\label{eq:capi}
p_i - p_0 = \pi_i(\bs^\ast), \qquad \forall i \in \{1,\dots, N\}. 
\ee
\end{subequations}
The capillary pressure functions $\bpi = \left(\pi_i\right)_{1\leq i \leq N}$ are assumed to derive 
from a strictly convex and $\varpi$-concave potential $\Pi: \bDelta^* \to \R_+$ for some $\varpi>0$, i.e., 
\be\label{eq:H-Pi}
0 <  \Pi(\wh \bs^*) - \Pi(\bs^*) - \bpi(\bs^*) \cdot (\wh \bs^* - \bs^*)
\leq \frac{\varpi}2|\wh \bs^* - \bs^*|^2, \qquad \forall \bs^*, \wh \bs^* \in \bDelta^*, \;\text{with}\; \bs^* \neq \wh \bs^*.
\ee
This implies that $\bpi:\bDelta^* \to \R^N$ is strictly monotone (thus one-to-one) and Lipschitz continuous:
\[
0 < \left( \bpi(\wh \bs^*) - \bpi(\bs^*)\right) \cdot\left(\wh \bs^* - \bs^* \right)\leq \varpi |\wh \bs^* - \bs^*|^2, 
\qquad \forall \bs^*, \wh \bs^* \in \bDelta^*, \;\text{with}\; \bs^* \neq \wh \bs^*,
\]
and thus
\[
0 \leq D^2 \Pi(\bs^*) \leq \varpi {\bf I}_N, \qquad \forall \bs^* \in \bDelta^*.
\]
The last inequalities have to be understood in the sense of the symmetric matrices. The function $\Pi$ is extended 
by $+\infty$ outside of $\bDelta^*$.

As established in~\cite{CGM17}, the problem~\eqref{eq:EL} can be interpreted as the 
Wasserstein gradient flow of the energy 
\be\label{eq:Ee}
\Ee(\bs) = \int_\O \left[\Pi(\bs^*) + \bs\cdot \bPsi + \chi_\bDelta(\bs)\right]\d\x, \qquad \forall \bs \in \bAa.
\ee
In \clem{f}ormula~\eqref{eq:Ee}, the exterior gravitational potential $\bPsi = \left( \Psi_i \right)_{0\leq i \leq N}$ 
is given by 
\be\label{eq:Psi_i}
\Psi_i(\x) = - \rho_i \g\cdot \x,  \qquad \forall \x \in \O.
\ee
\begin{rem}
In fact in \eqref{eq:Ee} one can consider a large class of arbitrary potential $\bPsi$, see \cite{CGM17} for details. 
\end{rem}
The constraint~\eqref{eq:cons} is incorporated in the energy 
rather than in the geometry
thanks to the term
\[
\chi_\bDelta(\bs) = \begin{cases}
0 & \text{if}\; \bs \in \bDelta, \\
+\infty & \text{otherwise}.
\end{cases}
\]
We refer to~\cite{BP02, BBG04, Loeschcke_PhD} for a presentation of the multiphase optimal transportation problem
for which the constraint~\eqref{eq:cons} is directly incorporated in the geometry.
In order to be more precise in our statements, we need to introduce some extra material 
concerning the Wasserstein distance to be used to equip $\bAa$. This is the purpose of the 
next section.

\begin{rem}
In the previous work~\cite{CGM17}, the uniform convexity of the capillary potential $\Pi$ was required. 
In~\eqref{eq:H-Pi}, we relax this assumption into a mere strict convexity requirement. This can be done 
by slightly adapting the proofs 
of~\cite{CGM17}. 
\end{rem}
\subsection{Wasserstein distance}\label{ssec:MK2}

For $i \in \{0,\dots, N\}$ we define 
$$
\Aa_i = \left\{ s_i \in L^1(\O; \R_+) \, \middle| \, \int_\O s_i \d\x = m_i \right\}.
$$
Given $s_{i}, \wh s_{i}\in\Aa_i$, the set of admissible transport plans between $s_i$ and $\wh s_i$ 
is given by
$$
\Gamma_i(s_{i}, \wh s_{i}) = \left\{ \gamma_i \in \Mm_+(\O\times\O)\, \middle| 
\,\gamma_i(\Omega\times\Omega)=m_i,\, 
\gamma_i^{(1)} = s_i \; \text{and} \; \gamma_i^{(2)} = \wh s_i \;
\right\},
$$
where $\Mm_+(\O\times\O)$ stands for the set of Borel measures on $\O\times \O$ 
and $\gamma_i^{(k)}$ is the $k^\text{th}$ marginal of the measure $\gamma_i$.
The quadratic Wasserstein distance $W_i$ on $\Aa_i$ is then defined as
\be
\label{eq:Wi}
W_i^2(s_{i}, \wh s_{i}) = 
\min\limits_{\gamma_i \in \Gamma(s_{i}, \wh s_{i})} \iint_{\O\times\O}
\frac{ \mu_i }{\kappa} |\x-\by|^2\d \gamma_i(\x,\by).
\ee
Equivalently, the continuity equation \eqref{eq:conv} allows to give the following dynamical characterization:
\begin{prop}[Benamou-Brenier formula \cite{BB00}]
For $s_{i,0},s_{i,1}\in \Aa_i$ we have
\be
\label{eq:Wi_BB}
W_i^2(s_{i,0}, s_{i,1})=\min\limits_{s_i,\bv}\int_0^1\int_\O \frac{ \mu_i }{\kappa}|\bv_t(\x)|^2\d s_{i,t}(\x) \d t,
\ee
where the minimum runs over curves of measures $t\mapsto s_{i,t}\in\Aa_i$ with endpoints $s_{i,0}, s_{i,1}$ and velocity fields $t\mapsto\bv_t\in\Mm^d(\O)$ such that
$$
\partial_t s_{i,t}+\div(s_{i,t}\bv_t)=0
$$
in the sense of distributions.
\end{prop}
\begin{rem}
\label{rem:momentum}
As originally developped in \cite{BB00}, the right variables to be used in the Benamou-Brenier formula \eqref{eq:Wi_BB} is not the velocity $\bv$, but in fact the momentum $\bm=s\bv$, since the action $A(s,\bm)=\frac{|\bm|^2}{s}=s|\bv|^2$ becomes then jointly convex in both arguments.
\end{rem}

A third equivalent formulation is the Kantorovich dual problem:
\begin{prop}
 There holds
 \be
 \label{eq:Wi_Kanto}
\frac{1}{2} W_i^2(s_{i}, \wh s_{i})=\max\limits_{\phi,\psi} \left\{
 \int\phi(\x)\d s_i(\x) +\int\psi(\by)\d \wh s_i(\by)
 \right\},
 \ee
 where the maximum runs over all pairs $(\phi,\psi)\in L^1(\d s_i)\times L^1(\d\wh s_i)$ such that $\phi(\x)+\psi(\by)\leq \frac{ \mu_i }{2\kappa} |\x-\by|^2 $.
 Any maximizer is called a (pair of optimal) Kantorovich potential.
 
\end{prop}

The viscosity $\mu_i$ and permeability $\kappa$ appear in \eqref{eq:Wi}--\eqref{eq:Wi_Kanto} as scaling factors in the cost function $\mu_i|\x-\by|^2/\kappa$, and this is required for consistency with Darcy's law \eqref{eq:Darcy}.
For more general heterogeneous permeability tensors $\mathbb K(\x)$ one could use instead the intrinsic distance $d^2_i(\x,\by)$ induced on $\Omega$ by the Riemannian tensor $\mu_i\bbK^{-1}(\x)$, see \cite{Lisini09} for a general approach of Wasserstein distances with variable coefficients and \cite{CGM17} in the particular context of multiphase flows in porous media.

\smallskip

With the phase Wasserstein distances $\left(W_i\right)_{0\le i \le N}$ at hand, 
we can define the global Wasserstein distance $\bW$ on 
$\bAa:=\Aa_0 \times \dots \times \Aa_N$ by setting 
$$
\bW(\bs, \wh \bs) = \left( \sum_{i=0}^N W_i(s_i, \wh s_i)^2\right)^{1/2}, \quad \forall \bs, \wh \bs \in \boldsymbol{\mathcal A}.
$$

\subsection{Approximation by minimization scheme}\label{ssec:JKO}

As already mentioned, the problem \eqref{eq:EL} 
is the Wasserstein gradient flow of our singular energy \eqref{eq:Ee}, see our earlier works \cite{CGM15,CGM17}.
Rather than discussing the meaning of gradient flows in the Wasserstein setting, we refer to the monograph \cite{AGS08} for an exposition of gradient flows in abstract metric spaces \cite{AGS08} and to \cite{Santambrogio_OTAM,Santambrogio17} for a detailed overview.
As is now well understood from the work of Jordan, Kinderlehrer, and Otto \cite{JKO98}, one possible way to formalize this gradient flow structure is to implement the JKO scheme (also referred to as DeGiorgi's minimizing movement, see~\cite{DeGiorgi93}).
Given an initial datum $\bs^0\in \boldsymbol{\mathcal A}$ with energy $\Ee(\bs^0)<\infty$ and a time step $\tau>0$, 
the strategy consists in:
\begin{enumerate}[(i)]
 \item 
construct a time discretization $\bs^n(.)\approx\bs(n\tau,.)$ by solving recursively
\begin{equation}
 \label{JKOscheme}
 \bs^{n+1} = \underset{\bs\in\boldsymbol{\mathcal A}}{\argmin}
 \left\{
 \frac{1}{2\tau}\bW^2(\bs,\bs^n) + \Ee(\bs)
 \right\};
\end{equation}
\item
define the piecewise-constant interpolation
$$
\bs_\tau(t):= \bs^{n+1}
\qquad
\mbox{if } t \in (n\tau,(n+1)\tau];
$$
\item
retrieve a continuous solution $\bs(t)=\lim\limits_{\tau\to 0}\bs_\tau(t)$ in the limit of small time steps.
\end{enumerate}
This is a variant in the Wasserstein space of the implicit variational Euler scheme: indeed, in Euclidean spaces $x\in \R^d$ and for smooth functions $E:\R^d\to \R$, the Euler-Lagrange equation corresponding to minimizing $x \mapsto \frac{1}{2\tau}|x-x^n|^2+E(x)$ is nothing but the finite difference approximation $\frac{x^{n+1}-x^n}{\tau}=-\nabla E(x^{n+1})$.
We refrain from giving more details at this stage and refer again to \cite{AGS08,Santambrogio_OTAM,Villani09}.

Due to lower semi-continuity and convexity, it is easy to prove that the minimization problem \eqref{JKOscheme} 
is well-posed, hence the discrete solution $\bs_{\tau}$ is uniquely and unambiguously defined.
But we still need to construct approximate phase pressures $\bp_\tau = (p_{1,\tau},p_{2,\tau})$. 
Their construction makes use of the backward Kantorovich potentials (see~\cite[Section 3]{CGM17}).

\begin{lem}
\label{lem:discrete_E-L_JKO}
There exist pressures $p_i^{n+1}$ and Kantorovich potentials $\phi_i^{n+1}$ (from $s_i^{n+1}$ to $s_i^{n}$) such that
\begin{equation}
\label{eq:discrete_E-L_JKO}
\frac{\phi_i^{n+1}}{\tau}=p_i^{n+1} + \Psi_i\quad \text{a.e. in}\; \{s_i^{n+1}>0\}\quad\mbox{for }i=0,\dots ,N,
\end{equation}
and
\begin{equation}
\label{eq:discrete_E-L_JKO_capi}
p_i^{n+1}-p_0^{n+1}=\pi_i(\bs^{{n+1},*})\quad  \text{a.e. in}\; \O\quad\mbox{for }i=1,\dots ,N.
\end{equation}

\end{lem}
From classical optimal transport theory \cite{Santambrogio_OTAM}, $\bv_i^{n+1}:=\frac{\kappa}{\mu_i}\frac{\nabla\phi_i^{n+1}}{\tau}$ should be interpreted as the discrete velocity driving the i-th phase, which will automatically give $\partial_t s_i+\div(s_i\bv_i)=0$ in the limit $\tau\to 0$. Hence \eqref{eq:discrete_E-L_JKO} is a discrete counterpart of Darcy law \eqref{eq:Darcy}. 
The capillary relation \eqref{eq:capi} hold as well at the discrete level thanks to relations~\eqref{eq:discrete_E-L_JKO_capi}, 
whereas the total saturation constraint \eqref{eq:cons} is automatically enforced in \eqref{JKOscheme} thanks to $\Ee(\bs^{n+1})<\infty$.
For the sake of brevity we omit the details and refer again to \cite{CGM17}.

\subsection{Main properties of the approximation}\label{ssec:Kapla}

Since our system \eqref{eq:EL} of PDEs is highly nonlinear, taking the limit $\bs(t)=\lim\limits_{\tau\to 0}\bs_\tau(t)$ will require sufficient compactness both in time and space.
In this section we sketch the main arguments leading to such compactness.

Compactness in time is derived from the classical \emph{total square distance estimate} below, which is a characteristic feature of any JKO variational discretization.
Testing $\bs=\bs^n$ as a competitor in \eqref{JKOscheme} gives first 
\[\frac{1}{2\tau}\bW^2(\bs^{n+1},\bs^n) + \Ee(\bs^{n+1}) \leq \Ee(\bs^{n}).\]
This implies of course the energy monotonicity $\Ee(\bs^{n+1}) \leq \Ee(\bs^{n})$, but summing over $n$, 
we also get the \emph{total square distance estimate} in the form
\begin{equation}
 \label{eq:tot_sq_dist}
 \frac{1}{\tau} \sum \limits_{n\geq 0} \bW^2(\bs^{n+1},\bs^n)
 \leq 2 \left( \Ee(\bs^0)-\inf\limits_{\Aa} \Ee\right).
\end{equation}
By definition of the piecewise-constant interpolation, an easy application of the Cauchy-Schwarz inequality gives then the \emph{approximate equicontinuity}
$$
\bW(\bs_\tau(t_1),\bs_\tau(t_2))\leq C|t_2-t_1+\tau|^{\frac 12},
\qquad \forall\,0\leq t_1\leq t_2,
$$
uniformly in $\tau$, which yields the desired compactness in time (see~\cite[Proposition 3.10]{AGS08} or
\cite[Theorem C.10]{Kangourou}).\\

Compactness in space will be obtained exploiting the \emph{flow interchange technique} from \cite{MMS09}.
Roughly speaking, this amounts to estimating the dissipation of the driving functional $\Ee$ along a well-behaved auxiliary gradient flow, driven by an auxiliary functional and starting from the minimizer $\bs^{n+1}$.
More explicitly, we define the $\eps$-perturbation $\tilde\bs_\eps=(\tilde s_{0,\eps},\dots,\tilde s_{N,\eps})$ as solutions to the independent heat equations
$$
\left\{
\begin{array}{ll}
\frac{\partial\tilde s_{i,\eps}}{\partial\eps}=\kappa\Delta \tilde s_{i,\eps} & \mbox{for small }\eps>0, \\
\left.\tilde s_i\right|_{\eps=0} = s^{n+1}_{i}.
\end{array}
\right. 
$$
The key observation is that, for each $i=0\dots N$, the above heat equation is a gradient flow in the Wasserstein space $(\Aa_i,W_i)$ with driving functional $\mu_i\mathcal H$, where the Boltzmann entropy
\begin{equation}\label{eq:H_omega}
\mathcal H(s) = \int_{\O} s(\x) \log\left(s(\x)\right) \d\x.
\end{equation}
In addition to the usual regularizing effects, this heat equation is particularly well-behaved here in the sense that it preserves the total saturation constraint $\sum\limits_{i=0}^N \tilde s_{i,\eps}=\sum\limits_{i=0}^N s_i^{n+1} =1$ and, since $\O$ is convex, the auxiliary driving functional $\mathcal H$ is displacement convex in $(\Aa_i,W_i)$ \cite{Villani09,McCann97}. 
If
$$
\mathcal F_\tau^n(\bs)=\frac{1}{2\tau}\bW^2(\bs,\bs^{n})+\Ee(\bs)
$$
denotes the JKO functional, then by optimality of the minimizer $\bs^{n+1}$ in \eqref{JKOscheme} we must have
$$
\limsup\limits_{\eps\to 0^+} \frac{d}{d\eps}\mathcal F_\tau^n(\tilde\bs_\eps)\geq 0.
$$
The energy term $\mathcal E(\tilde\bs_\eps)=\int_\Omega \Pi(\tilde \bs_\eps^*)$ can easily be differentiated under the integral sign (with respect to $\eps$), while the variation of the first $\bW^2(\tilde \bs_\eps,\bs^n)$ term can be estimated using the \emph{evolution variational inequality} \cite{AGS08} for the well-behaved $\Hh$-flow $\tilde \bs_\eps$ (this metric characterization precisely requires some displacement convexity of the auxiliary flow, see \cite[Theorem 2.23]{EKS15}).
Omitting again the details, one gets in the end the dissipation estimate
$$
\tau\sum\limits_{i=0}^N \|\nabla \pi_i((\bs^{n+1})^*)\|_{L^2(\O)}\leq C
\left(
\tau + \bW^2(\bs^{n+1},\bs^n) + \sum\limits_{i=0}^N \mathcal H(s^{n}_{i})-\mathcal H(s^{n+1}_{i})
\right),
$$
see \cite[Section 2.2]{CGM17} for the details.
Exploiting the previous total square distance estimate and summing over $n=0\dots\lfloor T/\tau\rfloor$ (or equivalently integrating in time), we control next
\begin{equation}
 \label{eq:pi_L2H1}
\|\bpi(\bs_\tau^*)\|_{L^2(0,T;H^1)}\leq C \left(T+\sum\limits_{i=0}^N\mathcal H(s^0_i) +1\right)=C_T
\end{equation}
for arbitrary $T>0$ and fixed initial datum $\bs^0$.
 It is worth recalling at this stage that, due to our assumption \eqref{eq:H-Pi}, $\bpi(\bs^*)=\nabla_{\bs^*}\Pi(\bs^*)$ 
 is a strictly monotone thus invertible map of $\bs^*$ due to the strict convexity of $\Pi$. 
 The compactness w.r.t. the space variable of $\left(\bs_\tau\right)_{\tau>0}$ then follows from 
\eqref{eq:pi_L2H1}.
 \begin{rem}
 \label{rem:log_test_function}
  A formal but more PDE-oriented explanation of the above flow-interchange simply consists in taking $\log(s_i)$ as a test function in the weak formulation of system \eqref{eq:EL}.
  The delicate technical part is to justify this computation and mimic this formal chain rule in the discrete time setting in order to retrieve enhanced regularity of the JKO minimizers.
   \end{rem}

 Exploiting the above compactness, one can argue as in \cite{CGM17} and finally prove the following convergence results. 
 The existence of a weak solution to the problem~\eqref{eq:EL} is a direct byproduct. 
 \begin{thm}
 For any discrete sequence $\tau_k\to 0$ and up extraction of a subsequence if needed, we have convergence
 $$
 \begin{array}{cl}
  \bs_{\tau_k}\to \bs &\mbox{strongly in all }L^q((0,T)\times\O),\\
  \bpi(\bs^*_{\tau_k}) \rightharpoonup \bpi(\bs^*) &\mbox{weakly in }L^2(0,T;H^1(\O)),\\
  \bp_\tau \rightharpoonup \bp & \mbox{weakly in }L^2(0,T;H^1(\O)),
 \end{array}
 $$
 and the limit $(\bs,\bp)$ is a weak solution of \eqref{eq:EL}.
 \end{thm}

\section{Numerical approximation of the flow}\label{sec:scheme}

We present here the ALG2-JKO scheme and the \clem{u}pstream mobility \clem{f}inite \clem{v}olume scheme.
The first method is based on the variational JKO scheme \eqref{JKOscheme} described in subsection \ref{ssec:JKO} whereas the second method is based on the PDE formulation of
the problem~\eqref{eq:EL} given by \eqref{eq:cons}-\eqref{eq:conv}-\eqref{eq:Darcy}-\eqref{eq:capi}. 
Both methods are well adapted for gradient flows equations, and more precisely we will check the following key properties for the numerical solutions:  
\begin{itemize} 
\item preservation of the positivity
\item conservation of the mass and saturation constraints,
\item energy dissipation along solutions.
\end{itemize}

\subsection{The ALG2-JKO scheme}\label{ssec:ALG2}

This algorithm relies on the seminal work of Benamou and Brenier~\cite{BB00} where an 
augmented Lagrangian approach was used to compute Wasserstein distances.
In~\cite{BCL15}, this approach was extended to the computation of 
Wasserstein gradient flows.
The method is very well suited for computing solutions to constrained gradient flows, as it will appear in the numerical 
simulations presented in Section~\ref{sec:numerics}. 

\subsubsection{The augmented Lagrangian formulation}
Roughly speaking, the ALG2-JKO scheme consists in rewriting the single JKO step~\eqref{JKOscheme} as a more fashionable (and effectively implementable) convex minimization problem.
In order to do so, let us first introduce the convex lower-semicontinuous $1$-homogeneous action function given, for all $(s,\boldsymbol{m}) \in \R \times {\R}^d$, by 
\begin{equation}
\label{eq:def_action_BB}
A( s,\boldsymbol{m}):= 
\left\{ \begin{array}{ll}
\frac{|\boldsymbol{m}|^2}{2s} & \text{ if } s>0,\\
0 & \text{ if } s=0 \text{ and } \boldsymbol{m}=\boldsymbol{0},\\
+\infty & \text{ otherwise.}
\end{array} \right. 
\end{equation}
We recall that $\bm=s\bv$ is the momentum variable in the continuity equation {$\partial_t s+\div(s\bv)=0$} and $|\bm|^2/s=s|\bv|^2$ is a kinetic energy, see Remark~\ref{rem:momentum}.
As originally observed in \cite{BB00}, the function $A$ can be seen as the support function
\begin{eqnarray}
\label{eq: dual Phi}
A(s,\boldsymbol{m})= 
\sup_{(a,\boldsymbol{b})\in K_2} \left\{ a s +\boldsymbol{b}\cdot \boldsymbol{m}\right\}
\end{eqnarray}
of the convex set $K_2$, where $K_\alpha$ is defined for $\alpha>0$ as
\begin{equation}
\label{eq:def_K_alpha}
 K_\alpha:=\left\{ (a,\boldsymbol{b}) \in \R \times {\R}^d \, : \, a +\frac{1}{\alpha} |\boldsymbol{b}|^2 \leq 0 \right\}.
 \end{equation}
Taking advantage of the Benamou-Brenier formula \eqref{eq:Wi_BB}, and given the previous JKO step $\bs^n$, \eqref{JKOscheme} can be recast as
\begin{eqnarray}
\label{eq:JKO BB}
\min_{\bs,\mathbf{m}}
\left\{
\sum_{i=0}^N \frac{\mu_i}{\kappa}\int_0^1 \int_\Omega A(s_{i,t}(\x),\boldsymbol{m}_{i,t}(\x)) \, \d\x \d t + \tau \Ee(\bs_{t=1})
\right\},
\end{eqnarray}
where the infimum runs over curves of measures $t\mapsto \bs_t=(s_{0,t},\dots,s_{N,t})\in\bAa$ and momenta $t\mapsto \bm_t=(\bm_{0,t},\dots,\bm_{N,t})\in\Mm^d(\O)^{N+1}$, subject to $N+1$ linear constraints 
\begin{equation}
\label{eq:linear constraint JKO BB}
\left\{
\begin{array}{ll}
\partial_t s_{i,t} + \div(\boldsymbol{m}_{i,t}) =0 & \mbox{in }\mathcal D',\\
\boldsymbol{m}_{i,t} \cdot \nu =0 & \mbox{on }\partial \O,\\
\left.{s_{i}}\right|_{t=0}= s_{i}^n,
\end{array}
\right.
\hspace{1cm} i=0,\dots, N.
\end{equation}
Note that only the initial endpoint $\bs_{t=0}=\bs^{n}$ is prescribed for the curve $(\bs_t)_{t\in[0,1]}$.
The terminal endpoint is free and contributes to the objective functional \eqref{eq:JKO BB} through the $\Ee(\bs_{t=1})$ term, and the JKO minimizer will be retrieved as $\bs^{n+1}=\bs_{t=1}$.
Note also that the minimizing curve $(\bs_t)_{t\in[0,1]}$ in \eqref{eq:JKO BB}--\eqref{eq:linear constraint JKO BB} will automatically be a Wasserstein geodesic between the successive JKO minimizers $\bs_{t=0}=\bs^n$ and $\bs_{t=1}=\bs^{n+1}$.

As a first step towards a Lagrangian formulation, we rewrite the constraint \eqref{eq:linear constraint JKO BB} as a sup problem with multipliers $\phi_i(t,\x)$
\begin{multline}
\label{eq:sup continuity equations}
\sup_{\bphi} \left\{ \sum_{i=0}^N \int_\Omega \phi_i(1,\cdot) s_{i,1} - \int_\Omega \phi_i(0,\cdot) s_{i}^n - \int_0^1 \int_\Omega (\partial_t \phi_i s_{i,t} + \grad \phi_i\cdot  \boldsymbol{m}_{i,t} ) \right\}
\\
=\left\{
\begin{array}{ll}
 0 & \mbox{if } \eqref{eq:linear constraint JKO BB} \mbox{ holds},\\
 +\infty & \mbox{else},
\end{array}
\right.
\end{multline}
and minimizing \eqref{eq:JKO BB} under the constraint \eqref{eq:linear constraint JKO BB} can thus be written $\inf\limits_{\bs,\bm}\sup\limits_{\boldsymbol\phi}\,\{\dots\}$.
Swapping $\inf\sup=\sup\inf$ as in \cite{BCL15} and using that the Legendre transform of $\frac{\mu_i}{\kappa}A$ is the charateristic function (convex indicator) of the convex set $K_{2\mu_i/\kappa}$ defined in \eqref{eq:def_K_alpha},
$$
\left(\frac{\mu_i}{\kappa}A\right)^*(a,\boldsymbol b)=\chi_{K_{2\mu_i/\kappa}}(a,\boldsymbol b)=
\left\{
\begin{array}{ll}
 0 & \mbox{if }(a,\boldsymbol b)\in K_{2\mu_i/\kappa},\\
 +\infty & \mbox{else,}
\end{array}
\right.
$$
the problem \eqref{eq:JKO BB}-\eqref{eq:linear constraint JKO BB} finally becomes after a few elementary manipulations
$$
 \inf_{\boldsymbol{\phi}} \left\{ \sum_{i=0}^N \int_\Omega \phi_{i}(0,\cdot) s_{i}^n  +  \Ee^*_\tau (-\boldsymbol{\phi}(1,\cdot)) \, : \qquad  (\partial_t \phi_i, \grad \phi_i) \in K_{2\mu_i/\kappa} \right\}.
 $$
Here $\Ee_\tau^*$ denotes the Legendre transform of $\Ee_\tau := \tau \Ee$. 
This dual problem can be reformulated as
$$
 \inf_{\boldsymbol{\phi}}\Big\{ F(\boldsymbol{\phi}) + G(\boldsymbol{q}) \, : \qquad \boldsymbol{q}=\Lambda \boldsymbol{\phi} \Big\},
 $$
where
$$
\Lambda \boldsymbol{\phi} = (\partial_t \boldsymbol{\phi} , \grad \boldsymbol{\phi} , - \boldsymbol{\phi}(1,\cdot) ) 
\qquad\mbox{ and }\qquad
\boldsymbol{q}=(\boldsymbol{a},\boldsymbol{b},\boldsymbol{c})
$$
are functions with values in $(\R \times {\R}^d \times \R)^{N+1}$,
$$
 F(\boldsymbol{\phi}) = \sum_{i=0}^N \int_\Omega \phi_i(0, \cdot) s_{i}^n,
 $$
$$
 G(\boldsymbol{q}) = \sum_{i=0}^N \int_0^1 \int_\Omega \chi_{K_{2\mu_i/\kappa}} (a_i,\boldsymbol{b}_i) + \Ee^*_\tau (\boldsymbol{c}),
 $$
 and $\chi_{K_{2\mu_i/\kappa}}$ stands again for the characteristic function of $K_{2\mu_i/\kappa}$.
 Introducing a Lagrange multiplier  
 $$
 \boldsymbol{\sigma} = ( \bs, \boldsymbol{m}, \tilde\bs_1)
 $$
for the constraint $\Lambda \boldsymbol{\phi}=\boldsymbol{q}$, finding a minimizer $\bs^{n+1}$ in the JKO scheme \eqref{JKOscheme} is thus equivalent to finding a saddle-point of the Lagrangian
\begin{eqnarray}
\label{lagrangian}
L(\boldsymbol{\phi}, \boldsymbol{q}, \boldsymbol{\sigma}) := F(\boldsymbol{\phi}) + G(\boldsymbol{q}) + \boldsymbol{\sigma} \cdot (\Lambda \boldsymbol{\phi} - \boldsymbol{q} ).
\end{eqnarray}
 Here we slightly abuse the notations: $\bs=(\bs_t)_{t\in[0,1]}$ and $\bm=(\bm_t)_{t\in[0,1]}$ are time-depending curves while $\tilde \bs_1\in \bAa$ is independent of time.
 The scalar product in \eqref{lagrangian} is
$$
  \boldsymbol{\sigma} \cdot (\Lambda \boldsymbol{\phi} - \boldsymbol{q} ) = \sum_{i=0}^N \Big( \int_0^1 \int_\Omega \left( s_{i} ( \partial_t \phi_i -a_i) +\boldsymbol{m}_{i} \cdot (\grad \phi_i -\boldsymbol{b}_i) \right) -\int_\Omega \tilde{s}_{1,i}( \phi_i(1,\cdot) +c_i) \Big).
  $$
We stress that the free variable $\tilde\bs_1$ is a priori independent of the curve $(\bs_t)_{t\in[0,1]}$, but that the saddle-point will ultimately satisfy $\bs_{t=1}=\tilde\bs_1$.
In the Lagrangian \eqref{lagrangian}, the original unknowns $(\bs, \boldsymbol{m}, \tilde\bs_1)$ become the Lagrange multipliers for the constraint $\boldsymbol{q}=\Lambda \boldsymbol{\phi}$, \clem{i.e.,} respectively
$$
  \boldsymbol{a}= \partial_t \boldsymbol{\phi},
 \qquad
 \boldsymbol{b}=\grad \boldsymbol{\phi},
 \qquad\text{and}
 \qquad
 \boldsymbol{c}=-\boldsymbol{\phi}(1,\cdot).
 $$ 

For some fixed regularization parameter $r>0$, we introduce now the augmented Lagrangian
 \begin{eqnarray}
\label{aug lagrangian}
L_r(\boldsymbol{\phi}, \boldsymbol{q}, \boldsymbol{\sigma}) := F(\boldsymbol{\phi}) + G(\boldsymbol{q}) + \boldsymbol{\sigma} \cdot (\Lambda \boldsymbol{\phi} - \boldsymbol{q} ) +\frac{r}{2}\|\Lambda \boldsymbol{\phi} - \boldsymbol{q}\|^2,
\end{eqnarray}
where the extra regularizing term is given by the $L^2$ norm
$$
 \frac{r}{2}\|\Lambda \boldsymbol{\phi} - \boldsymbol{q}\|^2 = \frac{r}{2}\sum_{i=0}^N\left( \int_0^1 \int_\Omega ( |\partial_t \phi_i -a_i |^2 + |\grad \phi_i -\boldsymbol{b}_i|^2) +\int_\Omega |\phi_i(1,\cdot) +c_i |^2 \right).
 $$
Observe that being a saddle-point of \eqref{lagrangian} is equivalent to being a saddle-point of \eqref{aug lagrangian}, see for instance \cite{FG83}.
Thus in order to solve one step of the JKO scheme~\eqref{JKOscheme}, it suffices to find a saddle-point of the augmented Lagrangian $L_r$.

\subsubsection{Algorithm and discretization} \label{eq:subsec_algo_ALG2JKO}

The augmented Lagrangian algorithm ALG2 aims at finding a saddle-point of $L_r$ and consists in a splitting scheme.
Starting from $(\boldsymbol{\phi}^0,\boldsymbol{q}^0,\boldsymbol{\sigma}^0)$, we generate a sequence $(\boldsymbol{\phi}^k,\boldsymbol{q}^k,\boldsymbol{\sigma}^k)_{k\geq 0}$ by induction as follows

\begin{itemize}
\item[\textbf{Step 1:}] minimize with respect to $\boldsymbol{\phi}$:
$$ 
\boldsymbol{\phi}^{k+1}= \argmin_{\boldsymbol{\phi}} \left( F(\boldsymbol{\phi}) + \boldsymbol{\sigma}^k \cdot \Lambda\boldsymbol{\phi} +\frac{r}{2}|\Lambda \boldsymbol{\phi} - \boldsymbol{q}^{k}|^2 \right),
$$

\item[\textbf{Step 2:}] minimize with respect to $\boldsymbol{q}$:
$$ 
\boldsymbol{q}^{k+1}= \argmin_{\boldsymbol{q}} \left( G(\boldsymbol{q}) - \boldsymbol{\sigma}^k \cdot \boldsymbol{q} +\frac{r}{2}|\Lambda \boldsymbol{\phi}^{k+1} - \boldsymbol{q}|^2 \right),
$$

\item[\textbf{Step 3:}] maximize with respect to $\boldsymbol \sigma$, which amounts here to updating the multiplier by the gradient ascent formula
$$ 
\boldsymbol{\sigma}^{k+1}= \boldsymbol{\sigma}^k +r (\Lambda \boldsymbol{\phi}^{k+1} - \boldsymbol{q}^{k+1}).
$$
\end{itemize}
Since step 3 is a mere pointwise update we only describe in details the first two steps.
In order to keep the notations light we sometimes write $s_i(t,\x)=s_{i,t}(\x)$, and likewise for any other variable depending on time.
\begin{itemize}
 \item 
The first step corresponds to solving $N+1$ independent linear elliptic problems in time and space, namely
$$
 -r \boldsymbol{\Delta}_{t,\x} \phi_i^{k+1} = \grad_{t,\x} \cdot( (s_{i}^k, \boldsymbol{m}_{i}^k) -r(a_i^k , \boldsymbol{b}_i^k) )
 \hspace{1cm}
 \text{ in } (0,1) \times \Omega
 $$
with the boundary conditions
$$
\left\{
\begin{array}{ll}
r\partial_t \phi_i^{k+1}(0, \cdot) = s_{i}^n(\cdot) - s_{i}^k(0,\cdot) + ra_i^k(0,\cdot) & \mbox{ in }\O,\\
r \, \Big( \partial_t \phi_i^{k+1}(1,\cdot) + \phi_i^{k+1}(1,\cdot) \Big )
= \tilde s_{1,i}^k(\cdot) - s_{i}^k(1,\cdot) + r\, \Big (a_i^k(1,\cdot) -c_i^k(\cdot)\Big )& \mbox{ in }\O,\\
\left(r\grad \phi_i^{k+1} + \boldsymbol{m}_{i}^k -r\boldsymbol{b}_i^k \right) \cdot \nu =0 & \text{ on } \partial \Omega.
\end{array}
\right.
$$
\item
The second step splits into two convex pointwise subproblems.
The first one corresponds to projections onto the parabolas $K_{2\mu_i/\kappa}$:
$$
(a_i^{k+1} , \boldsymbol{b}_i^{k+1})(t,\x)= P_{K_{2\mu_i/\kappa}}\left( (\partial_t \phi_i^{k+1}, \grad \phi_i^{k+1})(t,\x) +\frac{1}{r}(s_{i}^k, \boldsymbol{m}_{i}^k)(t,\x) \right),
\qquad \forall \,i=0,\dots , N.
$$
This projection $P_{K_{2\mu_i/\kappa}}$ onto $K_{2\mu_i/\kappa}$ is explicitly given by (see \cite{PPO14})
$$ 
P_{K_{2\mu_i/\kappa}}(\alpha, \boldsymbol{\beta})= \left\{ \begin{array}{ll}
(\alpha, \boldsymbol{\beta}), & \text{ if } (\alpha, \boldsymbol{\beta}) \in K_{2\mu_i/\kappa},\\
\left(\alpha -\lambda , \frac{\mu_i \boldsymbol{\beta}}{\kappa\lambda + \mu_i}\right), & \text{ otherwise,}
\end{array}\right.
$$
where $\lambda$ is the largest real root of the cubic equation 
$$ 
(\alpha -\lambda)(\mu_i/\kappa +\lambda )^2  + \frac{\mu_i}{2\kappa}|\boldsymbol{\beta} |^2 =0.
$$
The second subproblem should update $\boldsymbol{c}$.
To this end, we need to solve the pointwise proximal problem: for each $\x \in \Omega$
\begin{eqnarray}
\label{JKO step2 second subproblem}
\boldsymbol{c}^{k+1}(\x) = \argmin_{\boldsymbol{c} \in {\R}^{N+1}}
\left\{\frac{r}{2}\sum_{i=0}^N |\phi_i^{k+1}(1,\x) - \frac{1}{r}\tilde s_{1,i}^{k}(\x) + c_i |^2 + E_\tau^*(\x,\boldsymbol{c})\right\},
\end{eqnarray} 
where $E^*_\tau(\x,\cdot)$ is the Legendre transform of the energy density $E_{\tau}(\x,\cdot)=\tau E(\x,\cdot)$ in its second argument ($E$ being implicitly defined as $\Ee(\boldsymbol{s})= \int_\Omega E(\x,\boldsymbol{s}(\x)) \, \d\x$).
\end{itemize}
Notice that the energy functional $\mathcal E$ only plays a role in the minimization with respect to the internal $\boldsymbol{c}$ variable, namely the second subproblem \eqref{JKO step2 second subproblem} in {\bf Step 2}.
In Section~\ref{sec:numerics} we will try to make this step explicit for our two particular applications.\\

In order to implement this algorithm in a computational setting we use P2 finite elements in time and space for $\boldsymbol{\phi}$, and P1 finite elements for $\boldsymbol{\sigma}$ and $\boldsymbol{q}$.
The variables $\grad_{t,\x} \phi_i^{k+1}=(\partial_t  \phi_i^{k+1}, \grad \phi_i^{k+1})$ are understood as the projection onto P1 finite elements and the algorithm was implemented using \texttt{FreeFem++}~\cite{FreeFem}.
The convergence of this algorithm is known in finite dimension \cite{FG83}, i.e., the iterates $(\boldsymbol{\phi}^k,\boldsymbol{q}^k,\boldsymbol{\sigma}^k)$ are guaranteed to converge to a saddle point $(\boldsymbol{\phi},\boldsymbol{q},\boldsymbol{\sigma})$ as $k\to\infty$.
Once the saddle-point is reached, the output $\boldsymbol{\sigma}=(\bs, \boldsymbol{m}, \tilde\bs_1)$ is a minimizer for the problem \eqref{eq:JKO BB}-\eqref{eq:linear constraint JKO BB} and the solution of the JKO scheme \eqref{JKOscheme} is simply recovered as $\bs^{n+1}=\tilde\bs_1=\bs|_{t=1}$.

Numerically, the Benamou-Brenier formula involves an additional time dimension to be effectively discretized in each elementary JKO step, and this can be seen as a drawback.
However the successive JKO densities are close due to the small time step $\tau\to 0$ (indeed $\bW(\bs^{n+1},\bs^n)=\mathcal O(\sqrt\tau)$ from the total square distance estimate \eqref{eq:tot_sq_dist}) and, in practice, only a very few inner timesteps are needed.

\subsubsection{Some properties of the approximate solution} 
As previously mentioned, the above Lagrangian framework can be practically implemented by simply projecting the (infinite dimensional) problem onto P1/P2 finite elements.
Provided that the iteration procedure (Steps 1 to 3 in Section~\ref{eq:subsec_algo_ALG2JKO}) converges as $k\to\infty$, as guaranteed from \cite{FG83}, the saddle-point $\bsigma=(\bs,\bm,\tilde\bs_1)$ satisfies by construction:
\begin{enumerate}[(i)]
 \item 
$(s_i,\bm_i)$ remains in the domain $\operatorname{Dom}(A)$ of the action functional $A$ defined in \eqref{eq:def_action_BB}; 
\item
the continuity equation $\partial_t s_{i,t}+\div(\bm_{i,t})=0$ holds with zero-flux boundary condition.
\end{enumerate}
As a consequence of (i) the scheme preserves the positivity, i.e., $s_i^{n+1}\geq 0$, whereas (ii) ensures the mass conservation $\int_\O s_i^{n+1}=\int_\O s_i^n$.

Moreover, the fully discrete ALG2-JKO scheme preserves by construction the gradient flow structure,
hence the scheme is automatically energy diminishing.
Since the energy functional \eqref{eq:Ee} includes the $\chi_\bDelta$ term accounting for the saturation constraint $\sum s_i=1$, one can and should include this convex indicator term in the discretized energy.
This contraint is then passed on to the proximal operator to be used in the implementation, see Section \ref{sec:numerics} for details.
As a result the saturation constraint is satisfied.

\subsection{Upstream mobility Finite Volume scheme}\label{ssec:Upwind}

The ALG2-JKO scheme described in the previous section will be compared to the 
widely used upstream mobility Finite Volume scheme~\cite{Peaceman77, BJ91, EHM03}.
As a first step, let us detail how $\O$ is discretized. 

\subsubsection{The finite volume mesh}

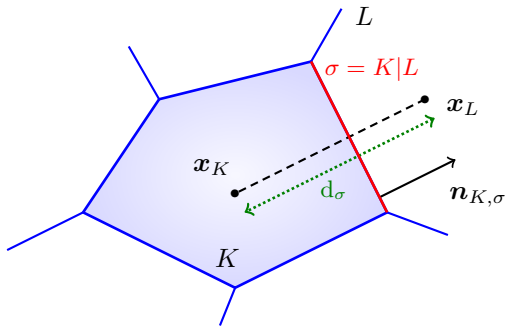
\begin{figure}[htb]
\begin{tikzpicture}[scale=1]
\shade[line width =1pt, outer color=blue!15, inner color=blue!3] (0,1) -- (2,-0.) -- (4,1) -- (3,3) -- (1,2.5) -- cycle;
\draw[line width=1pt, color=blue] (0,1) -- (2,-0.) -- (4,1) -- (3,3) -- (1,2.5) -- cycle;
\draw[line width =0.8pt, color=blue] (-1,0.5) -- (0,1);
\draw[line width =0.8pt, color=blue] (1.8,-0.5) -- (2,-0.);
\draw[line width =0.8pt, color=blue] (4.8,0.7) -- (4,1);
\draw[line width =0.8pt, color=blue] (3.4,3.7) -- (3,3);
\draw[line width =0.8pt, color=blue] (1,2.5) -- (0.6,3.2);
\draw[line width =0.8pt, densely dashed] (4.5,2.5) -- (2.,1.25);
\draw[<->,line width =1pt, densely dotted,color=green!50!black](4.625,2.25)--(2.125,1) ;
\draw[->,line width=0.8pt](3.9,1.2)--(4.9,1.7);
\draw[line width =1pt, color=red] (4,1) -- (3,3);
\path (4.5,2.5)  node{\scriptsize $\bullet$};
\path (2,1.25)  node{\scriptsize $\bullet$};
\path (5,2.4)  node{$\x_L$};
\path (1.7,1.6)  node{$\x_K$};
\path (1.9,0.4)  node{$K$};
\path (3.7,3.6)  node{$L$};
\path (3.8,2.9)  node{\small\color{red}$\sigma = {K|L}$};
\path (3.3,1.3)  node{\small\color{green!50!black}$\text{d}_{\sig}$};
\path (5.2,1.2) node {${\n}_{K,\sig}$};
\end{tikzpicture}
\caption{Here is an example of \emph{admissible mesh} in the sense of \cite{EGH00}}
\label{fig-am}
\end{figure}

The domain $\O$ is assumed to be polygonal. Then following \cite{EGH00}, an admissible 
mesh consists in a triplet $\left(\Ttt, \Eee,\left(\x_K\right)_{K\in\Ttt}\right)$. 
The elements $K$ of $\Ttt$ are open polygonal convex subsets of $\O$ called {\em control volumes}. 
Their boundaries are made of elements $\sig \in \Eee$ of codimension 1 (\emph{edges} if $d=2$ or \emph{faces} if $d=3$).
Let $K,L$ be two distinct elements of $\Ttt$, then $\ov K \cap \ov L$ is either empty, or reduced to a point (a vertex), 
or there exists $\sig \in \Eee$ denoted by $\sig = K|L$ such that $\ov K \cap \ov L = \ov \sig$. In particular, 
two control volumes share at most one edge. 
We denote by $\Eee_K = \left\{ \sig \in \Eee\; \middle| \; \bigcup_{\sig \in \Eee_K} \ov \sig = \p K \right\}$ 
the set of the edges associated to an element $K \in \Ttt$, and by 
$
\Nnn_K = \left\{ L \in \Ttt \; \middle| \; \text{there exists $\sig=K|L \in \Eee$}\; \right\}
$
the set of the neighboring control volumes to $K$. We also denote by 
\[
\Eee_{\rm ext} = \left\{ \sig \in \Eee \;\middle|\; \sig \subset \p\O\; \right\}, \qquad 
\Eee_{\rm int} = \Eee \setminus \Eee_{\rm ext}, \qquad 
\Eee_{\rm int,K} =  \Eee_{\rm int}\cap \Eee_K, \qquad \forall K \in \Ttt.
\]
The last element $\left(\x_K\right)_{K\in\Ttt}$ of the triplet 
corresponds to the so called \emph{cell-centers}. To each control volume $K \in \Ttt$, we associate an element 
$\x_K \in \O$ such that for all $L \in \Nnn_K$, the straight line $(\x_K,\x_L)$ is orthogonal to the edge ${K|L}$.
This implicitly requires that $\x_K$ and $\x_L$ are distinct, and we denote by $d_{\sig} = |\x_K - \x_L|$ 
for $\sig = K|L$ the distance between the cell centers of the neighboring control volumes $K$ and $L$. 
For $\sig \in  \Eee_K \cap \Eee_{\rm ext}$, we denote by $\x_\sig$ the projection of $\x_K$ on the hyperplane containing $\sig$, 
and by $d_\sig = |\x_K - \x_\sig|$. 
We also require that the vector $\x_L-\x_K$ is oriented 
in the same sense as the normal $\n_{K,\sig}$ to $\sig \in \Eee_K$ outward w.r.t. $K$.
We refer to Figure~\ref{fig-am} for an illustration of the notations used hereafter. 

Beyond cartesian grids, there are two classical ways to construct admissible meshes in the above sense when $d=2$. 
The first one consists in the classical Delaunay triangulation, the cell-center $\x_K$ of $K\in \Ttt$ 
being the center of the circumcircle of $K$. The second classical construction consists in choosing the cell centers $\left(\x_K\right)$ 
at first, and then to construct $\Ttt$ as the associated Vorono\"i diagram.

In what follows, we denote by $m_K$ the $d$-dimensional Lebesgue measure of the control volume $K\in\Ttt$, 
while $m_{\sig}$ denotes the $(d-1)$-dimensional Lebesgue measure of the edge $\sig \in \Eee$. 
We also denote by $a_{\sig} = \frac{m_{\sig}}{d_{\sig}}$ the {\em transmissivity} of the edge $\sig$.

In order to simplify the presentation, we restrict our presentation to the case of uniform time discretizations with time step $\tau>0$.
The extension to the case of time discretizations with varying time steps does lead to any particular difficulty. 

\subsubsection{Definition of the Finite Volume scheme}

The Finite Volume scheme relies on the discretization of the Euler-Lagrange equations~\eqref{eq:EL}
rather than on the minimizing movement scheme~\eqref{JKOscheme}. The main unknowns to the problems 
are located at the cell centers $\left(\x_K\right)_{K\in\Ttt}$. They consist in discrete saturations 
$s_{i,K}^n \simeq s_i(\x_K, n\tau)$ and discrete pressures $p_{i,K}^n \simeq p_i(\x_K,n\tau)$. 
In what follows, we denote by $\bs_K^n = \left(s_{i,K}^n\right)_{0 \leq i \leq N}$ (resp.
$\bp_K^n = \left(p_{i,K}^n\right)_{0 \leq i \leq N}$) and $\bs_\Ttt^n = \left(\bs_K^n\right)_{K\in\Ttt}$.

The first equation of the scheme is a straightforward consequence of~\eqref{eq:cons}, i.e., 
\begin{subequations}\label{eq:FVscheme}
\be\label{eq:cons_d}
\sum_{i=0}^N s_{i,K}^n = 1, \qquad \forall K \in \Ttt, \; \forall n \ge 1. 
\ee
This motivates the introduction of the discrete counterpart $\bXx_\Ttt$ of $\bXx$ defined by 
\[
\bXx_\Ttt = \left\{ 
\bs_{\Ttt}
\;\middle|\;
\bs_K \in \bDelta \;\text{for all}\; K \in \Ttt
\right\}, 
\]
so that \eqref{eq:cons_d} amounts to requiring that $\bs_\Ttt^n$ belongs to $\bXx_\Ttt$ for all $n$ (the nonnegativity 
of the saturations will be established later on).
The capillary pressure relations~\eqref{eq:capi} are discretized into 
\be\label{eq:capi_d}
p_{i,K}^n - p_{0,K}^n = \pi_i(s_K^n), \qquad \forall i \in \{1,\dots, N\}, \; \forall K \in \Ttt, \; \forall n \ge 1. 
\ee
Integrating~\eqref{eq:conv}  over the control volume $K \in \Ttt$ (recall here that the porosity $\omega$ was artificially 
set to $1$) and using Stokes' formula, one gets the natural approximation 
\be\label{eq:conv_d}
\frac{s_{i,K}^n - s_{i,K}^{n-1}}\tau m_K + \sum_{\sig \in \Eee_{K}}
s_{i,\sig}^n v_{i,K,\sig}^n = 0, 
\qquad \forall i \in \{0,\dots, N\}, \; \forall K \in \Ttt, \; \forall n \ge 1. 
\ee
Here, $v_{i,K,\sig}^n$ is an approximation of $\int_{\sig} \bv_i (\gamma, n\tau)\cdot \n_{K,\sig} \d\gamma$, where $\bv_i$ 
is related to $p_i$ through to Darcy law~\eqref{eq:Darcy}. Thanks to the orthogonality condition on the mesh, the choice
\be\label{eq:Darcy_d}
v_{i,K,\sig}^n = a_\sig \clem{\frac{\k}{\mu_i}} \left(p_{i,K}^n + \Psi_{i,K} - p_{i,L}^n - \Psi_{i,L} \right), \qquad \forall \sig = K|L \in \Eee_{\rm int},
\ee
is consistent --- we use the shortened notation $\Psi_{i,K} = \Psi_i(\x_K)$ ---.
In accordance with the no-flux boundary conditions, we impose that 
\[
v_{i,K,\sig}^n = 0, \qquad \forall \sig \in \Eee_K \cap \Eee_{\rm ext}, \; \forall n \ge 1. 
\]
It remains to define the approximate saturations $s^n_{i,\sig}$ for $\sig \in \Eee_{\rm int}$. We use here the 
very classical upwind choice~\cite{Peaceman77, BJ91, EHM03}, i.e., 
\be\label{eq:upwind}
s_{i,\sig}^n = \begin{cases}
\left(s_{i,K}^n\right)^+ & \text{if}\; v_{i,K,\sig}^n \ge 0, \\
\left(s_{i,L}^n\right)^+ & \text{otherwise}, 
\end{cases}
\qquad \forall \sig = K|L \in \Eee_{\rm int}.
\ee
Note that even though the mapping 
$
(\bs_\Ttt^n, \bp_\Ttt^n) \mapsto \bs_\Eee^n = \left(\left(s_{i,\sig}^n\right)_{0 \leq i \leq n} \right)_{\sig \in \Eee_{\rm int}}
$
is discontinuous, the quantity $s_{i,\sig}^n v_{i,K,\sig}^n$ depends in a continuous way of the main unknowns. 
\end{subequations}

The scheme~\eqref{eq:FVscheme} amounts to a nonlinear system of equations to be solved a 
each time step. This will be practically done thanks to Newton-Raphson method. But before, we 
establish some properties of the FV scheme, namely the energy decay, the entropy control, 
the non-negativity of the saturations, or the existence of a solution $\left(\bs_\Ttt^n, \bp_\Ttt^n\right)$ to the scheme.

\subsubsection{Some properties of the approximate solution}

The first key property of the FV scheme that we point out is the non-negativity of the saturations:
\[
s_{i,K}^n \ge 0 \qquad \forall i \in \{0,\dots, N\}, \; \forall K \in \Ttt, \; \forall n \ge 1. 
\]
In order to establish this estimate, it suffices to rewrite~\eqref{eq:conv_d} as 
\[
s_{i,K}^n  + \frac{\tau}{m_K} \sum_{\substack{\sig \in \Eee_{{\rm int}, K}\\\sig = K|L}} 
\left[\left(s_{i,K}^n\right)^+ \left(v_{i,K,\sig}^n\right)^+ - 
 \left(s_{i,L}^n\right)^+ \left(v_{i,K,\sig}^n\right)^- \right]= s_{i,K}^{n-1}
 \]
thanks to~\eqref{eq:upwind}. In the previous expression, we used the convention $a^- = \max(0, -a) \geq 0$. 
Assume for contradiction that $s_{i,K}^n$ is negative, then so does the left-hand side,
while the right-hand side is nonnegative by induction. \clem{Together with~\eqref{eq:cons_d},} this shows that 
\be\label{eq:bXx_d} 
\bs_\Ttt^n \in \bXx_\Ttt, \qquad \forall n \geq 1. 
\ee

The scheme is mass conservative for the $N+1$ phases since 
\[
v_{i,K,\sig}^n + v_{i,L,\sig}^n = 0 \quad \text{hence}\quad 
s_{i,\sig}^n v_{i,K,\sig}^n + s_{i,\sig}^n  v_{i,L,\sig}^n = 0, \qquad \text{for}\; \sig = K|L.
\] 
Together with the no-flux boundary conditions, this shows that the mass is conserved along time:
\be\label{eq:mass_d}
\sum_{K \in \Ttt} s_{i,K}^n m_K = \sum_{K \in \Ttt} s_{i,K}^{n-1} m_K = \sum_{K \in \Ttt} s_{i,K}^0 m_K, 
\qquad \forall n \geq 1, \; \forall i \in \{0, \dots, N\}. 
\ee
Then the discrete solution $\bs_\Ttt^n$ remains in the discrete counterpart $\bAa_\Ttt$ of $\bAa$ defined as 
the elements $\bs_\Ttt$ of $\R_+^\Ttt$ such that $\sum_{K \in \Ttt} s_{i,K} m_K =  \sum_{K \in \Ttt} s_{i,K}^0 m_K$ 
for all $i \in \{0,\dots, N\}$.

Multiplying the scheme~\eqref{eq:conv_d} by $\tau\left(p_{i,K}^n + \Psi_{i,K}\right)$ and summing over $K \in \Ttt$ yields
\begin{multline*}
\sum_{i=0}^N \sum_{K\in\Ttt} \left(s_{i,K}^n - s_{i,K}^{n-1} \right) \left(p_{i,K}^n + \Psi_{i,K}\right) m_K \\
+ \tau \sum_{i=0}^N \frac{\k}{\mu_i}\sum_{\substack{\sig \in \Eee_{\rm int}}} a_\sig  s_{i,\sig}^n \left(p_{i,K}^n + \Psi_{i,K} - p_{i,L}^n - \Psi_{i,L}\right)^2 = 0.
\end{multline*}
The second term in the above expression is clearly nonnegative. concerning the first term, one can use the constraint~\eqref{eq:cons_d} to rewrite as
\begin{align*}
\sum_{i=0}^N \sum_{K\in\Ttt} \left(s_{i,K}^n - s_{i,K}^{n-1} \right) p_{i,K}^n m_K 
= &
\sum_{i=1}^N \sum_{K\in\Ttt} \left(s_{i,K}^n - s_{i,K}^{n-1} \right) \left(p_{i,K}^n - p_{0,K}^n\right) m_K \\
\geq & \sum_{K\in\Ttt} \left(\Pi(\bs_K^{n,*}) - \Pi(\bs_K^{n-1,*}) \right) m_K,
\end{align*}
the last inequality being a consequence of the convexity of $\Pi$. 
This establishes that the scheme is energy diminishing: denoting by 
\[
\Ee(\bs_\Ttt^n) = \sum_{K\in \Ttt} \left(\Pi(\bs_K^{n,*}) + \sum_{i=0}^N s_{i,K}^n \Psi_{i,K} \right)m_K, \qquad n \geq 0, 
\]
one has 
\be\label{eq:NRJ_d}
\Ee(\bs_\Ttt^n) + \tau \sum_{i=0}^N \frac{\k}{\mu_i}\sum_{\substack{\sig \in \Eee_{\rm int}}} a_\sig  s_{i,\sig}^n \left(p_{i,K}^n + \Psi_{i,K} - p_{i,L}^n - \Psi_{i,L}\right)^2\leq \Ee(\bs_\Ttt^{n-1}), \qquad \forall n \geq 1. 
\ee

The last {\em a priori} estimate we want to point out is the discrete counterpart of the flow interchange estimate. 
It is obtained by multiplying \eqref{eq:conv_d} by $\tau \mu_i \log(s_{i,K}^n)$ and by summing over $i \in \{0,\dots,N\}$ 
and $K\in \Ttt$, leading to 
\be\label{eq:entro_d-0}
\sum_{i=0}^N  \mu_i \sum_{K\in\Ttt}\left(s_{i,K}^n - s_{i,K}^{n-1} \right)\log(s_{i,K}^n) m_K 
+\tau  \sum_{i=0}^N {\k}\sum_{\substack{\sig \in \Eee_{\rm int}}} a_\sig  s_{i,\sig}^n v_{i,K,\sig}^n \left( \log(s_{i,K}^n) - 
\log(s_{i,L}^n) \right) = 0.
\ee
As already discussed in Remark~\ref{rem:log_test_function} this corresponds to taking $\log s_i$ as a test-function in the weak formulation of the continuous PDEs.
The first term of~\eqref{eq:entro_d-0} can be estimated thanks to an elementary convexity inequality
\[
\sum_{i=0}^N  \mu_i \sum_{K\in\Ttt}\left(s_{i,K}^n - s_{i,K}^{n-1} \right)\log(s_{i,K}^n) m_K
\geq 
\Hh(\bs_\Ttt^n) - \Hh(\bs_\Ttt^{n-1}), \qquad \forall n \geq 1
\]
with 
\[
\Hh(\bs_\Ttt^n) = \sum_{i=0}^N \mu_i \sum_{K\in\Ttt} \left(h(s_{i,K}^n) - h(s_{i,K}^{n-1})\right) m_K, 
\qquad 
h(s) = s\log(s) - s + 1 \geq 0. 
\]
Note that the entropy functional $\Hh$ is bounded on $\bXx_\Ttt$.
The second term of~\eqref{eq:entro_d-0} can be estimated as follows. First, the concavity of $s \mapsto \log(s)$ 
yields
\[
s_{i,L}^n \left( \log(s_{i,K}^n) - \log(s_{i,L}^n) \right) \leq s_{i,K}^n - s_{i,L}^n \leq s_{i,K}^n \left( \log(s_{i,K}^n) - \log(s_{i,L}^n) \right), 
\qquad \sig = K|L, 
\]
so that the upwind choice~\eqref{eq:upwind} for $s_{i,\sig}^n$ ensures that 
\[ 
a_\sig  s_{i,\sig}^n v_{i,K,\sig}^n \left( \log(s_{i,K}^n) - \log(s_{i,L}^n) \right) \geq a_\sig v_{i,K,\sig}^n (s_{i,K}^n - s_{i,L}^n), 
\qquad \sig = K|L. 
\]
Using the expression~\eqref{eq:Darcy_d} of $v_{i,K,\sig}^n$ and the relation~\eqref{eq:cons_d} on the saturations, one gets that 
\[
\sum_{i=0}^N \sum_{\substack{\sig \in \Eee_{\rm int}}} a_\sig  s_{i,\sig}^n v_{i,K,\sig}^n \left( \log(s_{i,K}^n) - 
\log(s_{i,L}^n) \right) \geq A+B, 
\]
where
\begin{align*}
A = & \sum_{i=1}^N \sum_{\substack{\sig \in \Eee_{\rm int} \\ \sig = K|L}} 
	a_\sig (\pi_i(\bs_K^{n,*}) - \pi_i(\bs_L^{n,*}))(s_{i,K}^n - s_{i,L}^n),  \\
B = & \sum_{i=0}^N \sum_{\substack{\sig \in \Eee_{\rm int} \\ \sig = K|L}} 
	a_\sig (\Psi_{i,K}- \Psi_{i,L})(s_{i,K}^n - s_{i,L}^n) = 
	\sum_{i=0}^N \sum_{K \in \Ttt}  s_{i,K}^n \sum_{L\in\Nnn_K} a_\sig (\Psi_{i,K} - \Psi_{i,L}). 
\end{align*}
Recalling the definition~\eqref{eq:Psi_i} of the external potential and denoting by $\Psi_{i,\sig}= \Psi_i(\x_\sig)$, one has 
\[
\sum_{L \in \Nnn_K} a_\sig (\Psi_{i,K} - \Psi_{i,L}) + \sum_{\sig \in \Eee_{{\rm ext}} \cap \Eee_K}a_\sig (\Psi_{i,K}- \Psi_{i,\sig}) = 0. 
\]
Since $0 \leq s_{i,K}^n \leq 1$, this implies that 
\[
B \geq \tau \k |\p\O| |\g|. 
\]
On the other hand, the assumption~\eqref{eq:H-Pi} on the capillary pressure potential ensures that 
\[
A \geq 
\frac1{\varpi} \sum_{i=1}^N \sum_{\substack{\sig \in \Eee_{\rm int} \\ \sig = K|L}} 
	a_\sig \left(  \pi_i(\bs_K^{n,*}) - \pi_i(\bs_L^{n,*})\right)^2, \qquad \forall \sig = K|L \in \Eee_{\rm int}. 
\]
Hence collecting the previous inequalities in~\eqref{eq:entro_d-0} provides the following discrete $L^2_{\rm loc}(H^1)$-estimate 
one the capillary pressures
\be\label{eq:pi_L2H1_d}
\sum_{n=1}^M \tau  \sum_{i=1}^N \sum_{\substack{\sig \in \Eee_{\rm int} \\ \sig = K|L}} 
a_\sig \left(  \pi_i(\bs_K^{n,*}) - \pi_i(\bs_L^{n,*})\right)^2 \leq C (1+M \tau). 
\ee
Clearly, \eqref{eq:pi_L2H1_d} is the discrete counterpart of the estimate~\eqref{eq:pi_L2H1} obtained thanks to the 
flow interchange technique. 
The derivation of a discrete $L^2_{\rm loc}(H^1)$ estimate on the phase pressures 
from~\eqref{eq:pi_L2H1_d} and~\eqref{eq:NRJ_d} requires one additional assumption on the 
capillary pressure functions $\left(\pi_i\right)_{1\leq i \leq n}$. More precisely, we assume that 
\be\label{eq:H-pi-strong}
\text{
$\pi_i$ only depends on $s_i$:
}\; \frac{\p}{\p s_j} \pi_i (\bs^*) = 0\; \text{if $i\neq j$}. 
\ee
Since $\Pi$ is convex, the functions $\pi_i$ are increasing. Assumption~\eqref{eq:H-pi-strong} is needed to establish that, 
at least for fine enough grids, there holds 
\[
\sum_{i = 0}^N s_{i,\sig}^n \geq \alpha >0, \qquad \forall n \ge 1, \; \forall \sig \in \Eee_{\rm int},
\]
for some uniform $\alpha$. Thanks to this estimate, one can follow the lines of \cite[Proposition 3.4 \& Corollary 3.5]{CGM17} 
(see also~\cite{CN_FVCA8})
to derive the estimate 
\be\label{eq:L2H1_p_d}
\sum_{n=1}^M \tau  \sum_{0=1}^N \sum_{\substack{\sig \in \Eee_{\rm int} \\ \sig = K|L}} 
a_\sig \left(  p_{i,K}^{n} - p_{i,L}^{n})\right)^2 \leq C (1+M \tau). 
\ee
The phase pressures being defined up to an additive constant (recall that they are related to Kantorovich potentials), 
one has to fix this degree of freedom. This can be done by enforcing 
\[
\sum_{K\in\Ttt} p_{0,K}^n m_K = 0, \qquad \forall n \geq 1. 
\]

Based on the {\em a priori} estimates~\eqref{eq:bXx_d} and~\eqref{eq:L2H1_p_d}, 
we can make use of a topological degree argument (see for instance~\cite{Dei85}) to claim 
that there exists (at least) one solution to the scheme. Moreover, assuming some classical regularity on the mesh $\Ttt$ 
(see for instance~\cite{Ahmed_intrusion}), one can prove the piecewise constant approximate solutions converge 
towards a weak solution when the size of the mesh $\Ttt$ and the time step $\tau$ tend to $0$. 
This convergence results together with the properties~\eqref{eq:bXx_d}--\eqref{eq:L2H1_p_d}
as well as the wide popularity of this scheme in the engineering community makes this scheme 
a reference for solving~\eqref{eq:EL}. In the next section, we show that the ALG2-JKO scheme 
presented in Section~\ref{ssec:ALG2} produces very similar results: same qualitative results, 
conservation of the mass of each phase  and preservation of the positivity.

\section{Numerical experiments}\label{sec:numerics}

In this section, we compare the numerical results produced by the ALG2-JKO scheme 
presented in Section~\ref{ssec:ALG2} with the upstream mobility Finite Volume scheme of Section~\ref{ssec:Upwind}. 
In the sequel the regularization parameter $r$ introduced in the augmented Lagrangian formulation~\eqref{aug lagrangian} is fixed to $r =1$ for simplicity, which gives satisfactory numerical results.
The case of a three phase flow (typically water, oil and gas) is 
presented in Section~\ref{ssec:3phases}, whereas a two-phase flow is simulated in Section~\ref{ssec:2phases}. 
In both cases, we do not have analytical solutions at hand and the results are compared 
thanks to snapshots. 

\clem{
Note the both time discretizations are of order 1. The extension to order two methods is a challenging task. Concerning the ALG2-JKO scheme, 
one possibility could be to use the order 2 approximation based on the midpoint rule proposed in~\cite{LT17}, but there is no rigorous 
foundation to this work  up to now as far as we know. An alternative approach would be to use the variational BDF2 approach 
proposed in~\cite{MP17_arXiv}. But the variational problem to be solved at each time step is no longer convex-concave, 
so that its practical resolution becomes more involving. 
Concerning the finite volume scheme, there is (up to our knowledge) no time integrator of order 2 that ensures the decay of a general energy. 
Going to higher order time discretizations yields also difficulties concerning the preservation of the positivity. This explains why the 
backward Euler scheme is very popular in the context of the simulation of multiphase porous media flows. 
}

\subsection{Two-phase flow with Brooks-Corey capillarity}\label{ssec:2phases}

As a first example we consider a two-phase flow, where water ($s_0$) and oil ($s_1$) are competing within the background porous medium.
For the capillary pressure, we choose the very classical Brooks-Corey (or Leverett) model
\begin{equation}
\label{eq:B-C capillarity}
 p_1 - p_0 = \pi_1(s_1)= \alpha(1- s_1)^{-1/2}.
\end{equation}
We refer to~\cite{BB90} for an overview of the classical capillary pressure relation for two-phase flows.

As in Section~\ref{ssec:model}, the corresponding energy reads explicitly
$$
\Ee(s_0,s_1)= \int_\Omega \Psi_0 s_0 + \int_\Omega \Psi_1 s_1 - 2\alpha \int_\Omega (1- s_1)^{1/2} + \int_\Omega \chi_\bDelta (s_0,s_1).
$$
As already mentioned, only the second subproblem \eqref{JKO step2 second subproblem} in step 2 of the ALG2-JKO algorithm depends on the choice of the energy functional.
For the above particular case, this reads: for each $\x \in\Omega$ and setting $\overline{\boldsymbol{c}}:= -\boldsymbol{\phi}^{k+1}(1,\x) + \tilde{\bs}_1^k(\x)$, solve
$$ 
\boldsymbol{c}^{k+1}(\x) = \argmin_{\boldsymbol{c}\in {\R}^3}\,\left\{ \frac{1}{2}  |\boldsymbol{c} -  \overline{\boldsymbol{c}} |^2 + E_\tau^*(\x,\boldsymbol{c})\right\},
$$
where $E_\tau^*(\x,\cdot)$ is the Legendre transform of $E_\tau(\x,\cdot)$ defined by
$$
  E_\tau(\x,c_0,c_1) = \tau \Psi_0(\x) c_0 + \tau \Psi_1(\x) c_1 -2 \tau\alpha (1- c_1)^{1/2} + \chi_\bDelta (c_0,c_1) \text{ for all } c_0,c_1 \in \R.
 $$
This minimization problem is equivalent to computing
$$ 
\boldsymbol{c}^{k+1}(\x) = \text{Prox}_{E_\tau^*(\x,\cdot)} ( \overline{\boldsymbol{c}}),
$$
where the proximal operator $\text{Prox}_f$ of a given convex, lower semicontinuous function $f \, : \, {\R}^{N+1} \rightarrow \R \cup  \{+\infty\}$ is defined by
$$ 
\text{Prox}_f(\overline{\by}) := \argmin_{\by \in {\R}^{N+1} } \left\{
\frac{1}{2}|\by -\overline{\by}|^2 +f(\by) \right\},
\qquad \forall\, \overline{\by} \in {\R}^{N+1}.
$$
Thanks to Moreau's identity
\begin{equation}
\label{Moreau's identity}
\text{Prox}_{f^*}(\overline{\by}) = \overline{\by} - \text{Prox}_{f}(\overline{\by}) \qquad \forall\, \overline{\by} \in {\R}^{N+1},
\end{equation} 
it suffices to compute $\text{Prox}_{E_\tau}$ in order to determine $\text{Prox}_{E_\tau^*}$, and we never actually compute the Legendre transform $E_\tau^*(\x,\cdot)$.
Computing the proximal operator $\boldsymbol{c}^{k+1}(\x) =\text{Prox}_{E_\tau^*(\x,\cdot)} ( \overline{\boldsymbol{c}})$ thus amounts to evaluating
$$
(c_0^{k+1}(\x),c_1^{k+1}(\x)) = ( \overline{c}_0,\overline{c}_1) - \text{Prox}_{E_\tau(\x, \cdot)} (\overline{c}_0,\overline{c}_1).
$$
Finally, $(\tilde{c}_0,\tilde{c}_1):=\text{Prox}_{E_\tau(\x, \cdot)} (\overline{c}_0,\overline{c}_1) $ is computed by solving
$$ 
\tilde{c}_1=\argmin_{0 \leq c_1 \leq 1}\,
\left\{\frac{1}{2}|c_1 + \overline{c}_0 - \tau \Psi_0(\x) -1 |^2 + \frac{1}{2}|c_1 - \overline{c}_1 + \tau \Psi_1(\x) |^2 -2\tau \alpha (1 -c_1)^{1/2}
\right\}
$$
and then setting $\tilde{c}_0=1-\tilde{c}_1$.
More explicitly, $\tilde{c}_1$ is the positive part of the root on $(-\infty, 1)$ of 
$$ 
2c - \overline{c}_1 + \tau \Psi_1(\x) +\overline{c}_0 - \tau \Psi_0(\x) -1 +\frac{\tau \alpha}{(1 -c )^{1/2}} =0.
$$
To conclude, we set $(c_0^{n+1}(\x),c_1^{n+1}(\x)) = ( \overline{c}_0 - \tilde{c}_0,\overline{c}_1-\tilde{c}_1)$.

\begin{figure}[!htbp]
\centering
\begin{subfigure}{\textwidth}
\centering
\includegraphics[width=5cm]{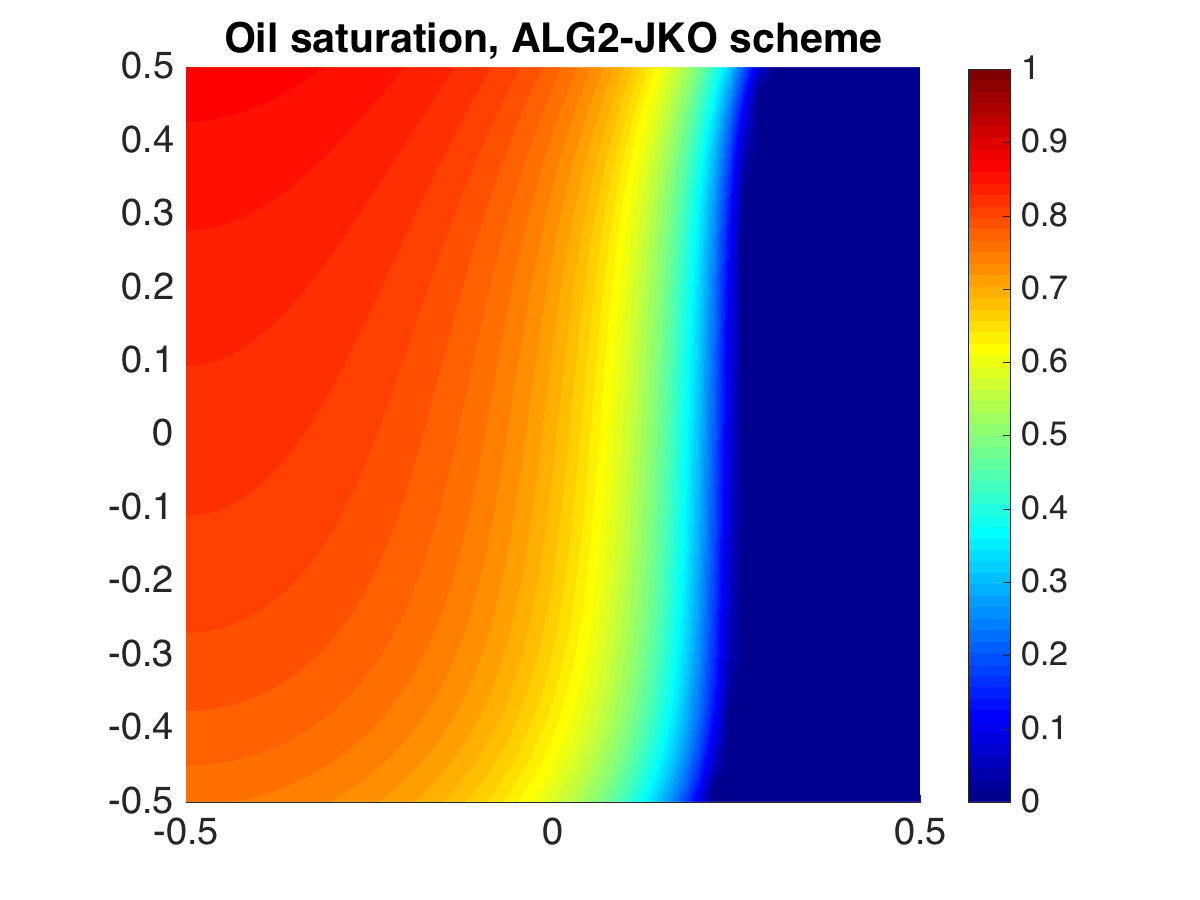}
\hspace{.1cm}
\includegraphics[width=5cm]{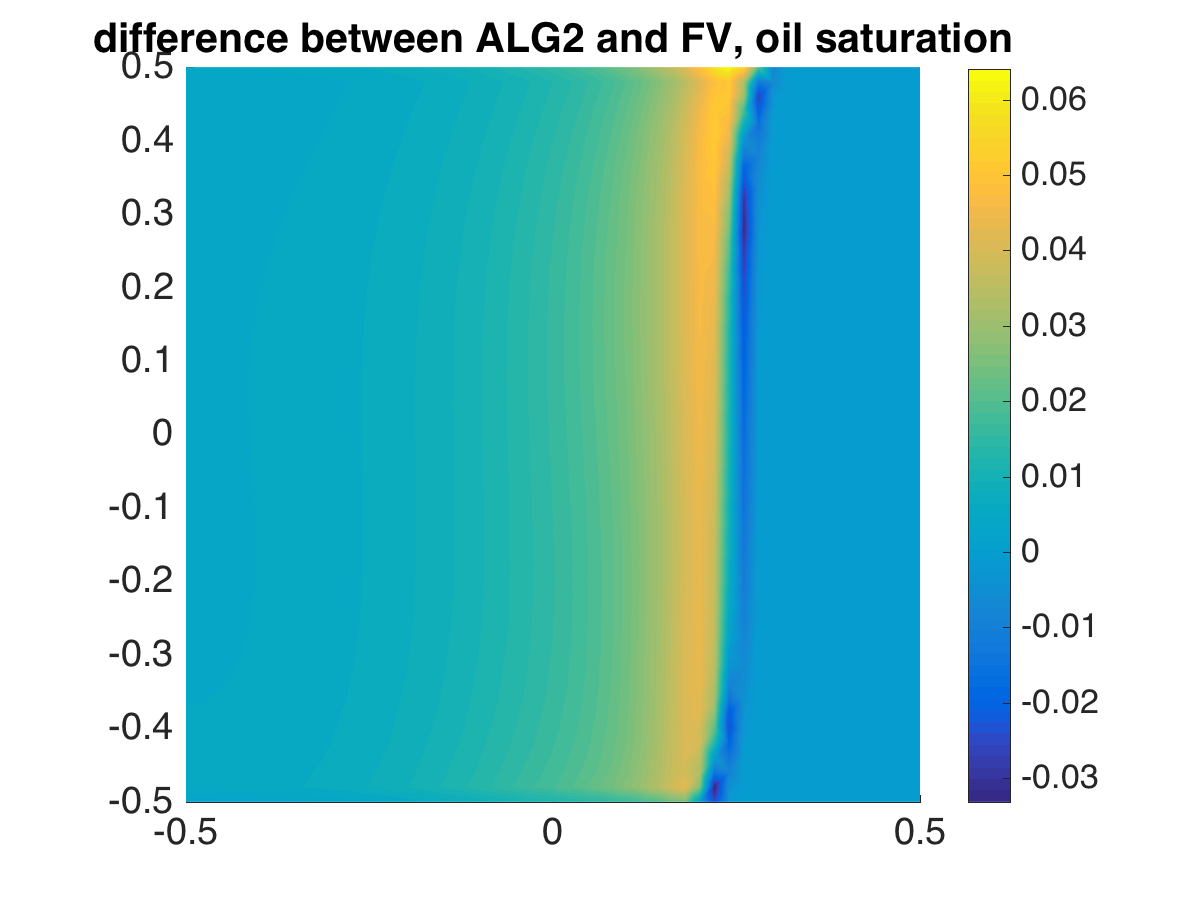}
\caption{$t=2.5$}
\end{subfigure}
\medskip

\begin{subfigure}{\textwidth}
\centering
\includegraphics[width=5cm]{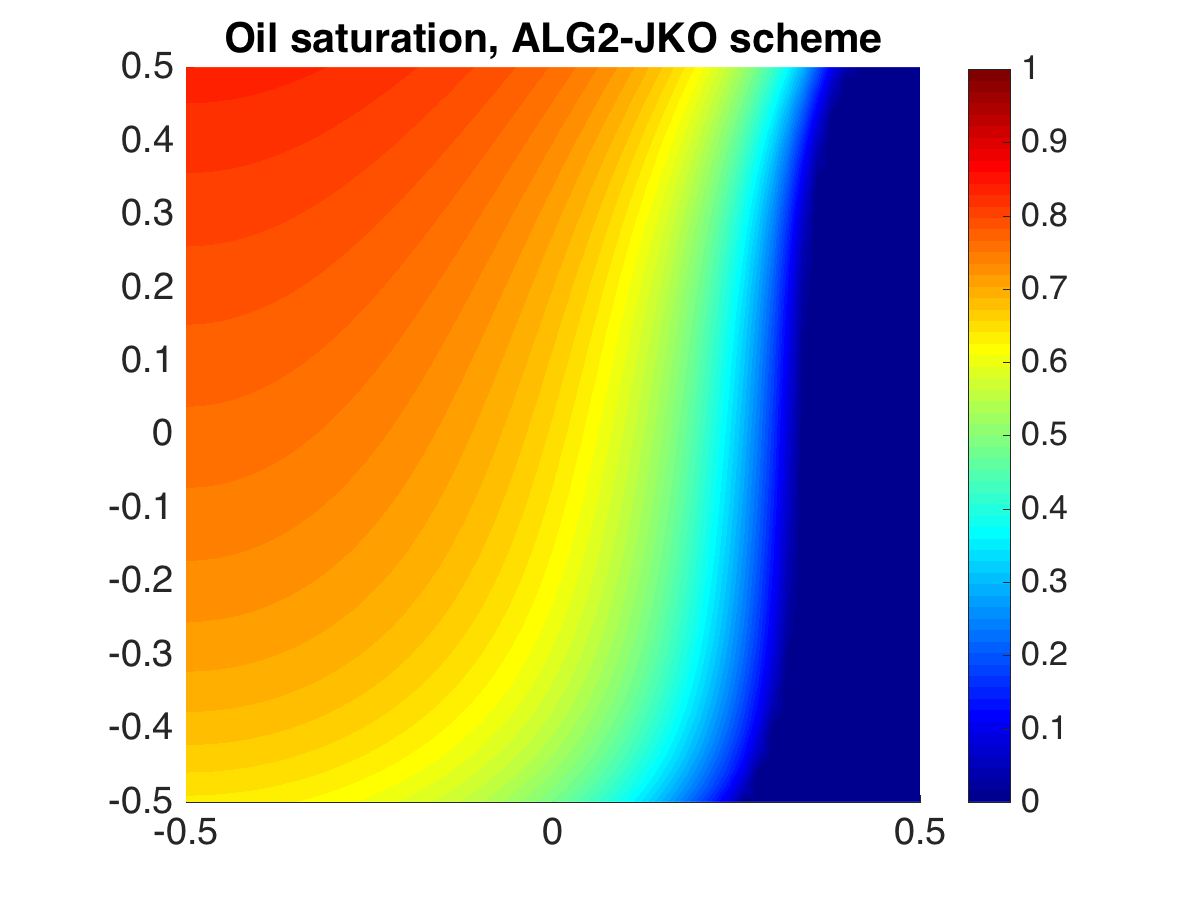}
\hspace{.1cm}
\includegraphics[width=5cm]{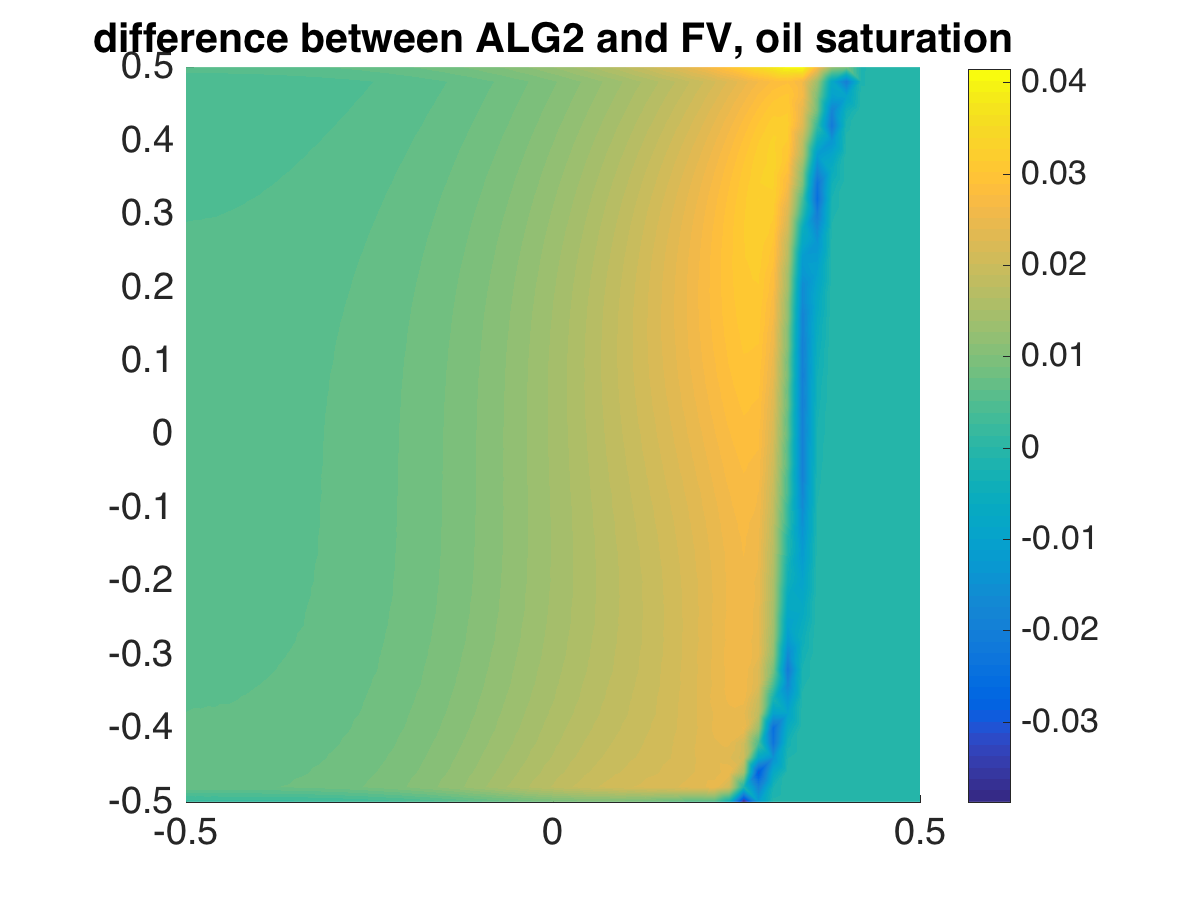}
\caption{{$t=5$}}
\end{subfigure}
\medskip

\begin{subfigure}{\textwidth}
\centering
\includegraphics[width=5cm]{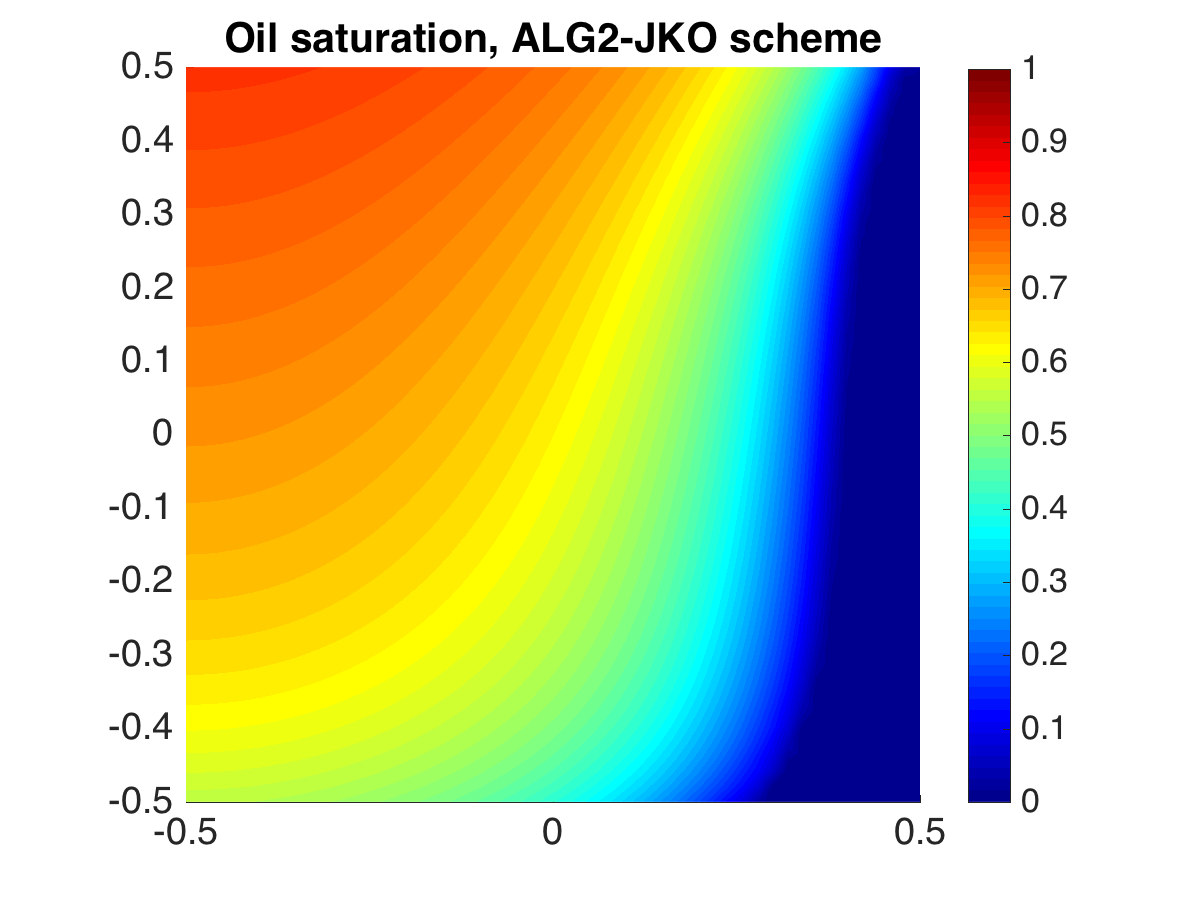}
\hspace{.1cm}
\includegraphics[width=5cm]{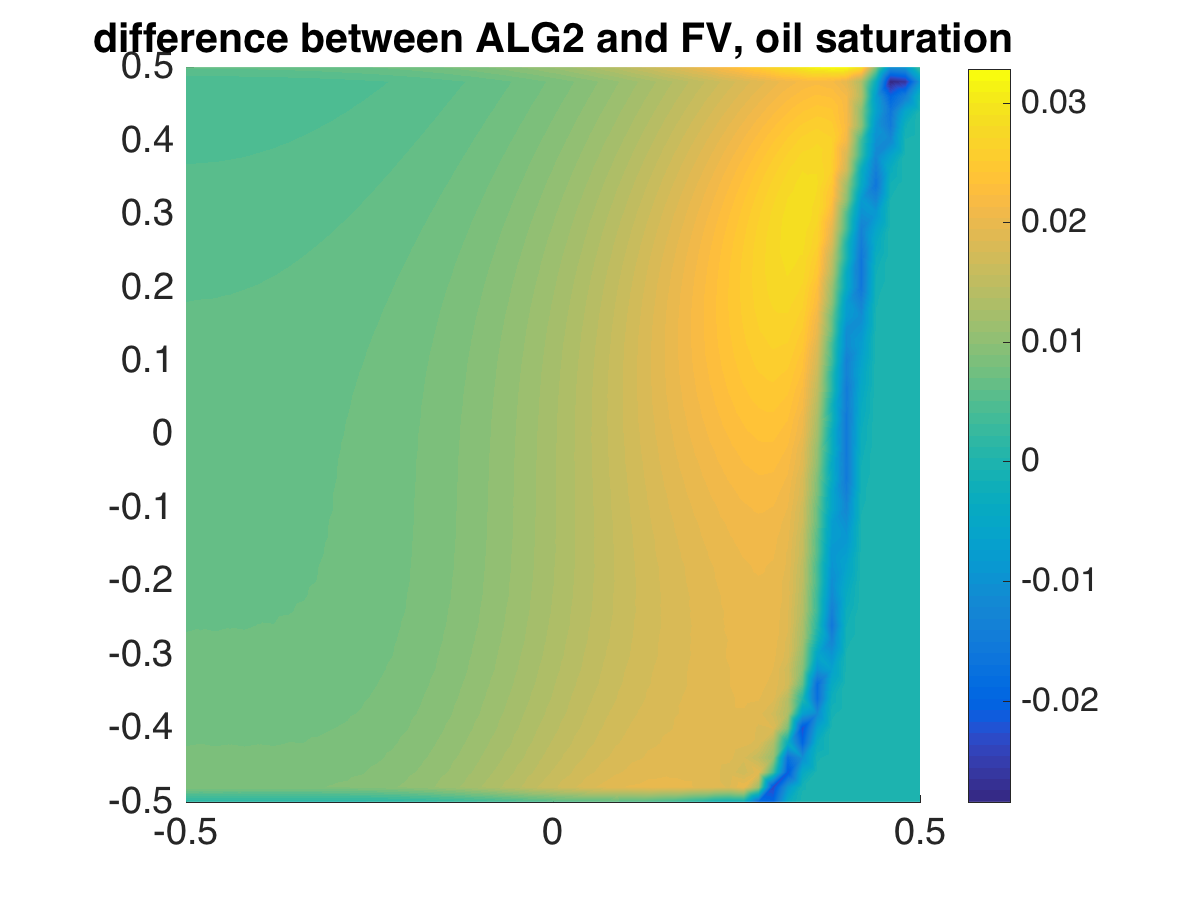}
\caption{$t=7.5$}
\end{subfigure}
\medskip

\begin{subfigure}{\textwidth}
\centering
\includegraphics[width=5cm]{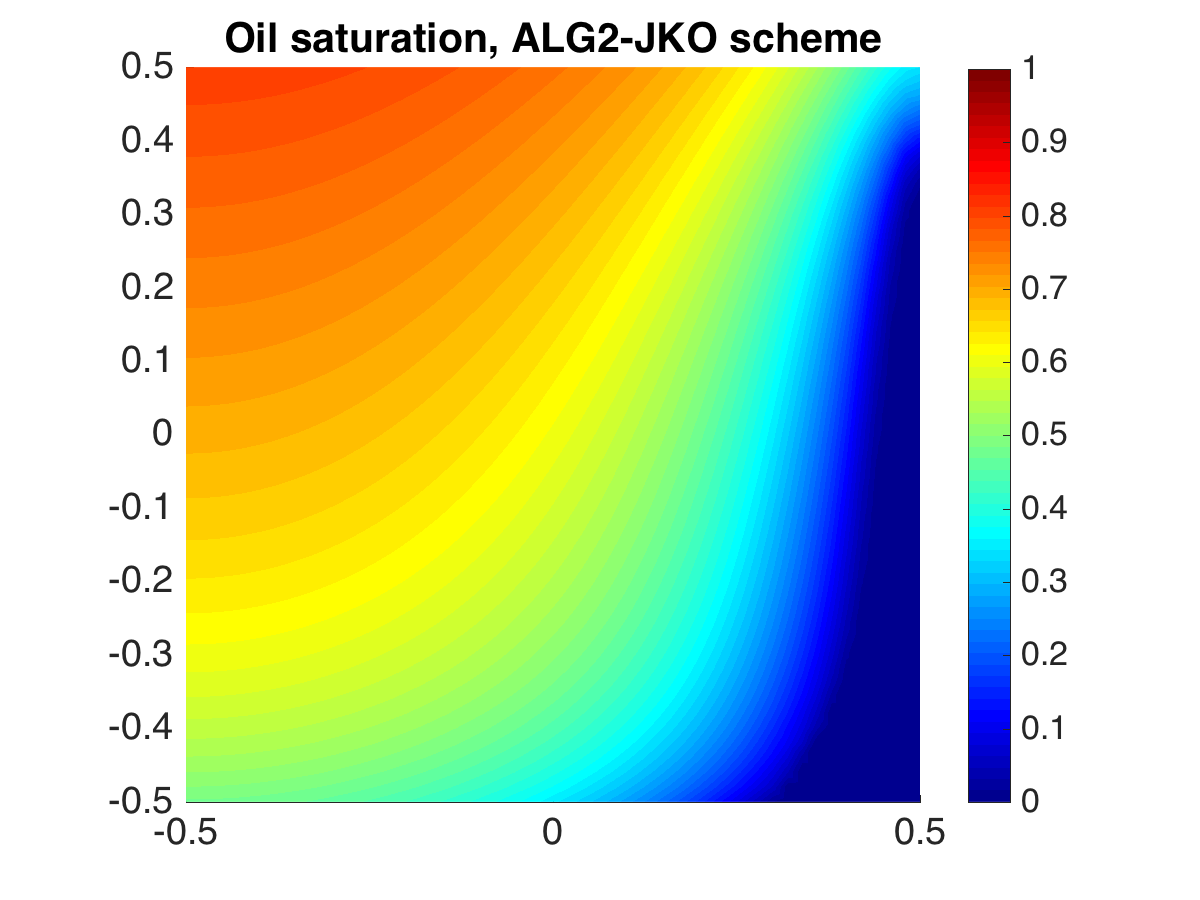}
\hspace{.1cm}
\includegraphics[width=5cm]{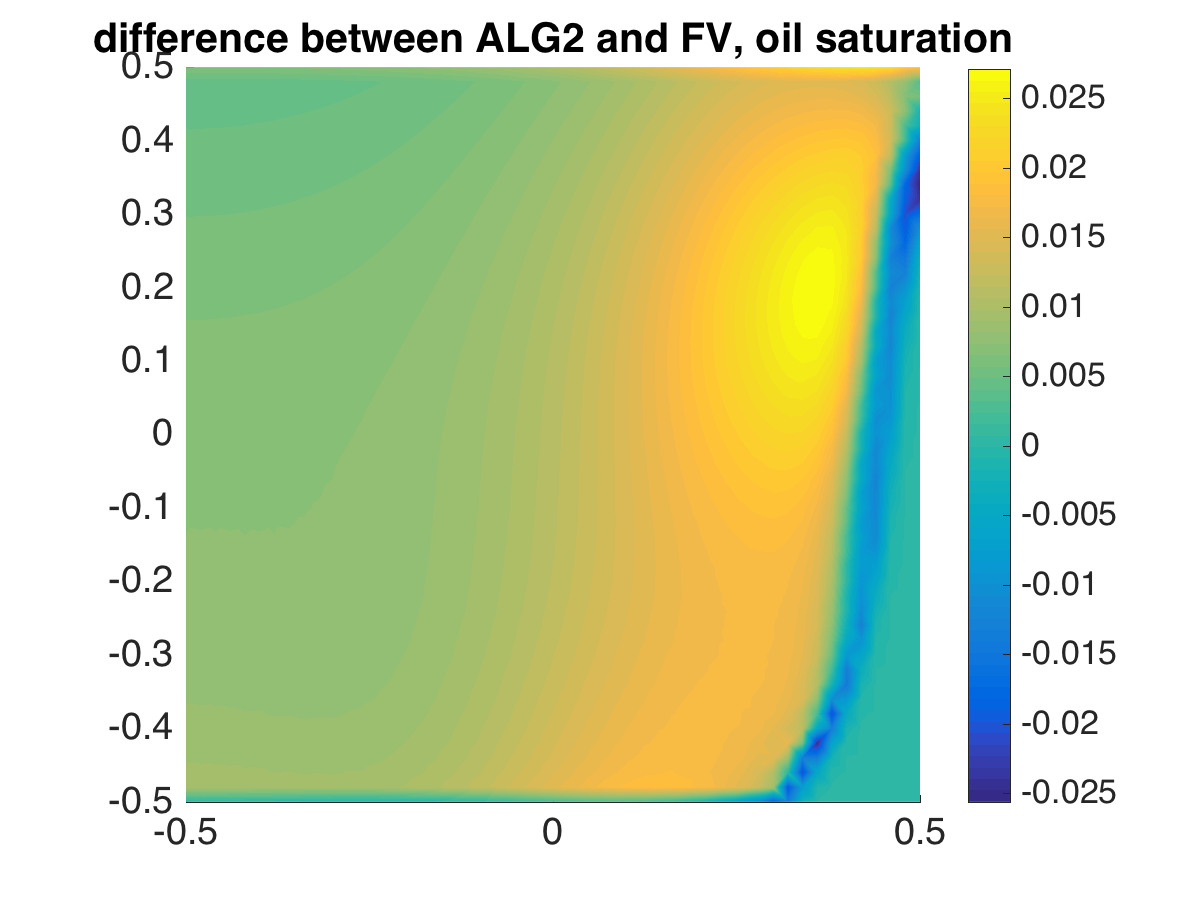}
\caption{$t=10$}
\end{subfigure}
\caption{Oil saturation for the two-phase flow problem with Brooks-Corey capillary pressure function~\eqref{eq:B-C capillarity}, $\alpha = 1$: 
numerical solution provided by the ALG2-JKO scheme (left) and difference between the ALG2-JKO approximate solution and the upstream mobility finite volume approximation solution (right).}
\label{fig:2phases-BC}
\end{figure}

On Figure \ref{fig:2phases-BC}, we compare the numerical solutions of problem~\eqref{eq:EL} with Brooks-Corey capillarity~\eqref{eq:B-C capillarity} obtained thanks to the ALG2-JKO scheme and to the upstream mobility finite volume scheme. Simulations with the ALG2-JKO scheme are carried using a structured grid with 5000 triangles and 2601 vertices in space and a single inner time step, and with $200$ JKO steps ($\tau=0.05$).
Simulations with the upstream mobility finite volume scheme are performed on \clem{the corresponding Cartesian grid with 2500 squares}. 
The time step $\tau$ appearing in~\eqref{eq:conv_d} can be also set to $0.05$ here since Newton's method converges rather easily in this test case. 

%

\begin{figure}[htb]
\centering
\includegraphics[width=5cm]{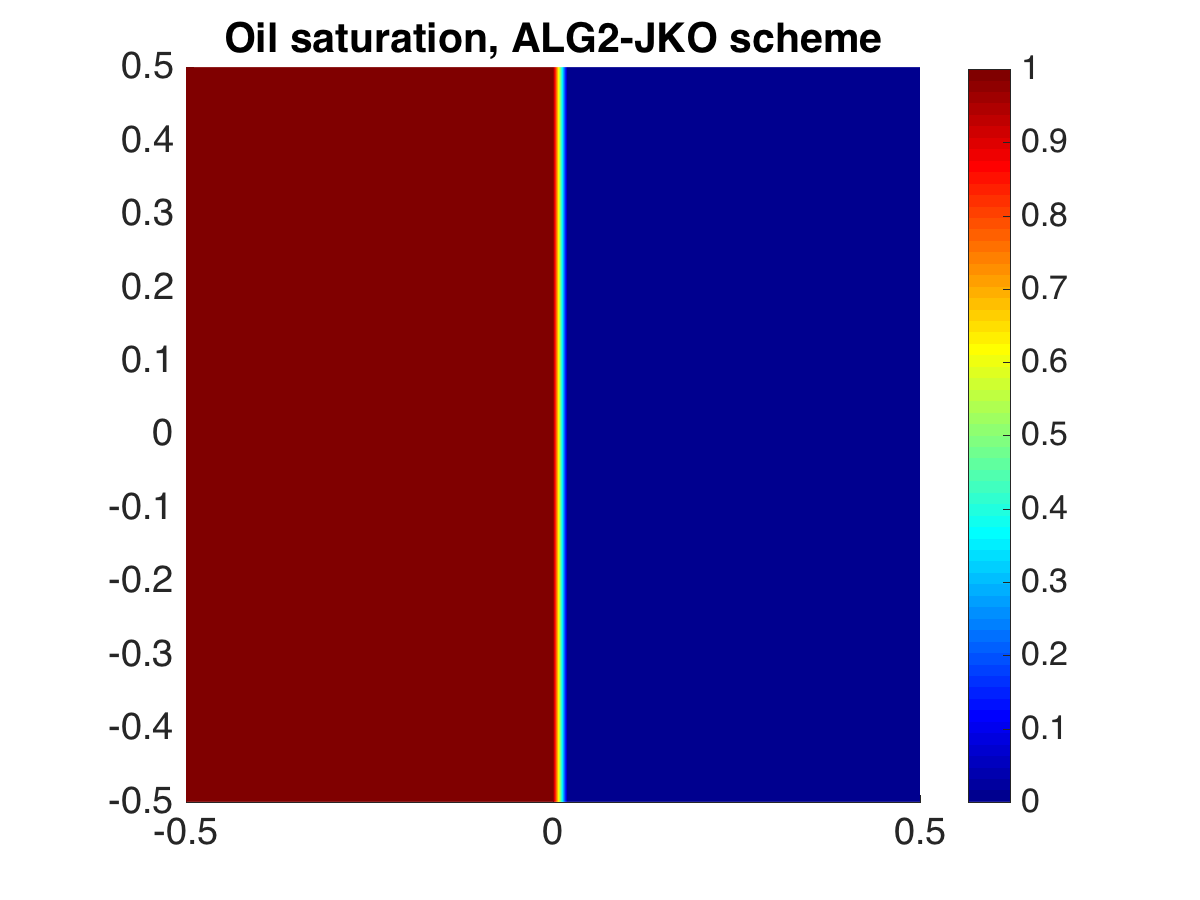}
\hspace{.1cm}
\includegraphics[width=5cm]{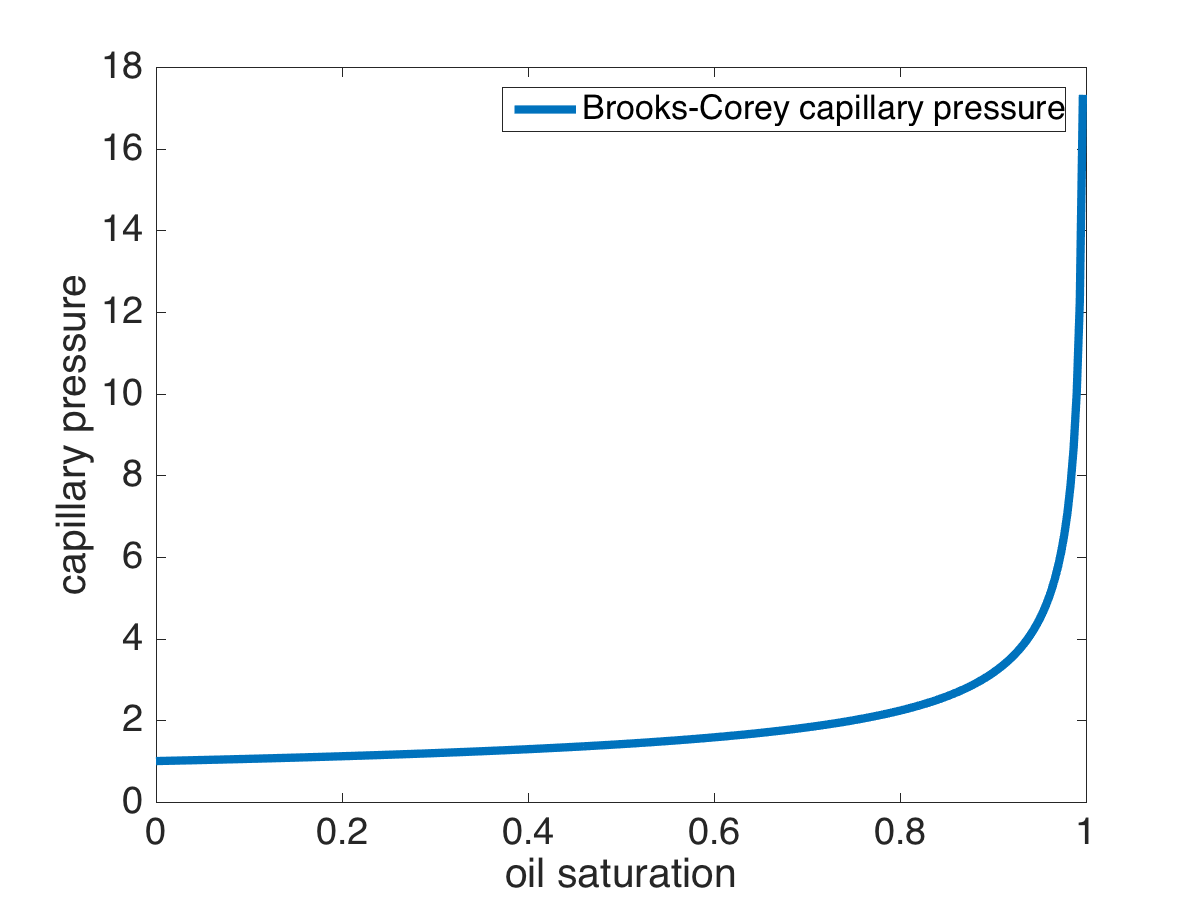}
\caption{{Two-phase flow: initial oil saturation profile (left) and Brooks-Corey capillary pressure function~\eqref{eq:B-C capillarity} with $\alpha = 1$ (right).}}
\end{figure}

As expected, the results produced by the two schemes are very similar. The dense phase (the water) is instantaneously diffused 
in the whole domain because of the singularity of $\pi_1$ near $1$. When time goes, oil slowly moves to the top because of 
buoyancy.

\subsection{Three-phase flow with quadratic capillary potential}\label{ssec:3phases}

In the second test case, we consider the case of a three-phase flow where 
water ($s_0$), oil ($s_1$), and gas ($s_2$) are in competition within the porous medium.
Here we assume that the capillary pressure functions $\pi_1$ and $\pi_2$ are linear,
$$
p_1 -p_0 = \pi_1(s_1) = \alpha_1 s_1
\qquad\text{ and }\qquad
p_2 -p_0 = \pi_2(s_2) = \alpha_2 s_2.
$$
The corresponding capillary potential $\Pi$ is then given by 
\[
\Pi(\bs^*) = \frac{\a_1}2 (s_1^2) +  \frac{\a_2}2 (s_2^2).
\]
The Assumption~\eqref{eq:H-Pi} and \eqref{eq:H-pi-strong} are fulfilled, so that 
we are in the theoretical framework of our statements, i.e., convergence of the minimizing movement scheme
and of the finite volume scheme.
However, the problem is difficult to simulate because of the rather large ratios on the viscosities.  
Indeed, the phase $0$ represents water, the phase $1$ corresponds to oil and the phase $2$ corresponds to gas, 
and we set 
\[
\mu_0 = 1, \quad \mu_1 = 50, \quad \mu_2 = 0.1, 
\qquad \text{and} \qquad
\rho_0 = 1, \quad \rho_1 = 0.87, \quad \rho_2 = 0.1. 
\]

The resulting energy in the JKO scheme~\eqref{JKOscheme} is given by 
$$ 
\Ee(s_0, s_1, s_2) := \sum_{i=0}^2 \int_\Omega \Psi_i s_i + \frac{\alpha_1}{2}\int_\Omega s_1^2 + \frac{\alpha_2}{2}\int_\Omega s_2^2 +\int_\Omega \chi_{\bDelta}(s_0, s_1, s_2),
$$
and we denote accordingly, for $\x\in\O$ and $\boldsymbol{c}=(c_0,c_1,c_2)\in\R^3$
$$
 E_\tau(\x,\boldsymbol{c}):=\sum_{i=0}^2 \tau\Psi_i(\x) c_i + \frac{\tau\alpha_1}{2} c_1^2 + \frac{\tau\alpha_2}{2} c_2^2 +\chi_{\bDelta}(\boldsymbol{c}).
 $$
Setting again $\overline{\boldsymbol{c}}= -\boldsymbol{\phi}^{k+1}(1,\x) + \tilde{\bs}_1^k(\x)$ and taking advantage of Moreau's identity \eqref{Moreau's identity}, the second subproblem \eqref{JKO step2 second subproblem} of step 2 is equivalent to, for all $\x \in {\R}^d$, 
$$ 
\boldsymbol{c}^{k+1}(\x) = \overline{\boldsymbol{c}} -  \text{Prox}_{E_\tau(\x, \cdot)} (\overline{\boldsymbol{c}} ).
$$
Evaluating the proximal operator $\tilde{\boldsymbol{c}}:=\text{Prox}_{E_\tau(\x, \cdot)}(\overline{\boldsymbol{c}})$ is equivalent to solving
\begin{multline}
\label{eq:prox 3d}
(\tilde{c}_1, \tilde{c}_2)= \argmin_{0\leq c_i \leq 1, \, 0 \leq c_1+c_2 \leq 1}
\Bigg\{\sum_{i=1}^2 \left( \frac{1}{2}|c_i - \overline{c}_i +\tau \Psi_i(\x)|^2 + \tau \frac{\alpha_i}{2}c_i^2 \right) \\
+ \frac{1}{2}|c_1 + c_2 +\overline{c}_0- \tau \Psi_0(\x) -1 |^2
\Bigg\},
\end{multline}
with $\tilde{c}_0=1-\tilde{c}_1- \tilde{c}_2$.
The solution $(u_1,u_2)$ of the unconstrained version of \eqref{eq:prox 3d} is explicitly given by
$$ 
u_1 = \frac{(2+ \tau  \alpha_2)\gamma_1 -\gamma_2}{(2+  \tau\alpha_1)(2+\tau \alpha_2) -1} \text{ and }u_2 = \frac{(2+ \tau \alpha_1)\gamma_2 -\gamma_1}{(2+\tau  \alpha_1)(2+ \tau \alpha_2) -1},
$$
where $\gamma_i:= \overline{c}_i - \tau\Psi_i(\x) - \overline{c}_0 + \tau \Psi_0(\x) +1$.
If $(u_1,u_2) \in \boldsymbol{\Delta}^*$ then $(\tilde{c}_1, \tilde{c}_2)=(u_1,u_2)$ is the true solution of \eqref{eq:prox 3d}, and $\tilde{c}_0=1-u_1 -u_2$.
Otherwise, one should seek for the minimizer of \eqref{eq:prox 3d} on the boundary $\partial\boldsymbol{\Delta}^*=\{s_1=0,0\leq s_2\leq 1\}\cup\{0\leq s_1\leq 1,s_2=0\}\cup\{s_1+s_2=1\}$.
This leads to three easy minimization problems that can be again solved explicitly, and we omit the details.
To conclude, the update of $\boldsymbol{c}^{k+1}(\x)$ is given by $ \boldsymbol{c}^{k+1}(\x)= \overline{\boldsymbol{c}} - \tilde{\boldsymbol{c}}$.

\begin{figure}[!htbp]
\centering
\begin{subfigure}{\textwidth}
\centering
\includegraphics[width=5cm]{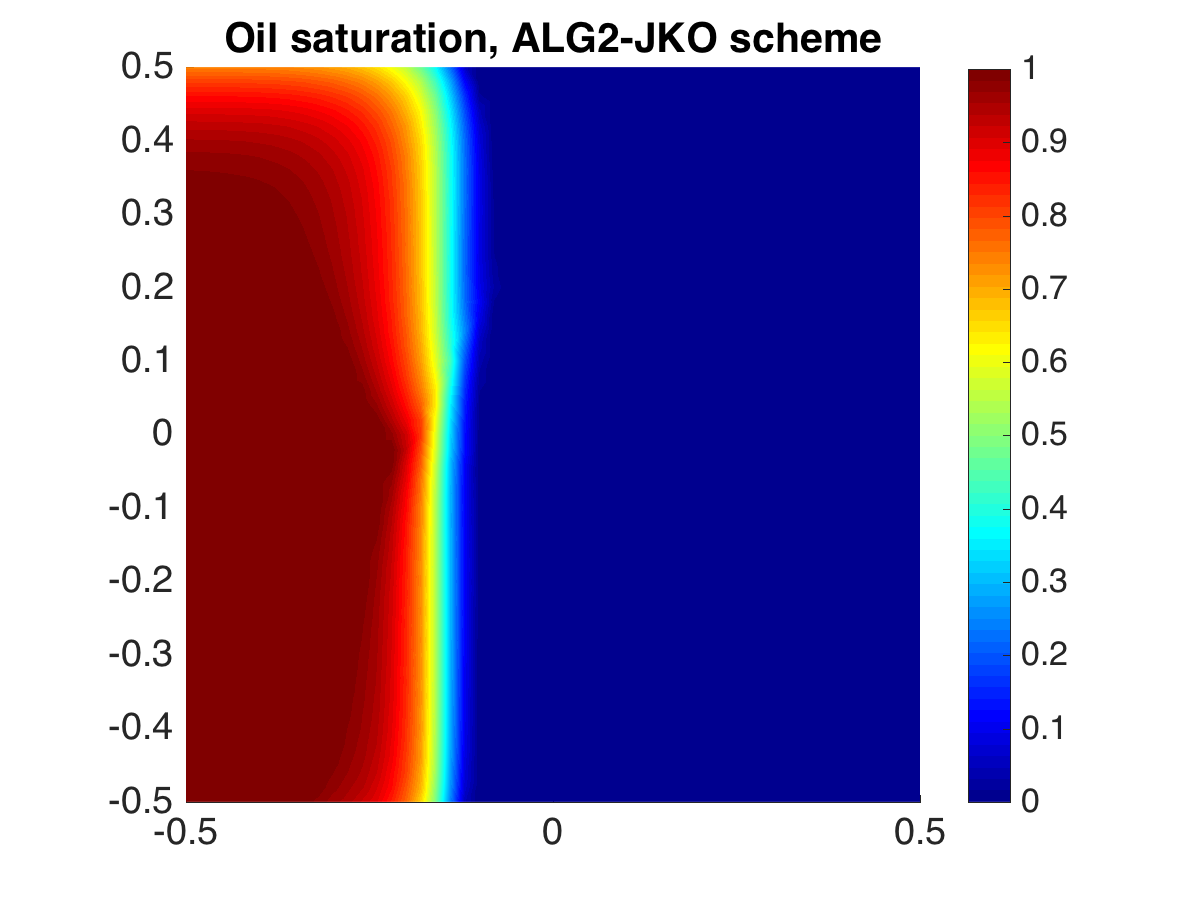} 
\hspace{.1cm}
\includegraphics[width=5cm]{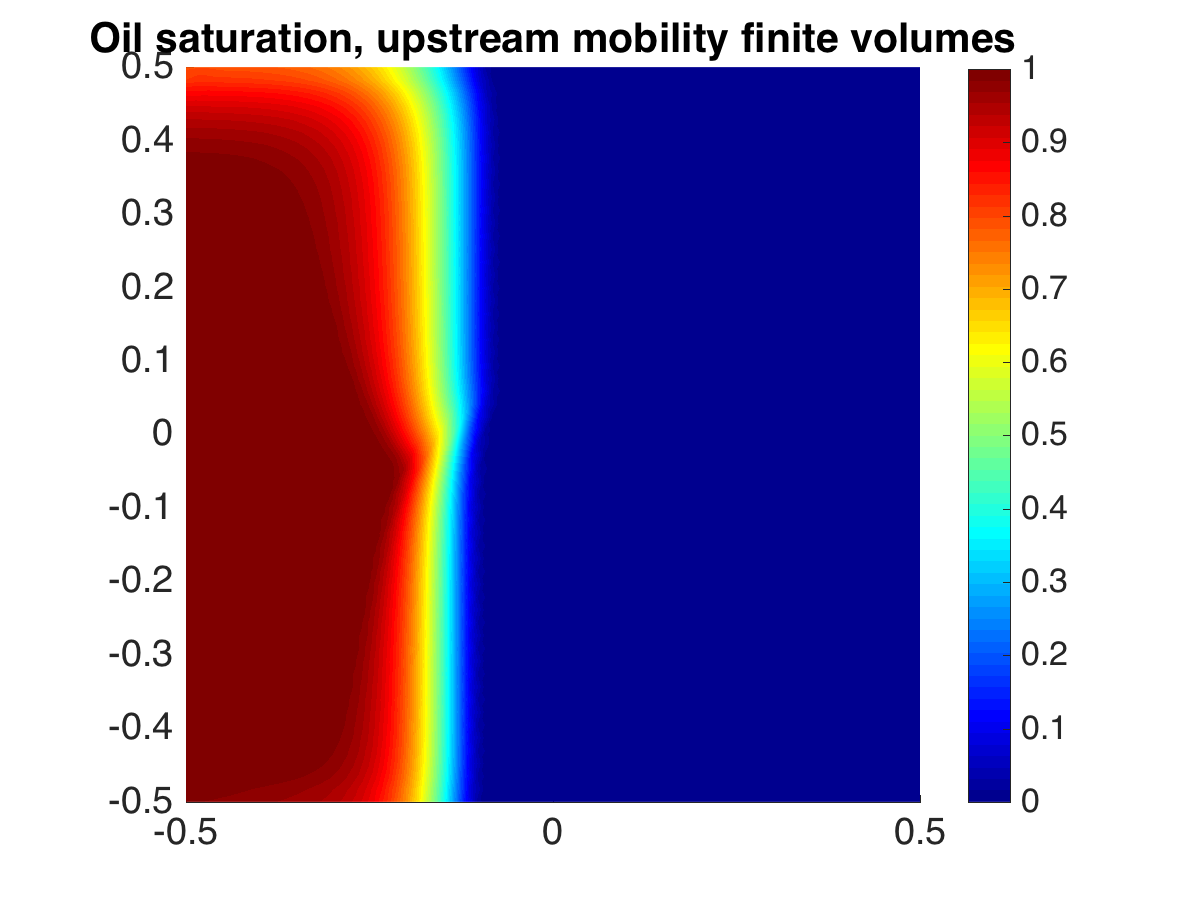} 
\caption{t=0.1}
\end{subfigure}
\medskip

\begin{subfigure}{\textwidth}
\centering
\includegraphics[width=5cm]{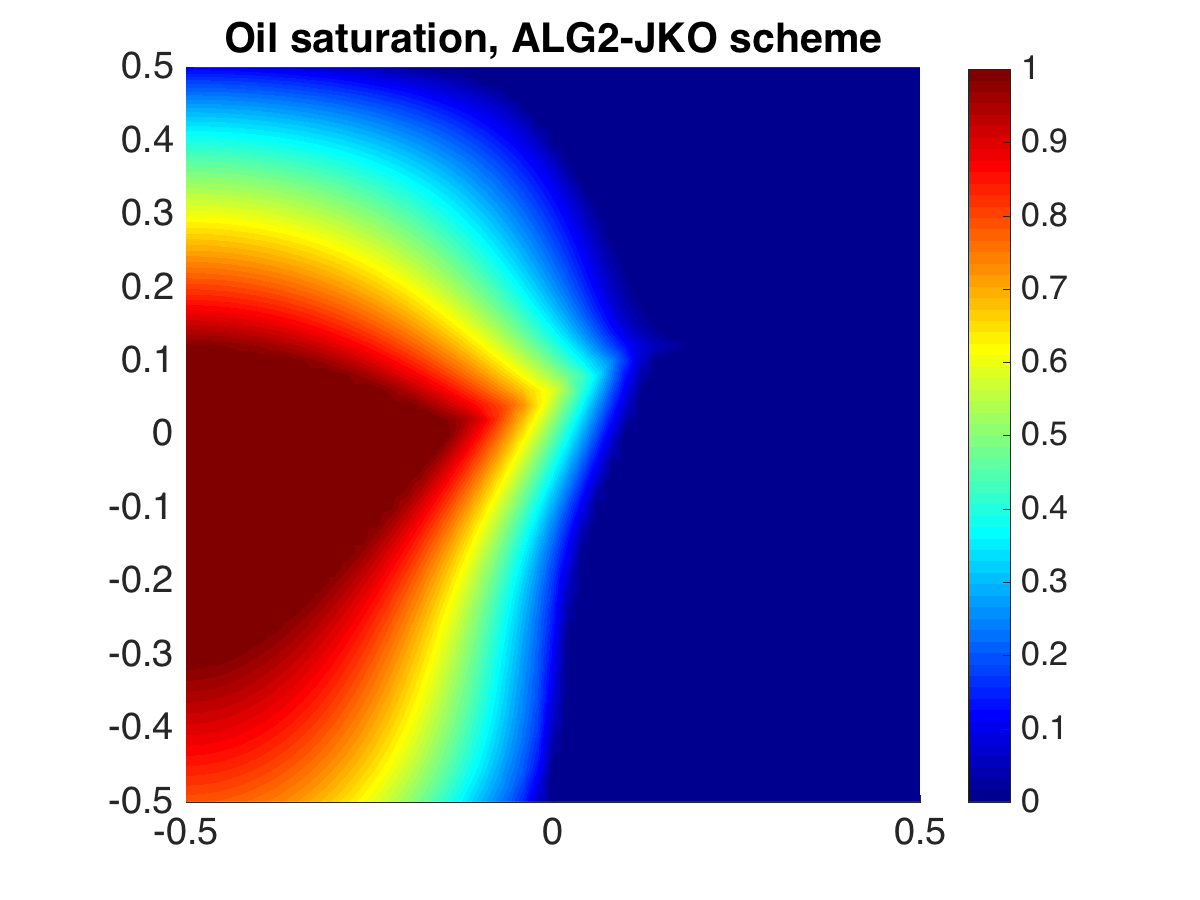} 
\hspace{.1cm}
\includegraphics[width=5cm]{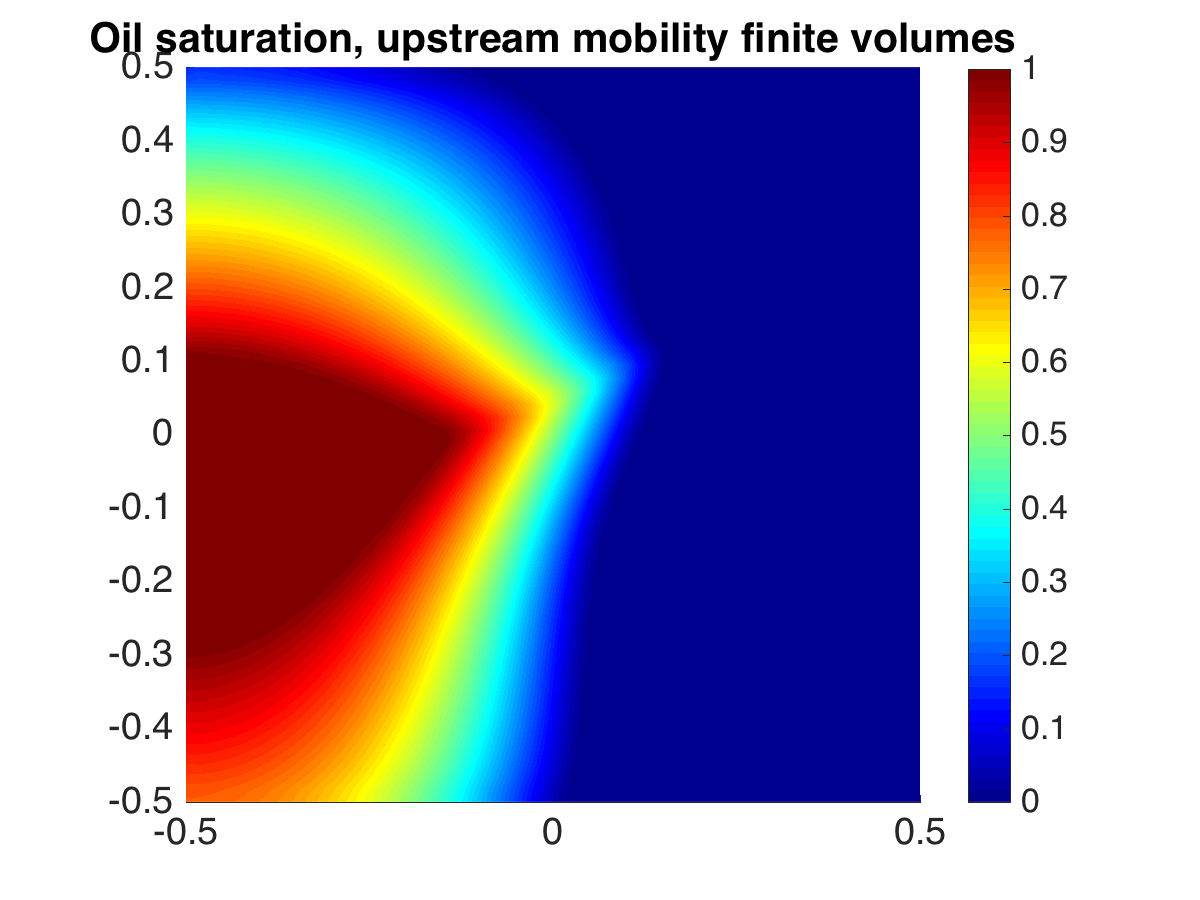} 
\caption{t=1.25}
\end{subfigure}
\medskip

\begin{subfigure}{\textwidth}
\centering
\includegraphics[width=5cm]{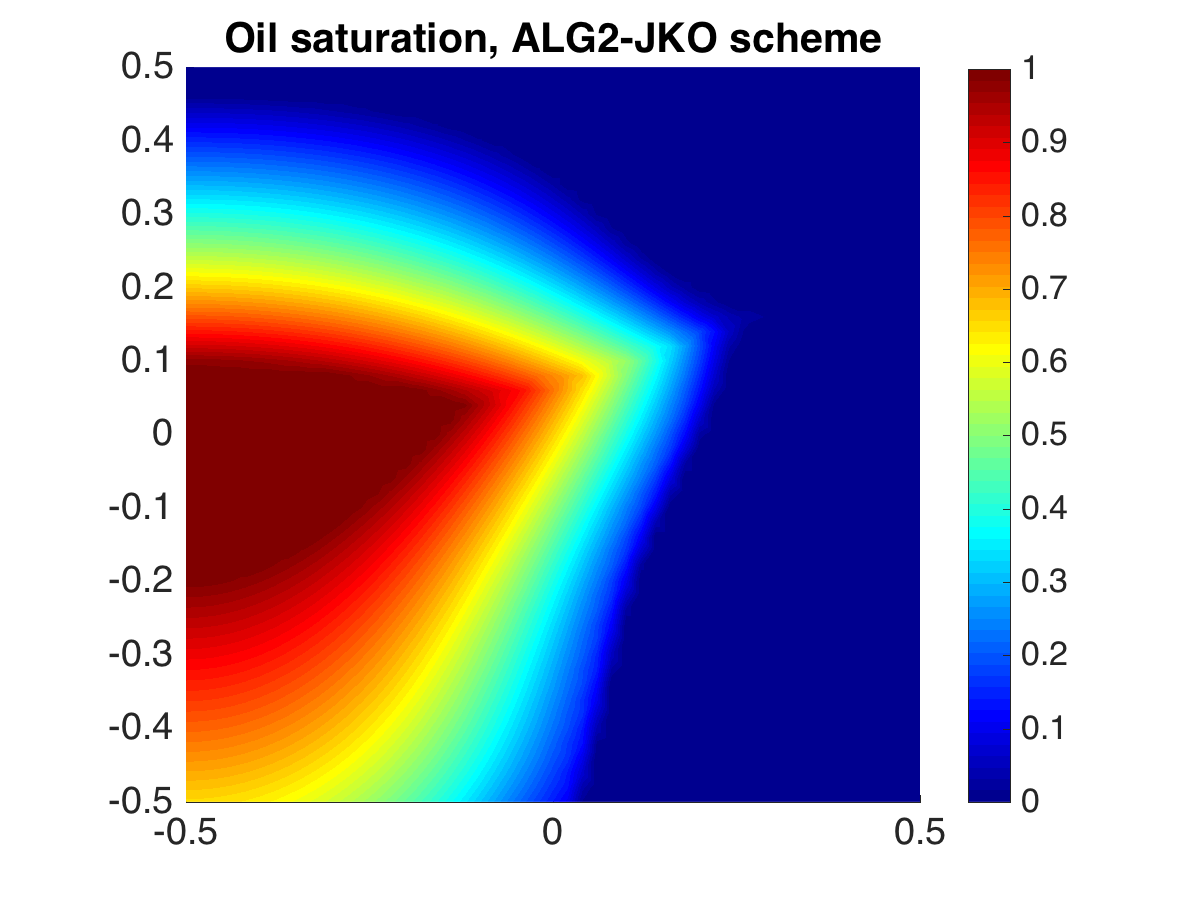} 
\hspace{.1cm}
\includegraphics[width=5cm]{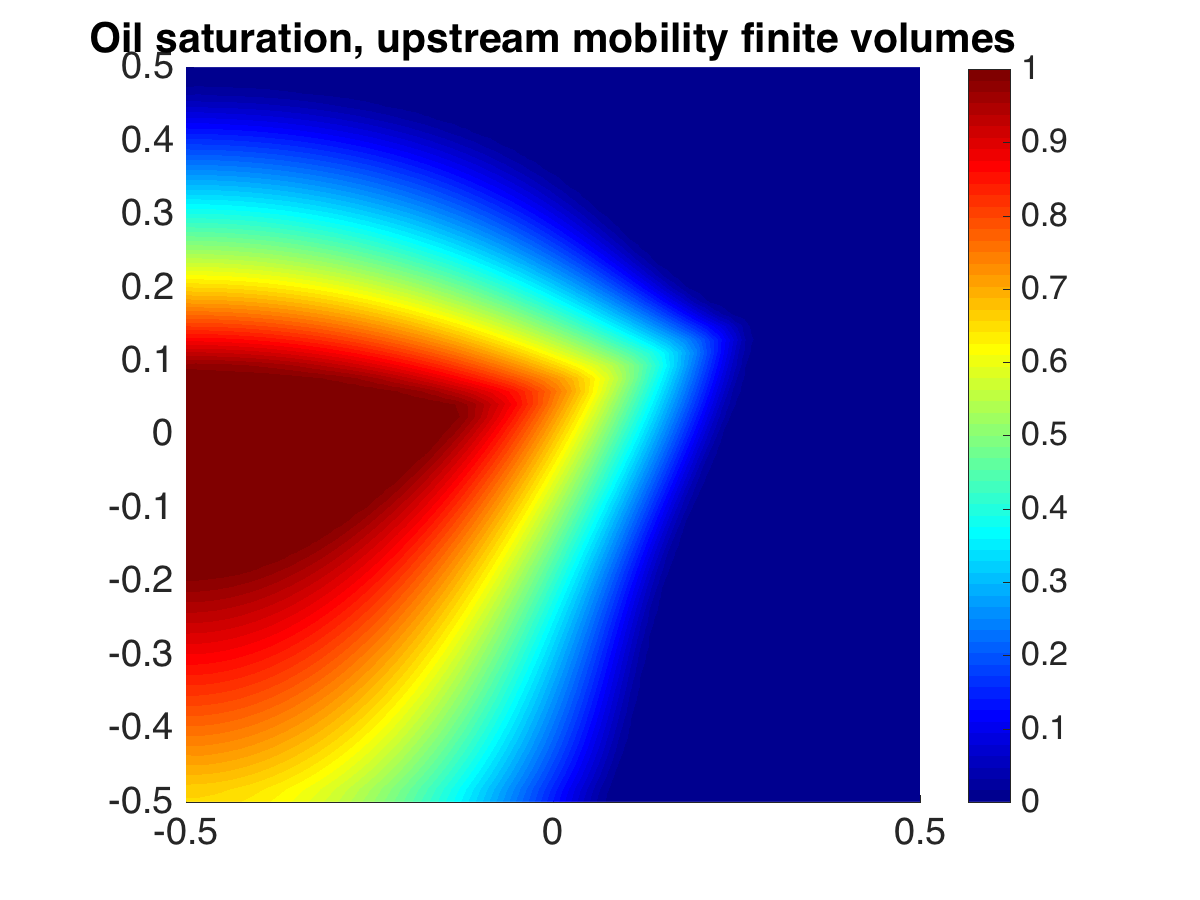} 
\caption{t=2.5}
\end{subfigure}
\medskip

\begin{subfigure}{\textwidth}
\centering
\includegraphics[width=5cm]{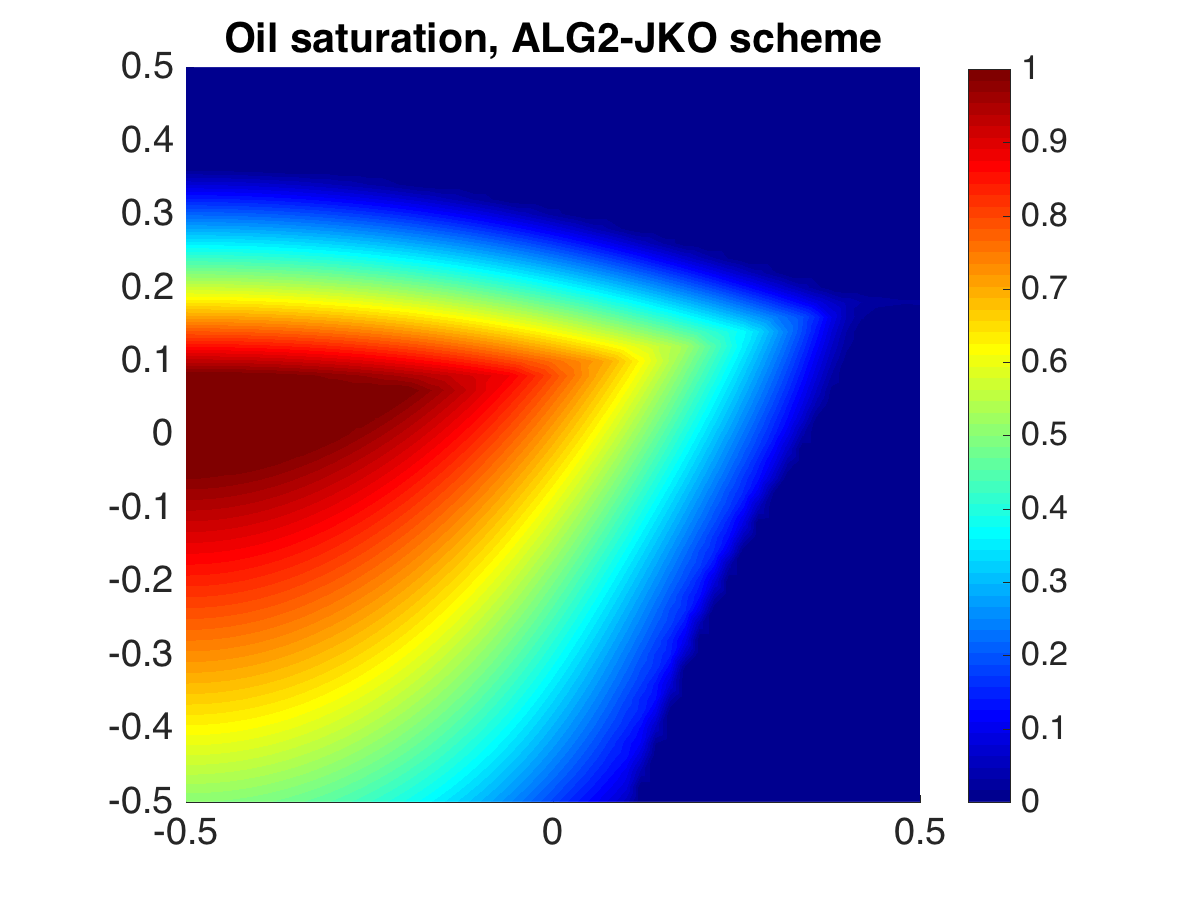} 
\hspace{.1cm}
\includegraphics[width=5cm]{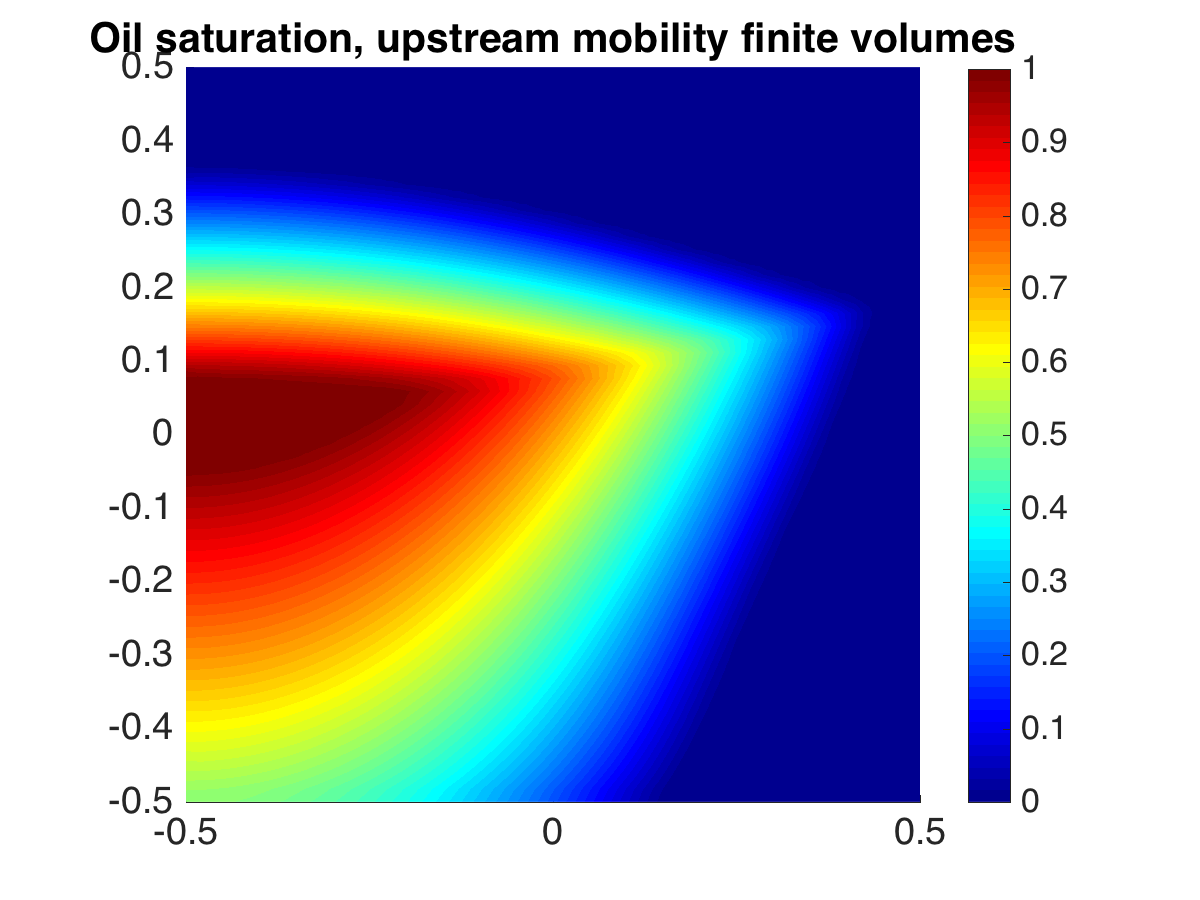} 
\caption{t=5}
\end{subfigure}
\medskip

\begin{subfigure}{\textwidth}
\centering
\includegraphics[width=5cm]{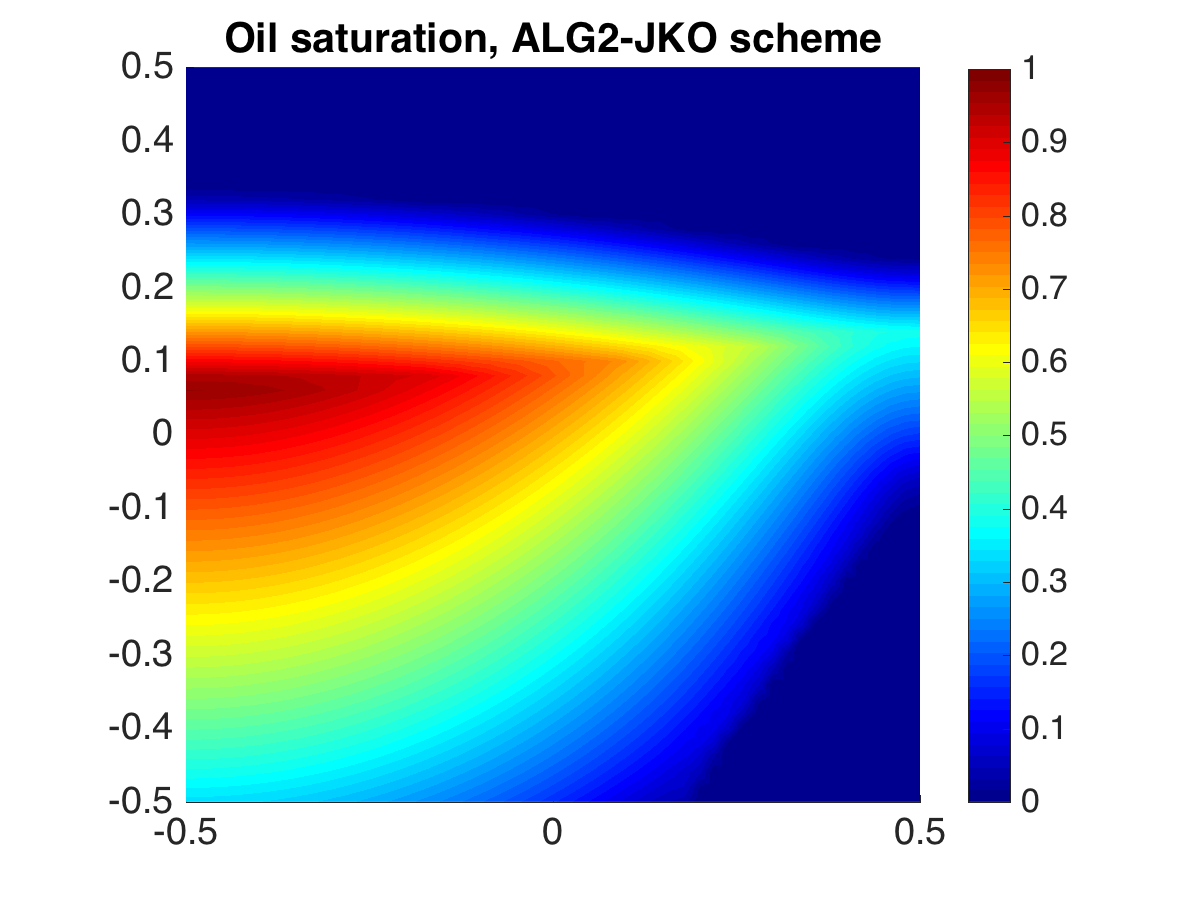} 
\hspace{.1cm}
\includegraphics[width=5cm]{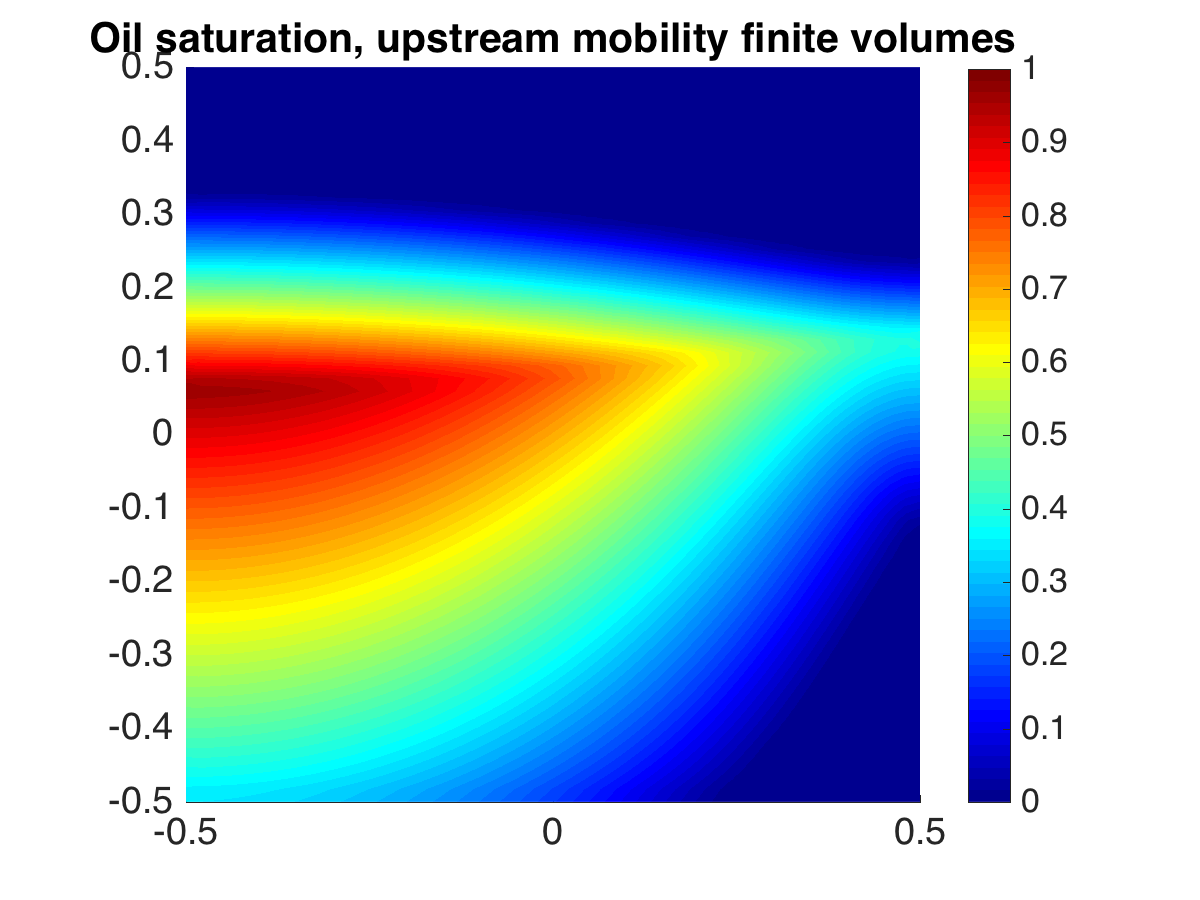} 
\caption{t=10}
\end{subfigure}
\caption{Three-phase flow, snapshots of the oil saturation profiles at different times: ALG2-JKO scheme (left) and upstream mobility finites volumes (right).}
\label{fig:3phases_o}
\end{figure}

\begin{figure}[!htbp]
\centering
\begin{subfigure}{\textwidth}
\centering
\includegraphics[width=5cm]{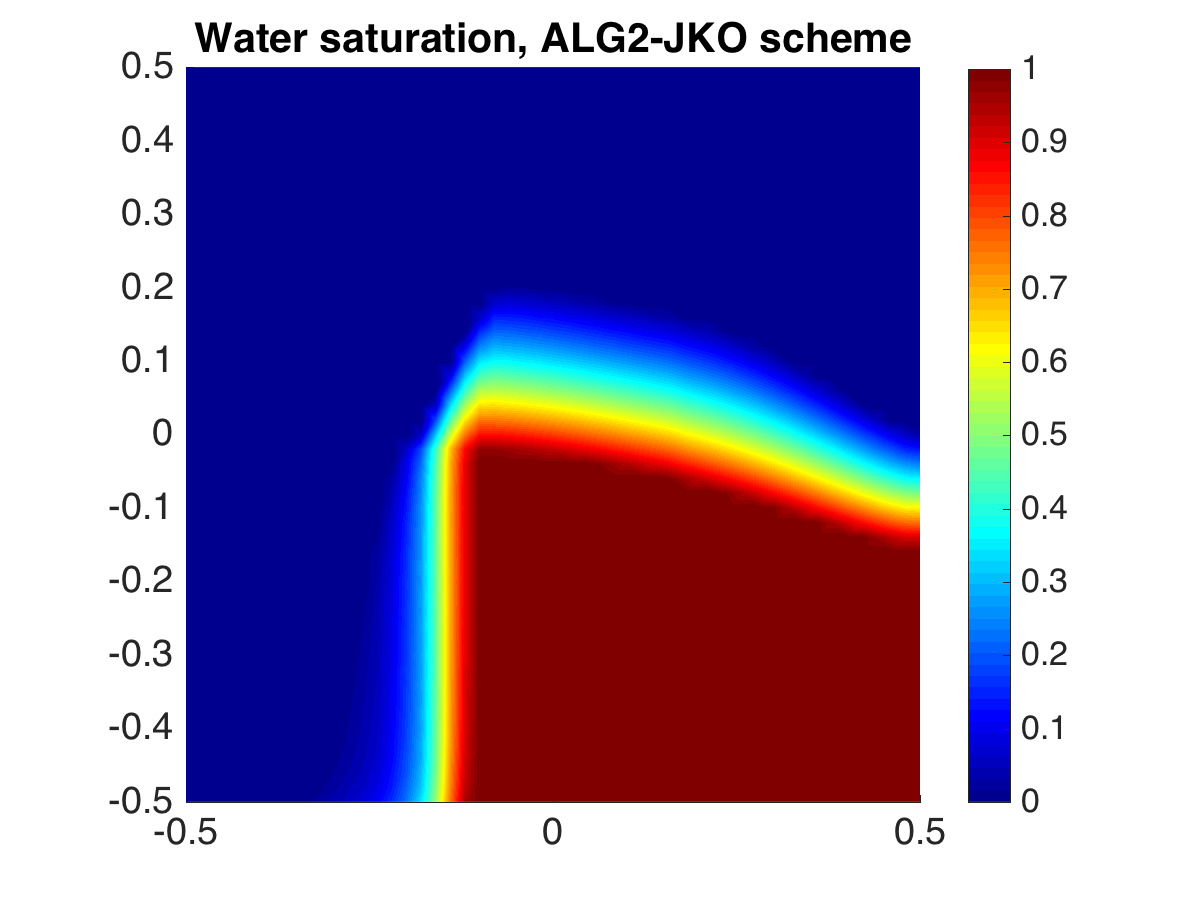} 
\hspace{.1cm}
\includegraphics[width=5cm]{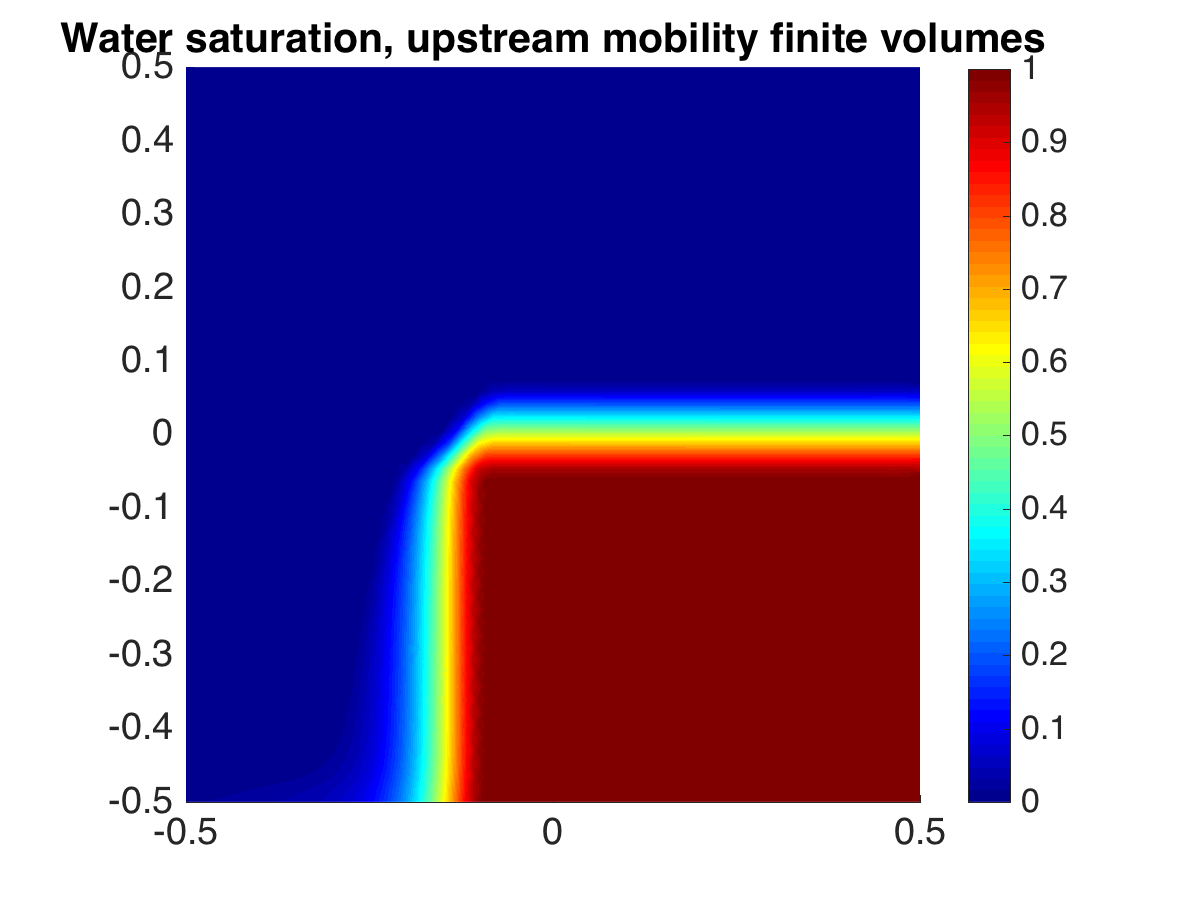} 
\caption{t=0.1}
\end{subfigure}

\medskip

\begin{subfigure}{\textwidth}
\centering
\includegraphics[width=5cm]{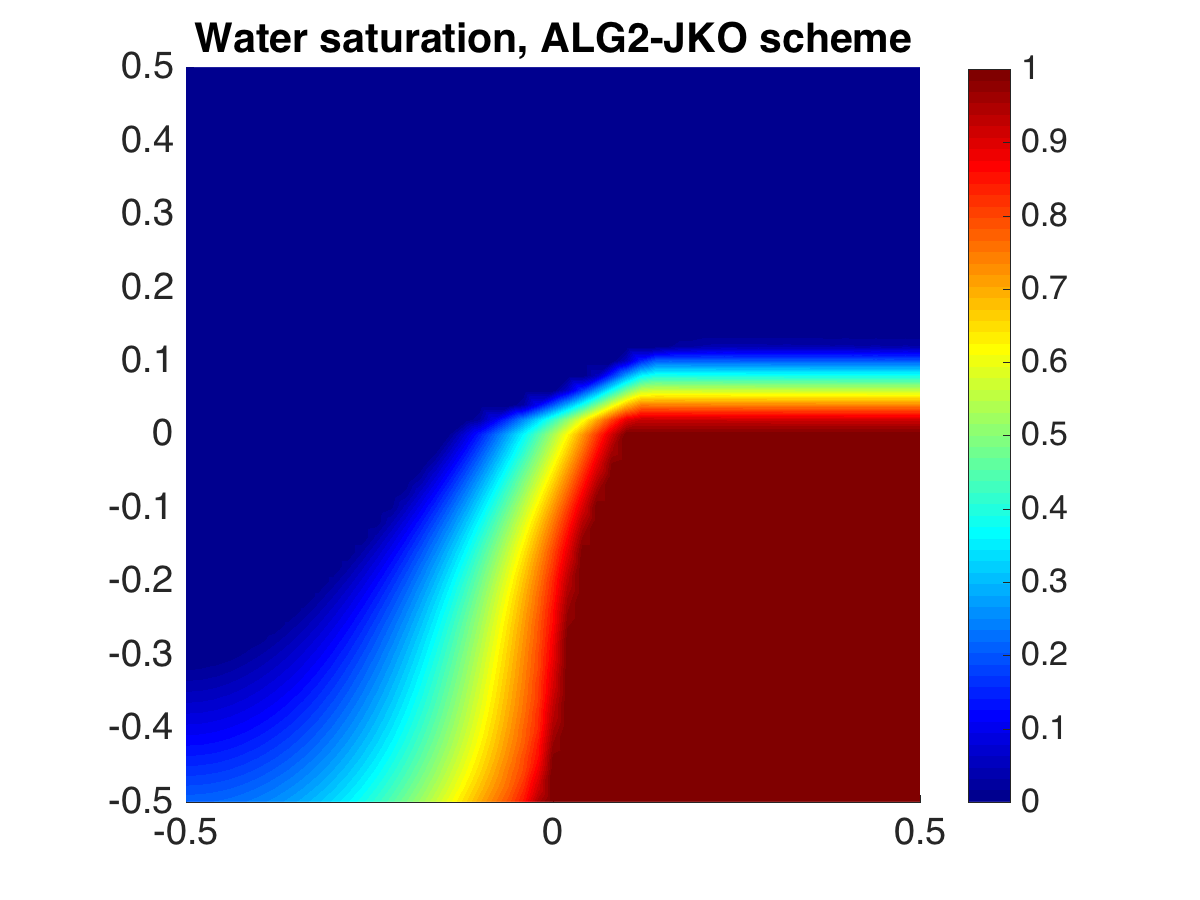} 
\hspace{.1cm}
\includegraphics[width=5cm]{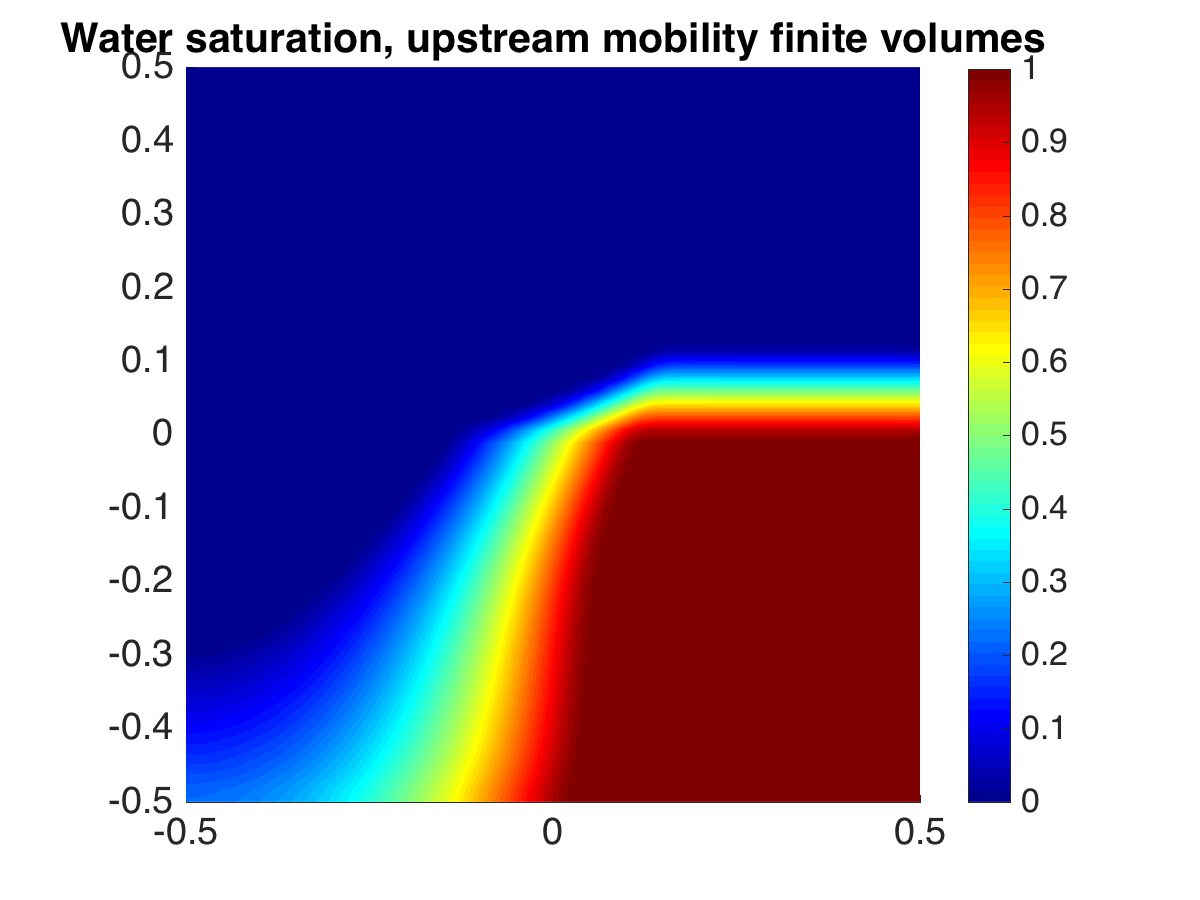} 
\caption{t=1.25}
\end{subfigure}

\medskip

\begin{subfigure}{\textwidth}
\centering
\includegraphics[width=5cm]{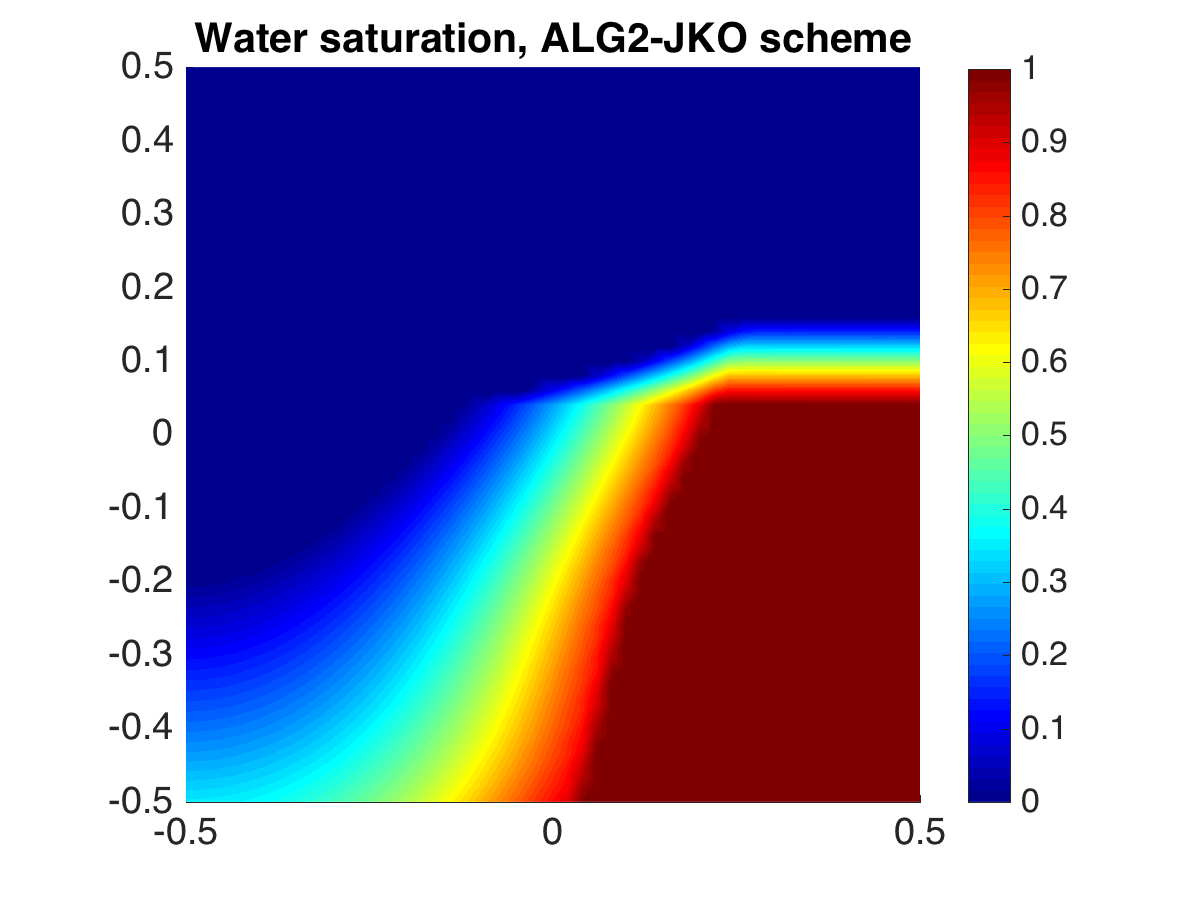} 
\hspace{.1cm}
\includegraphics[width=5cm]{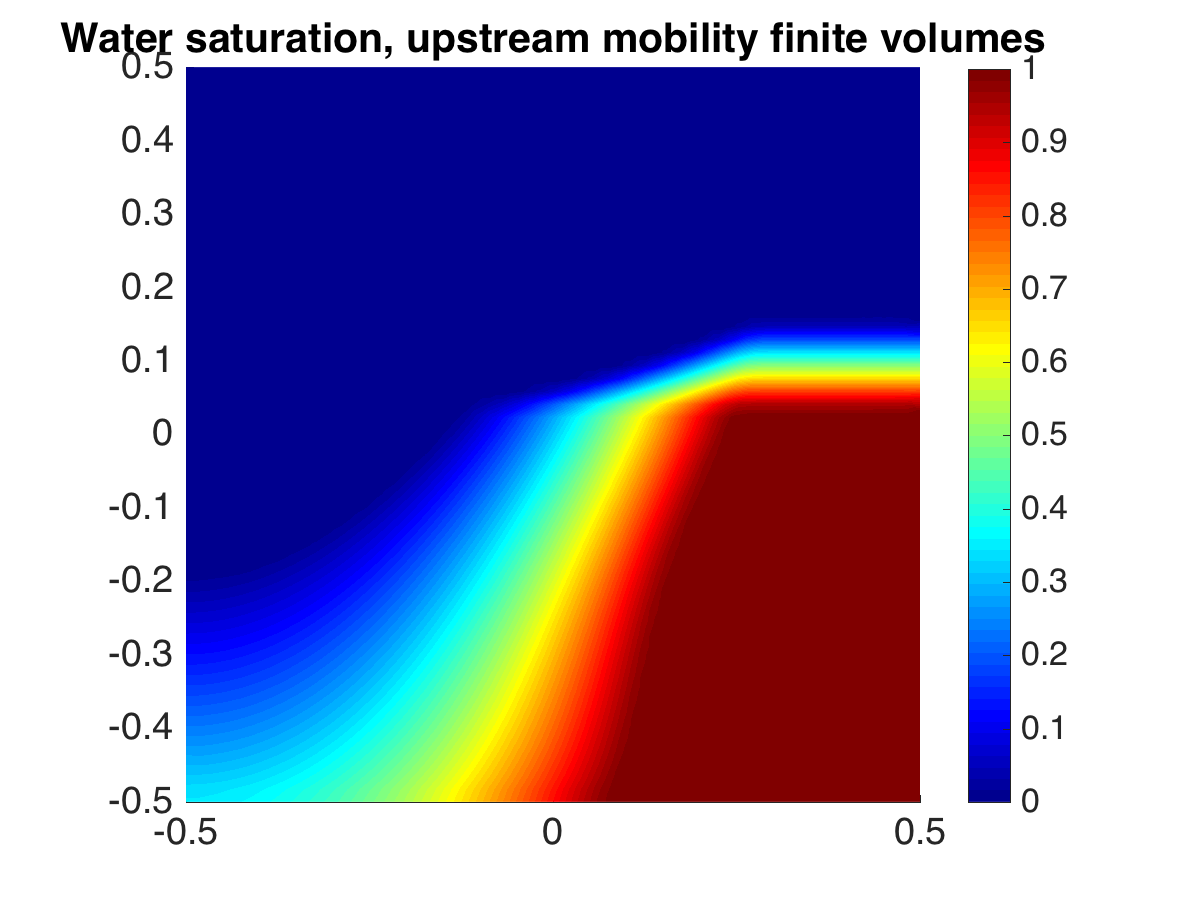} 
\caption{t=2.5}
\end{subfigure}

\medskip

\begin{subfigure}{\textwidth}
\centering
\includegraphics[width=5cm]{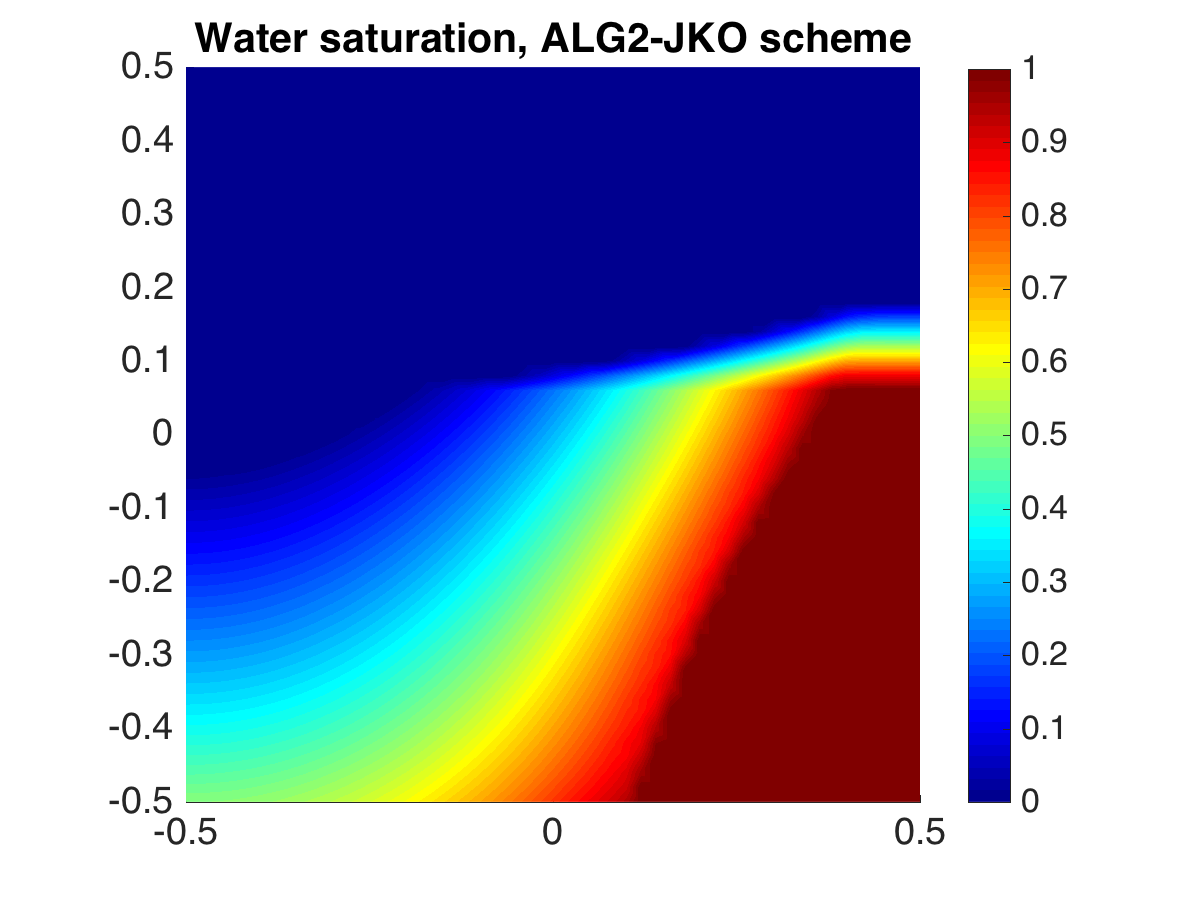} 
\hspace{.1cm}
\includegraphics[width=5cm]{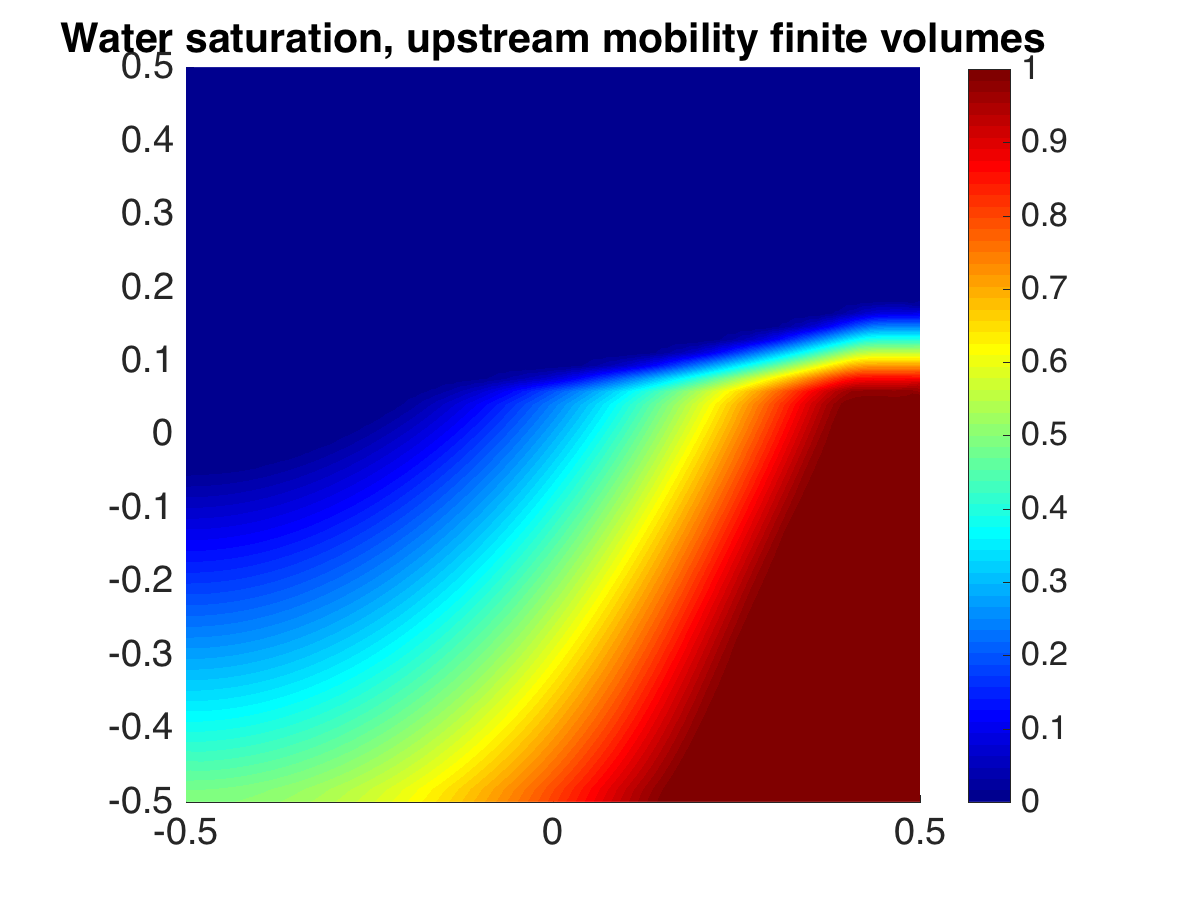} 
\caption{t=5}
\end{subfigure}

\medskip

\begin{subfigure}{\textwidth}
\centering
\includegraphics[width=5cm]{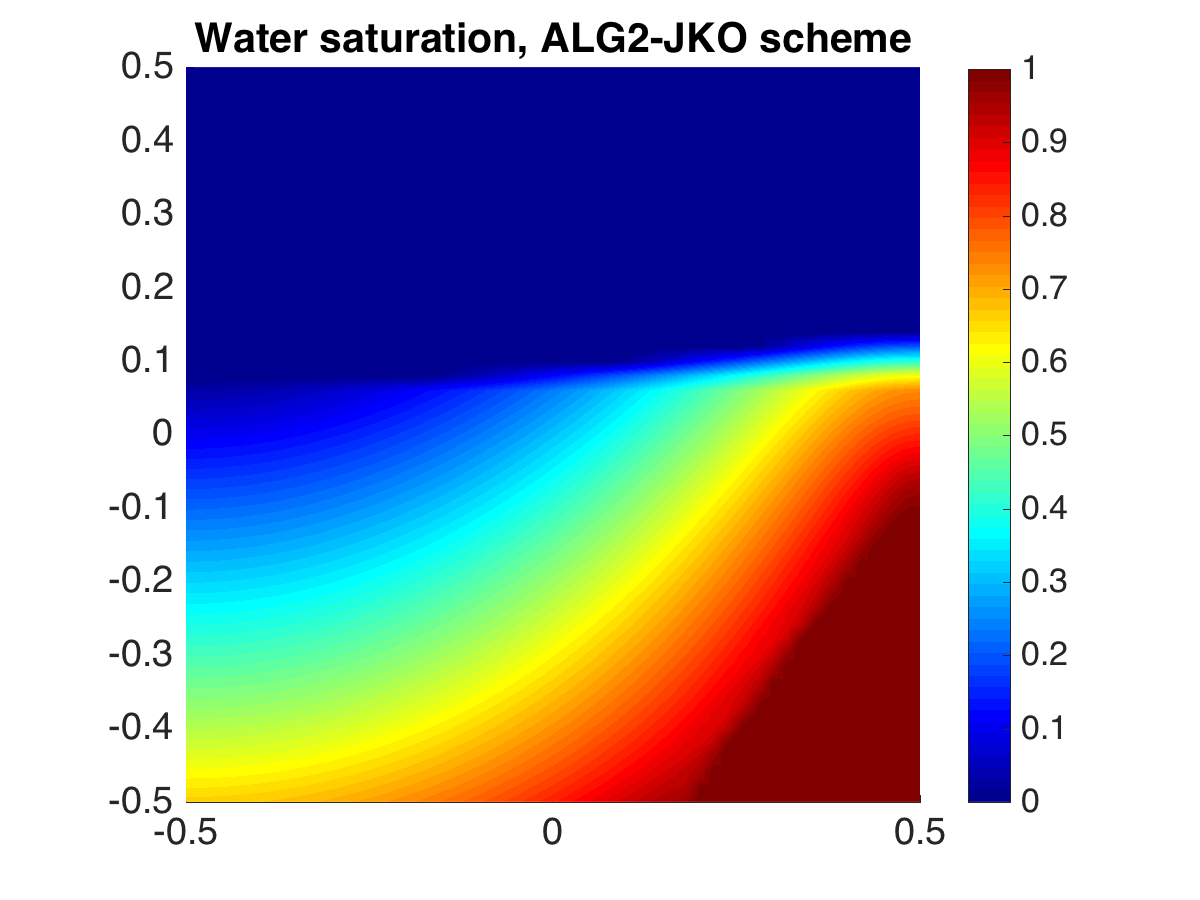} 
\hspace{.1cm}
\includegraphics[width=5cm]{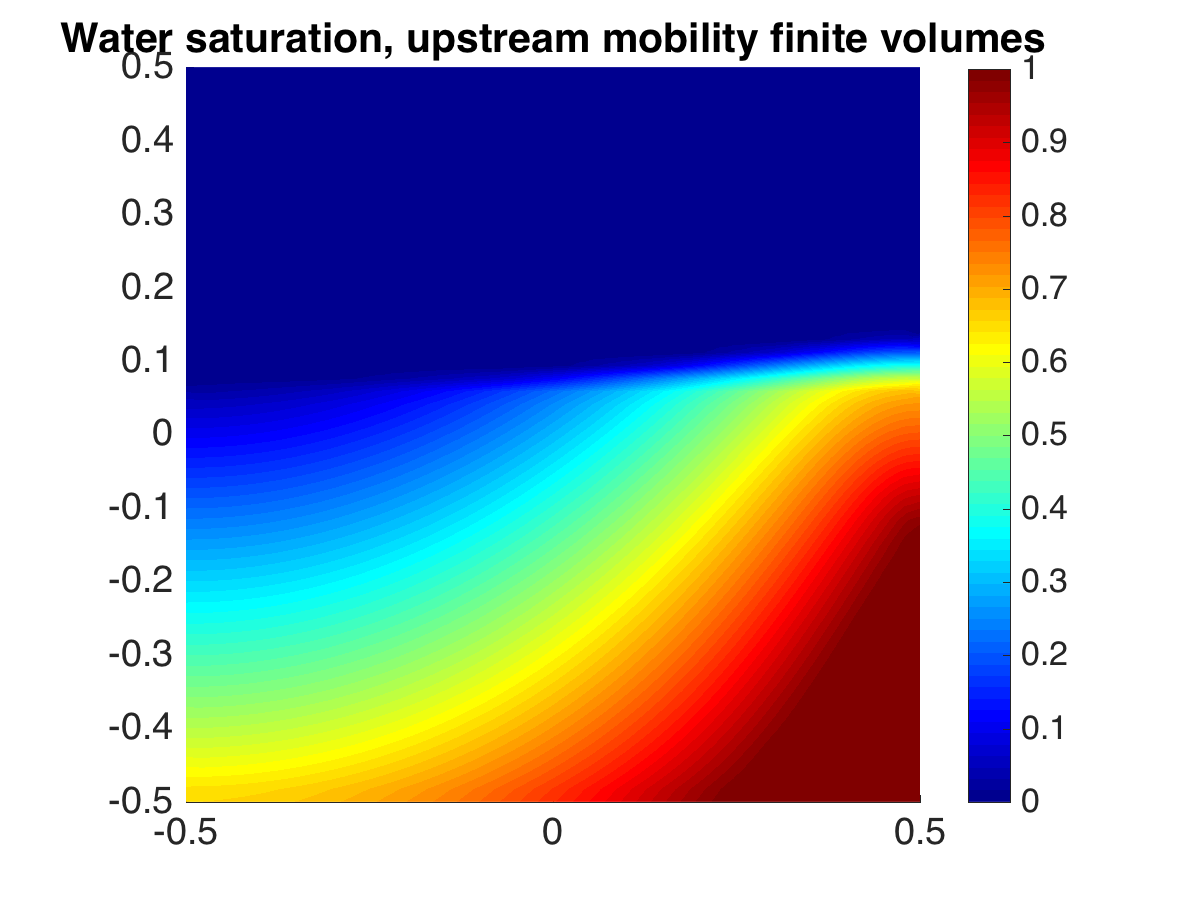} 
\caption{t=10}
\end{subfigure}
\caption{Three-phase flow, snapshots of the water saturation profiles at different times: ALG2-JKO scheme (left) and upstream mobility finites volumes (right).}
\label{fig:3phases_w}
\end{figure}

\begin{figure}[!htbp]
\centering
\begin{subfigure}{\textwidth}
\centering
\includegraphics[width=5cm]{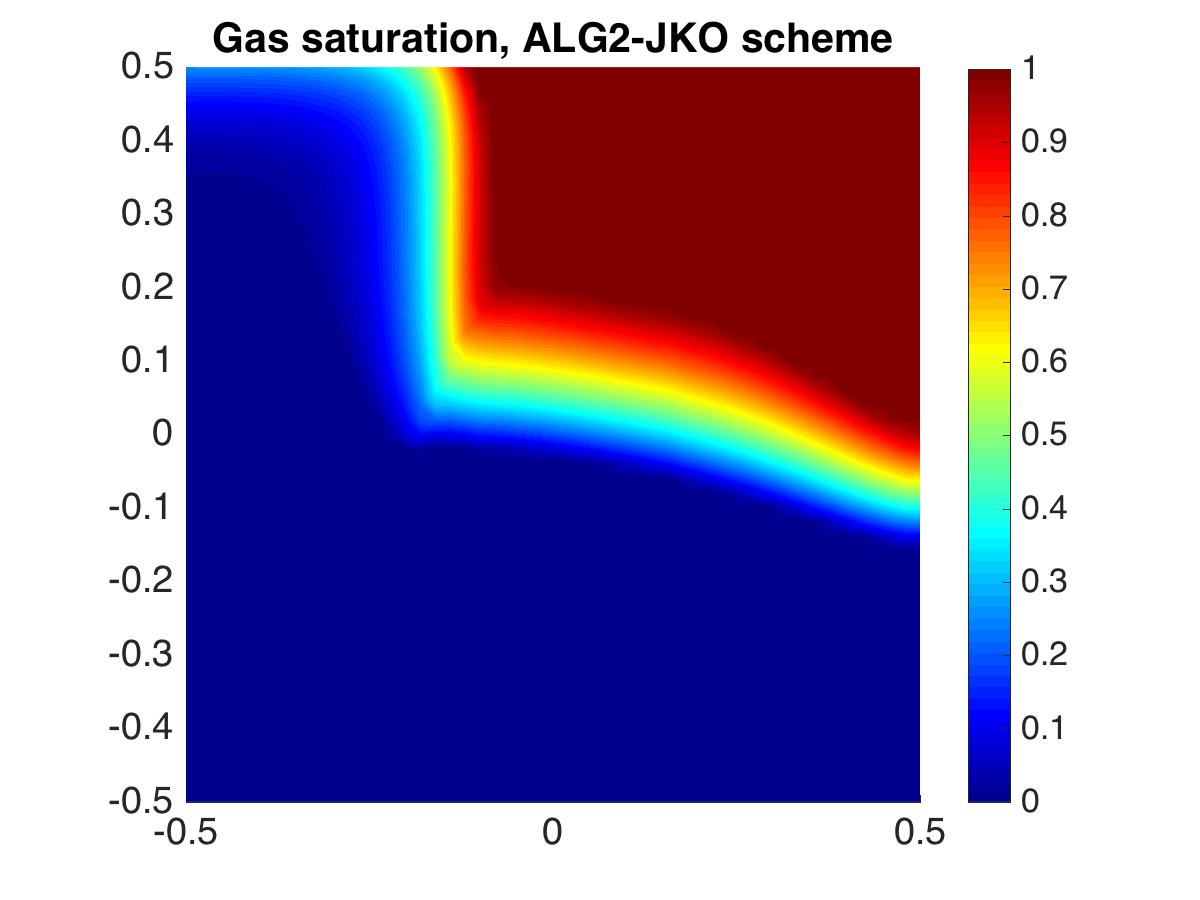} 
\hspace{.1cm}
\includegraphics[width=5cm]{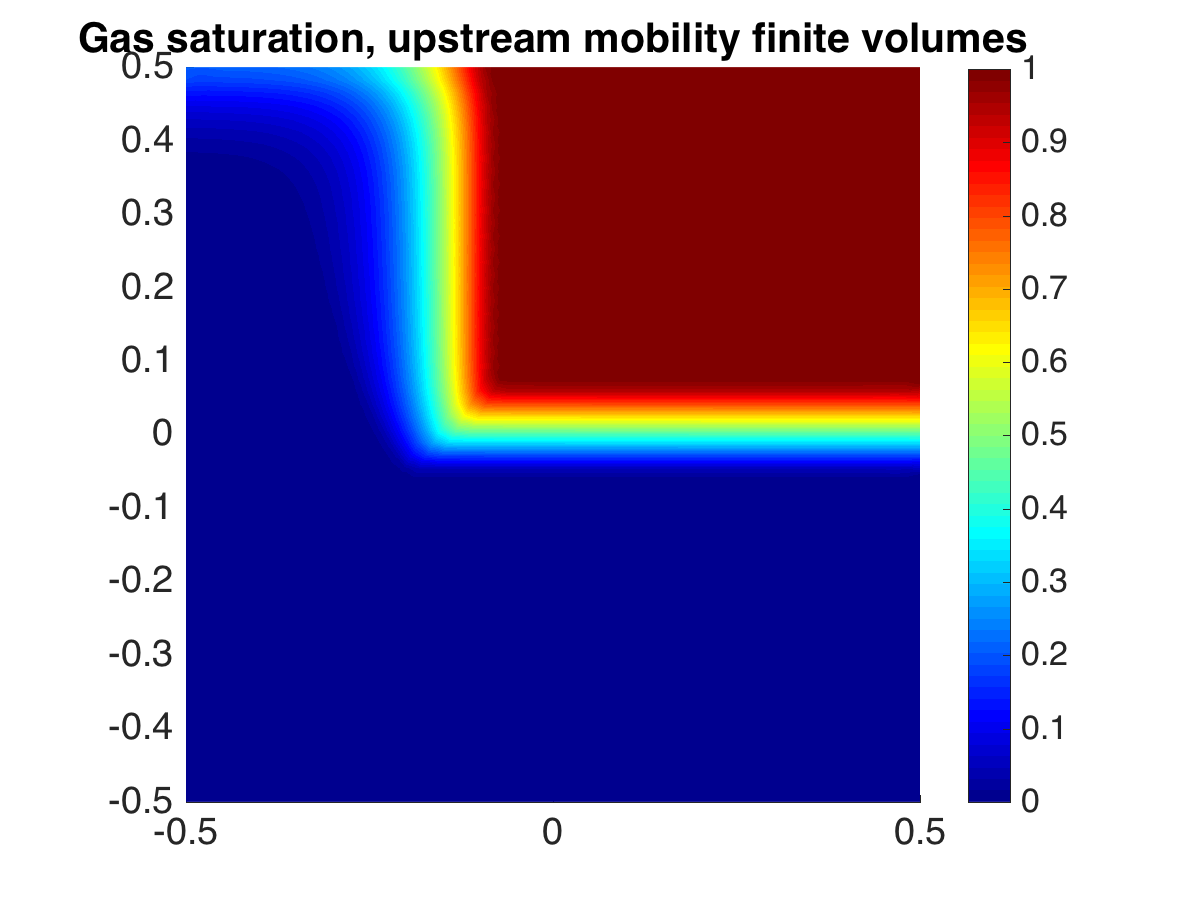} 
\caption{t=0.1}
\end{subfigure}

\medskip

\begin{subfigure}{\textwidth}
\centering
\includegraphics[width=5cm]{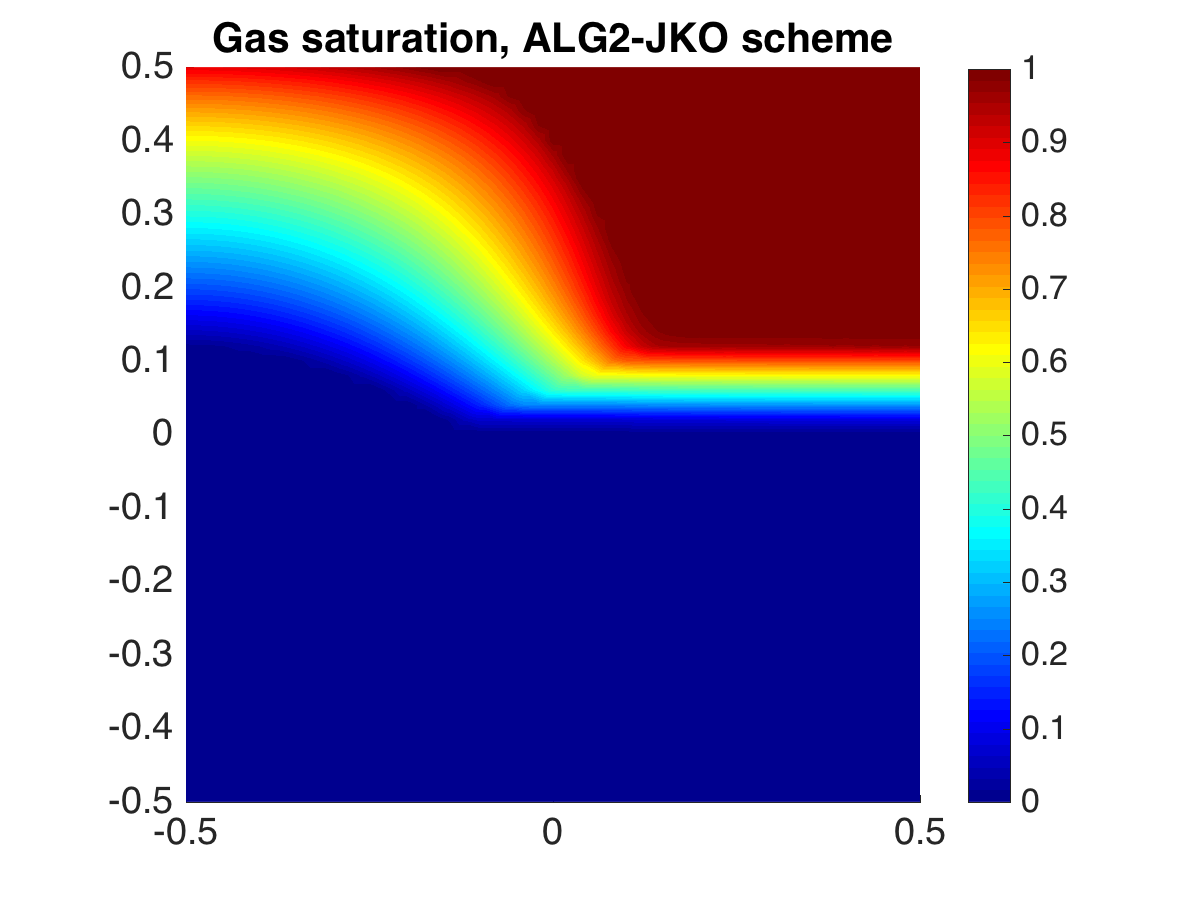} 
\hspace{.1cm}
\includegraphics[width=5cm]{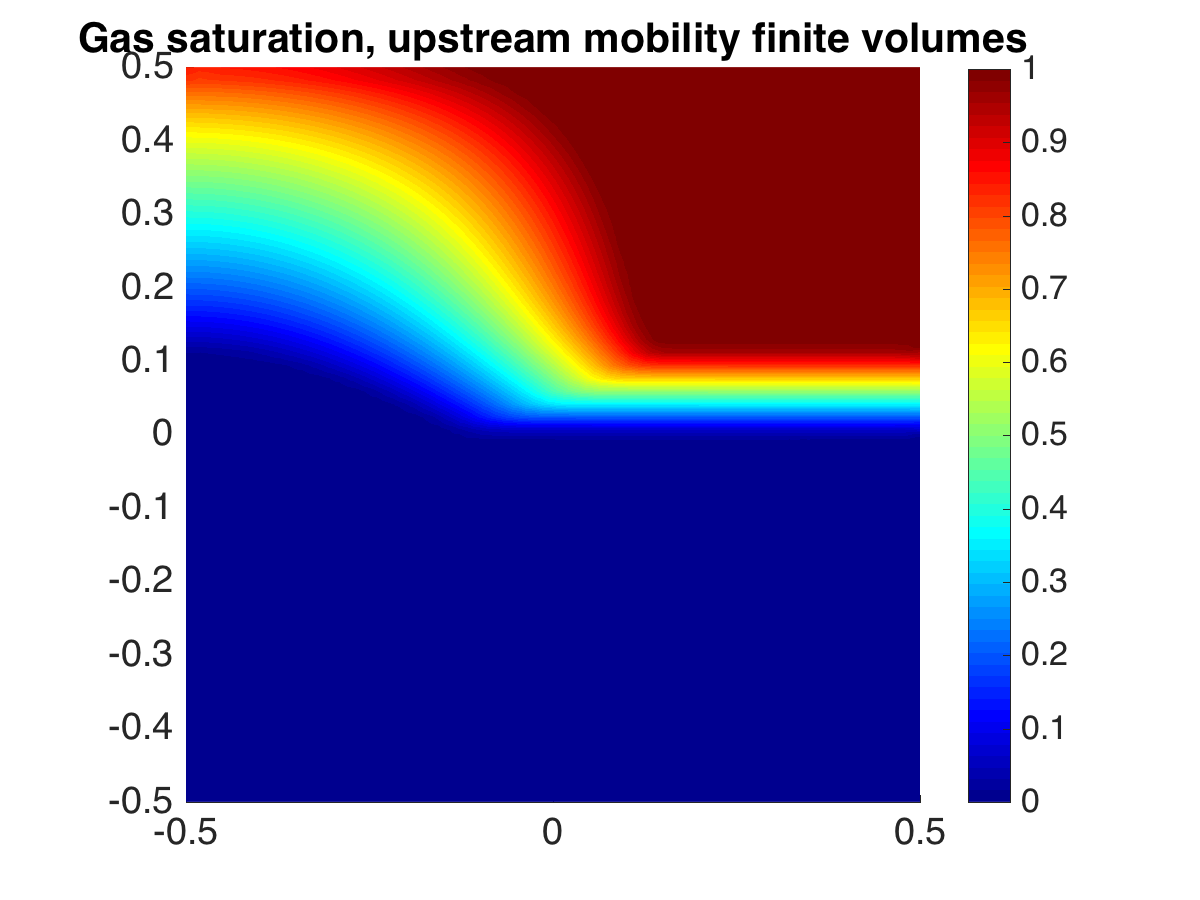} 
\caption{t=1.25}
\end{subfigure}

\medskip

\begin{subfigure}{\textwidth}
\centering
\includegraphics[width=5cm]{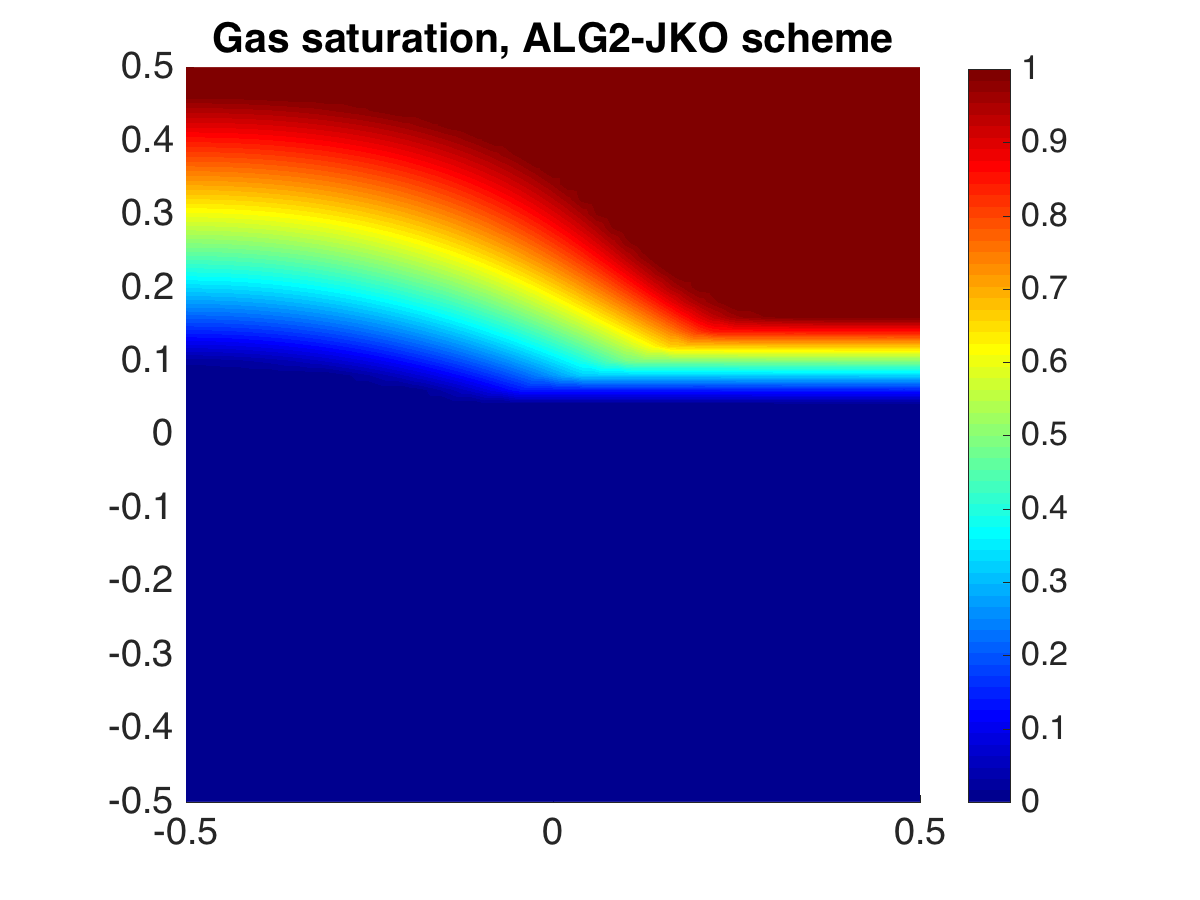} 
\hspace{.1cm}
\includegraphics[width=5cm]{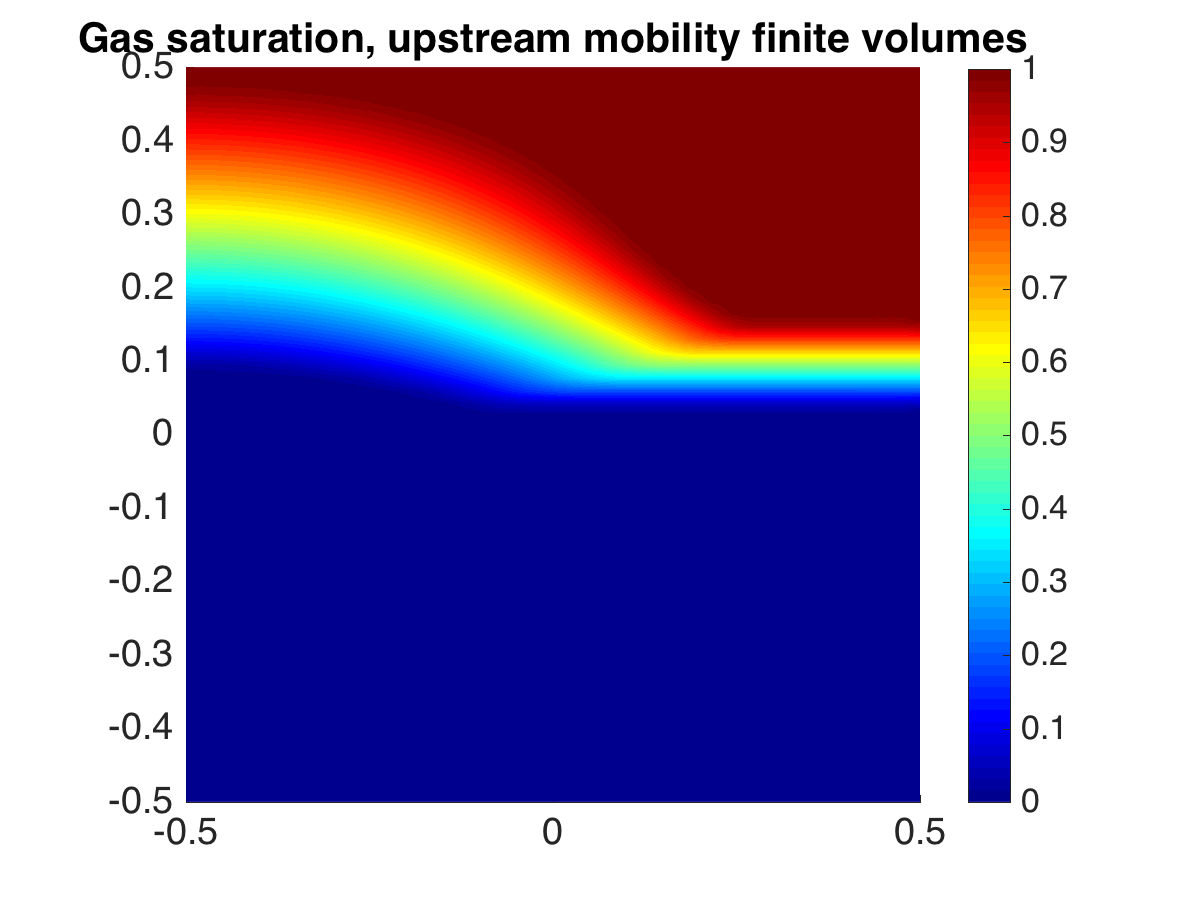} 
\caption{t=2.5}
\end{subfigure}

\medskip

\begin{subfigure}{\textwidth}
\centering
\includegraphics[width=5cm]{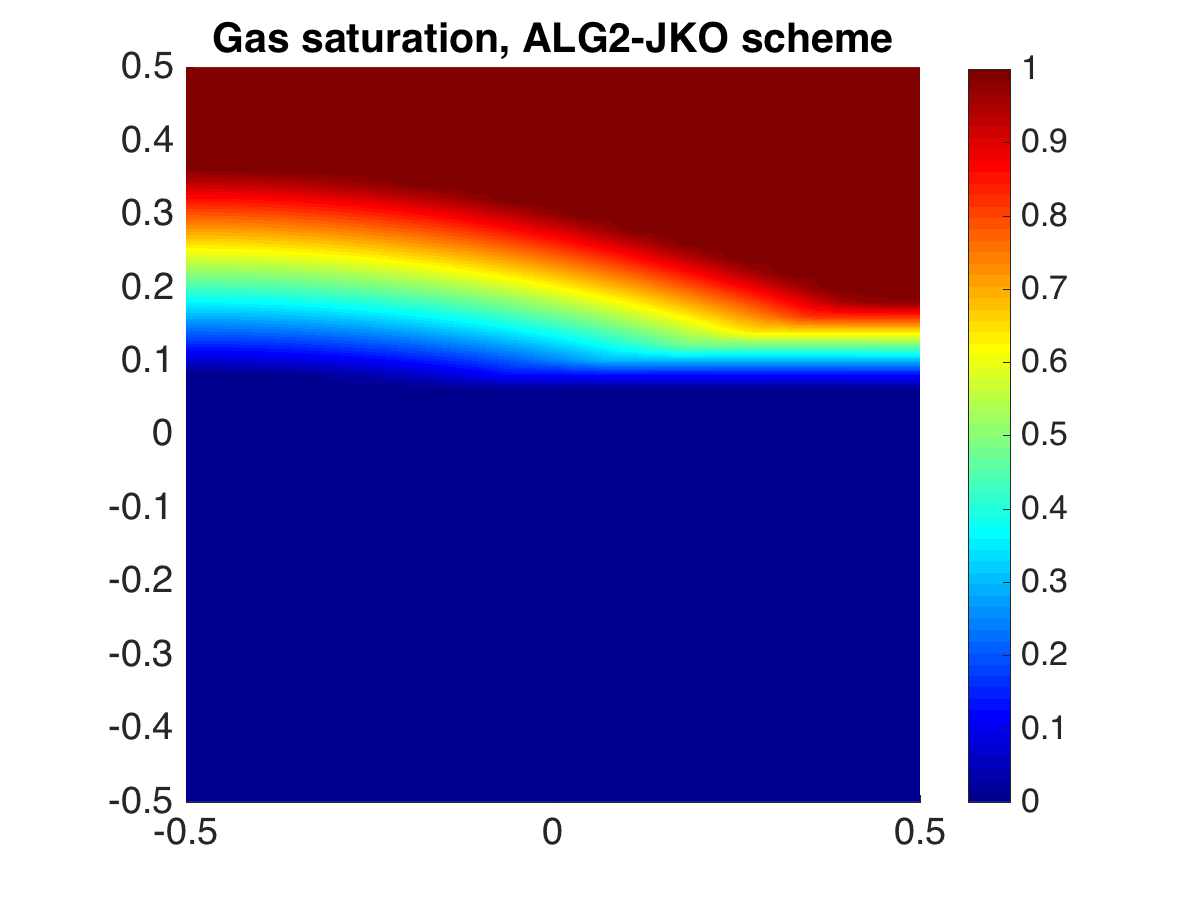} 
\hspace{.1cm}
\includegraphics[width=5cm]{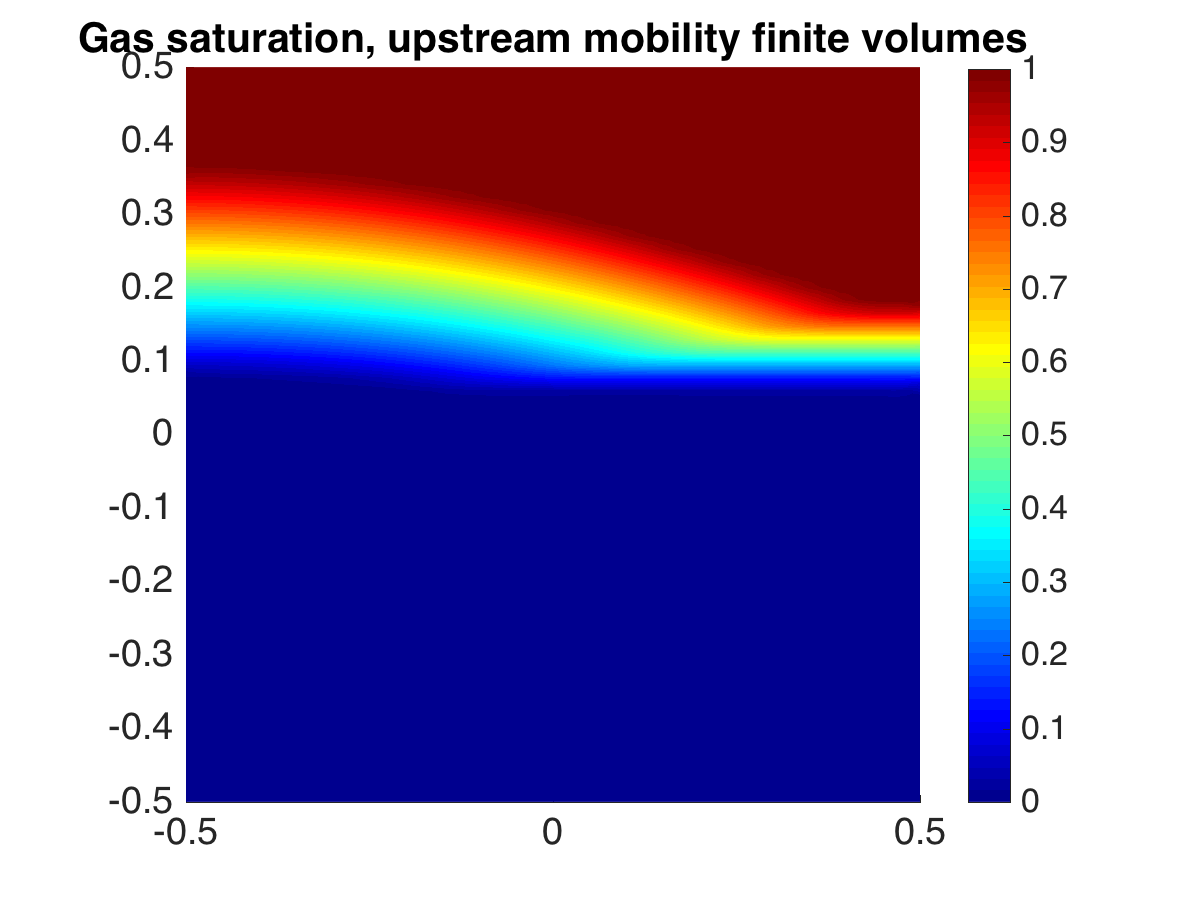} 
\caption{t=5}
\end{subfigure}

\medskip

\begin{subfigure}{\textwidth}
\centering
\includegraphics[width=5cm]{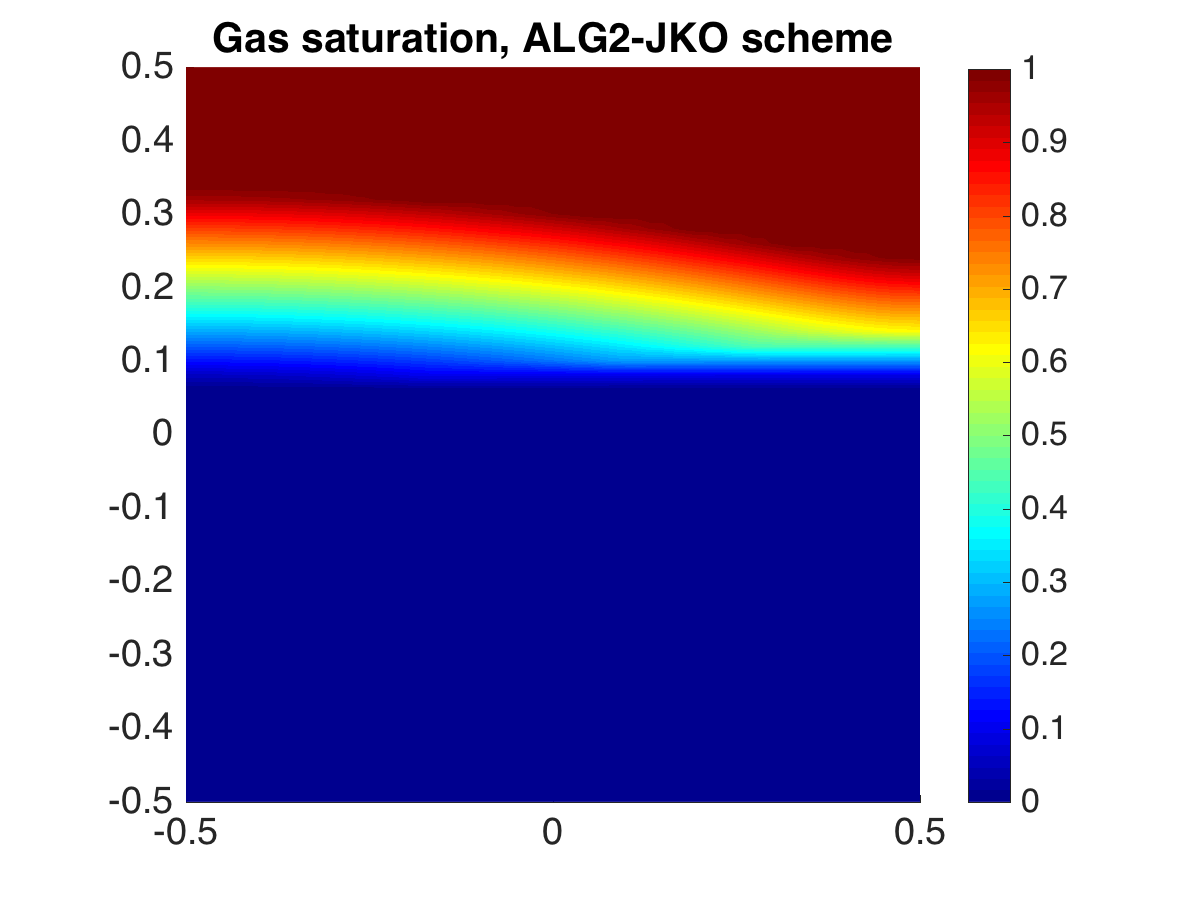} 
\hspace{.1cm}
\includegraphics[width=5cm]{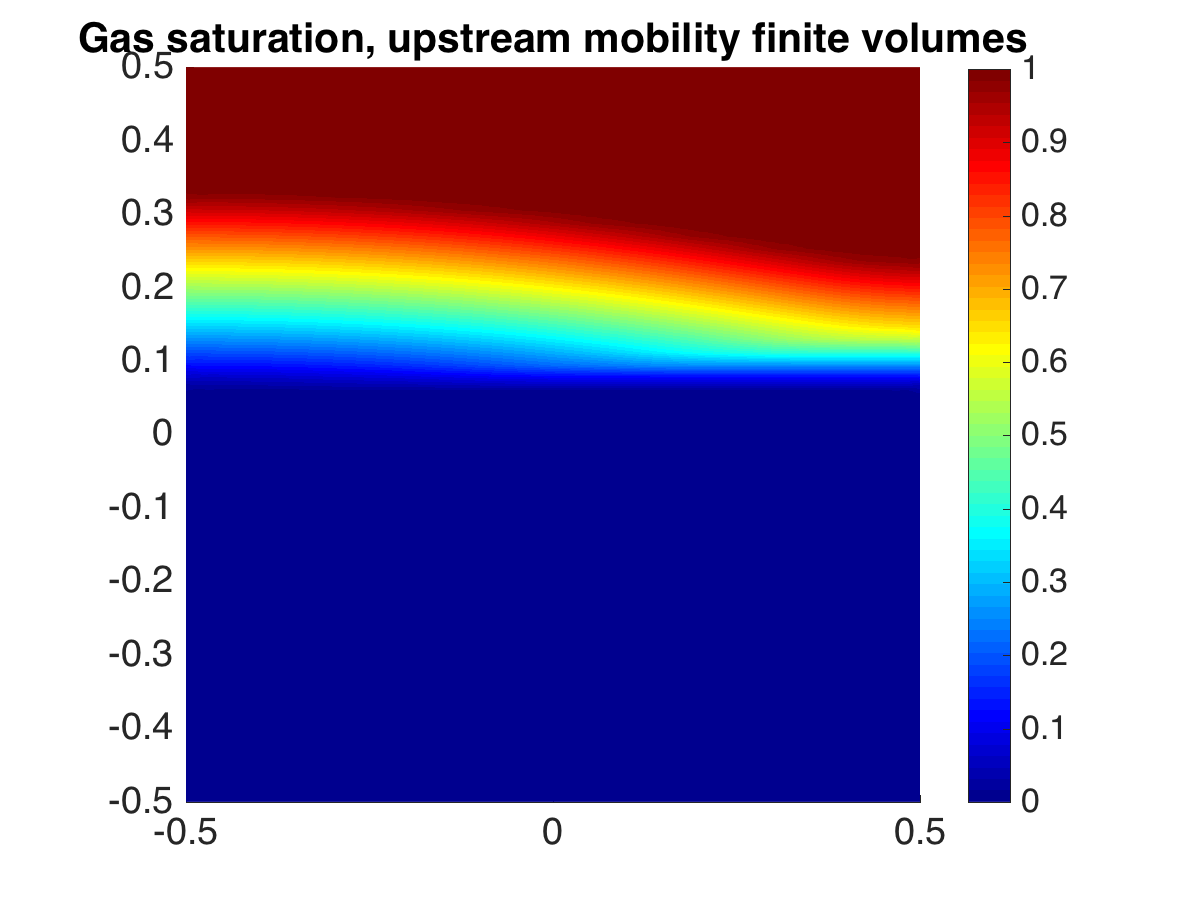} 
\caption{t=10}
\end{subfigure}
\caption{Three-phase flow, snapshots of the gas saturation profiles at different times: ALG2-JKO scheme (left) and upstream mobility finites volumes (right).}
\label{fig:3phases_g}
\end{figure}

Figures \ref{fig:3phases_o}--\ref{fig:3phases_g} show the evolution of the three phases with quadratic capillarity potential.
Again, the simulation with the ALG2-JKO scheme is carried out using a $50 \times 50$ discretization in space, with a single inner time step.
There are $200$ JKO steps ($\tau=0.05$). The convergence of the augmented Lagrangian 
iterative method is rather slow: it took around 10 hours on a laptop to produce the results with {\tt FreeFem++}. 
But because of the large viscosity ratio, Newton's method had severe difficulties to converge for the upstream mobility scheme. 
A very small time step ($\tau=10^{-4}$) was needed, so that more that 2 days of computation on a cluster were needed 
to produce the results with {\tt Matlab}. \clem{Concerning the upstream mobility finite volume scheme, we run the scheme on an unstructured Delaunday triangulation 
made of $5645$ triangles.}
Once again, both methods produce similar results, as highlighted on the figures~\ref{fig:3phases_o}--\ref{fig:3phases_g} below.
\begin{figure}[!htbp]
\centering
\includegraphics[width=0.31\textwidth]{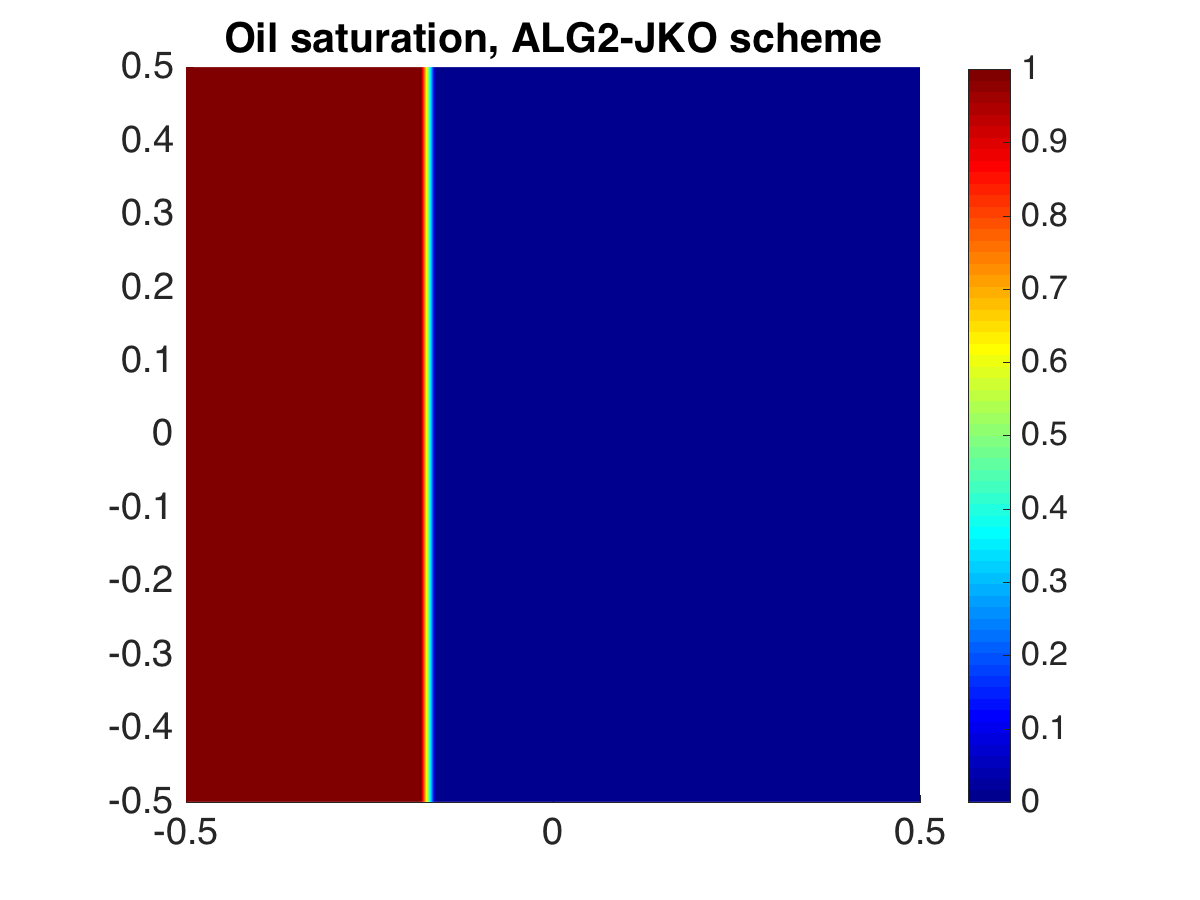} 
\includegraphics[width=0.31\textwidth]{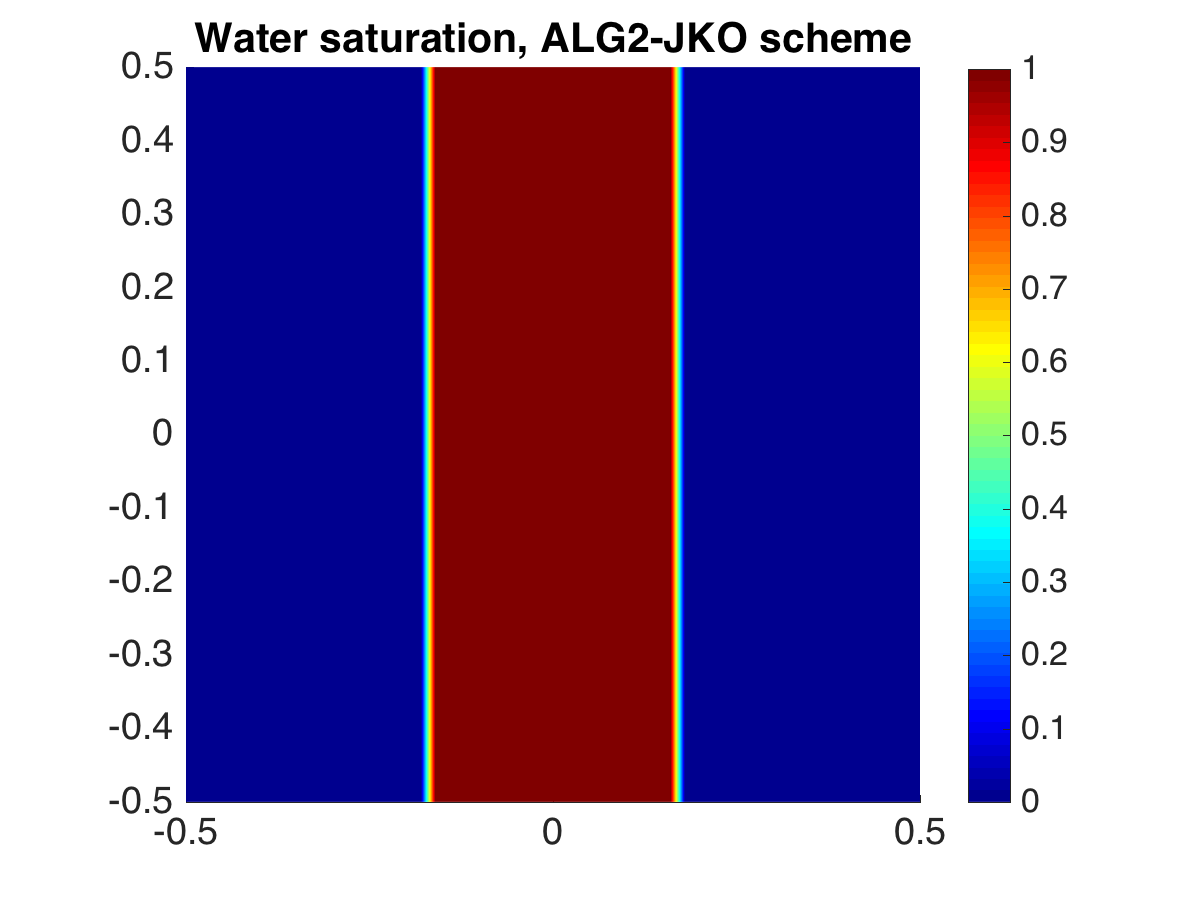} 
\includegraphics[width=0.31\textwidth]{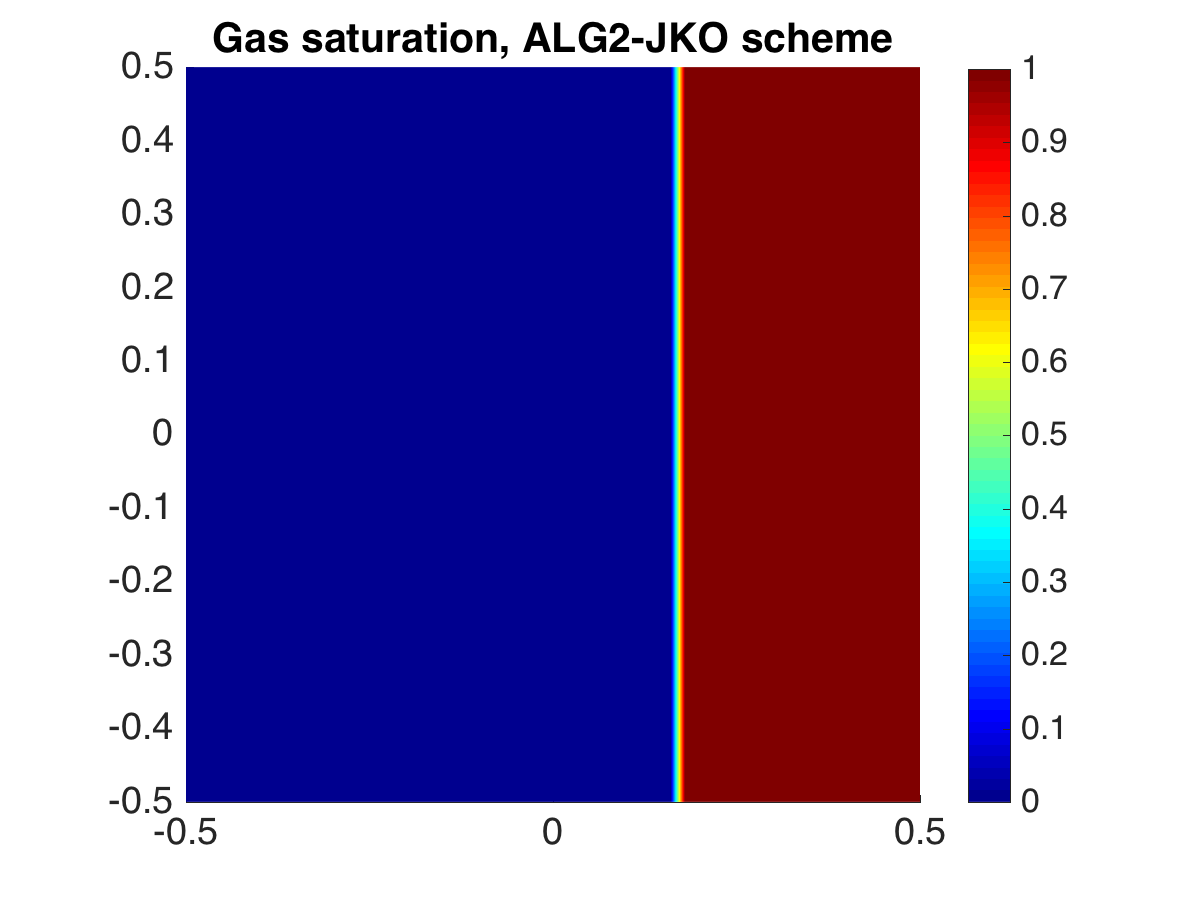} 
\caption{Initial oil (left), water (center) and gas (right) saturation profiles.}
\label{fig:3phases_ini}
\end{figure}

Due to the large viscosity ratios, two distinct time scale appear in the numerical results. 
Since water and gas have smaller mobilities, they move much faster than oil. 
\clem{This quick phenomenon 
is not well captured by the ALG2-JKO scheme. 
The interface between oil and gas is already almost horizontal at $t=0.1$. 
This horizontal interface is captured by the finite volume scheme but 
not by the ALG2-JKO scheme that encounters difficulties to converge for the early time steps. The finite volume scheme 
also has difficulties to converge, enforcing us to consider very small time steps}. 
Oil is much less mobile and its interface with the two other 
phases remains almost vertical at that time. Then oil evolves slowly towards its equilibrium state, 
that consists in a horizontal layer trapped between gas above and water below. This long time 
equilibrium is not yet reached for $t=10$.

\clem{
\subsection{Energy dissipation}
As already highlighted, both schemes dissipate the energy along time. 
The goal of this test case is to compare the energy dissipation. To this end, 
we consider a test case proposed in~\cite{Can_OGST}. We consider a two-phase flow 
with oil ($i=1$) and water ($i=0$) with $\rho_1 = 0.87$, $\rho_0 = 1$, $\mu_1 = 10$ and $\mu_0 =1$, 
while $\k = 1$ and $\omega = 1$.
The capillary pressure law is given by 
\[
p_1 - p_0 = \pi_1(s_1) =  \frac{s_1}{2}, 
\]
so that the energy is defined by 
\[
\Ee(s_1) = \int_\O \left(\frac{(s_1)^2}{4} + s_1 (\rho_0-\rho_1) \g \cdot \x \right).
\]
We consider the initial data $s_1^{0}(\x)= e^{-4|\x|^2}$. At equilibrium, the saturation $s_1^\infty$ minimizes 
$\Ee$ under the constraints $s_1^\infty \in [0,1]$ and 
\be\label{eq:mass_eq}
\int_\O s_1^\infty = \int_\O s_1^{0}.
\ee It is therefore given by
\be\label{eq:steady}
\text{either}\; s_1^\infty\in \{0,1\} \; \text{or}\; \pi_1(s_1) = (\rho_1 - \rho_0)\g\cdot \x + \gamma, 
\ee
the constant $\gamma$ being fixed thanks to~\eqref{eq:mass_eq}. Similar calculations can be performed in the discrete settings, 
both for the ALG2-JKO scheme and the finite volume scheme. 
\begin{figure}[htb]
\centering
\includegraphics[width = 5cm]{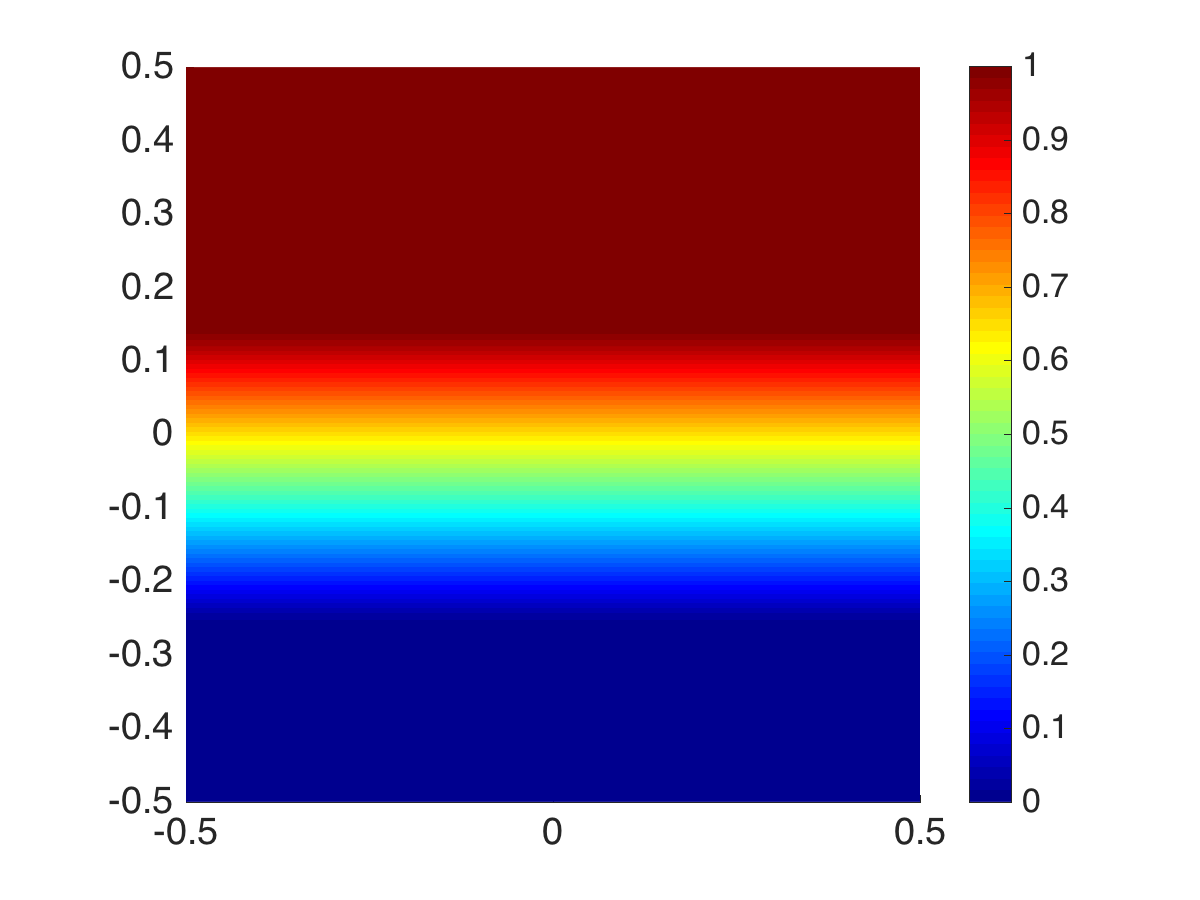}
\hspace{.5cm}
\includegraphics[width = 5cm]{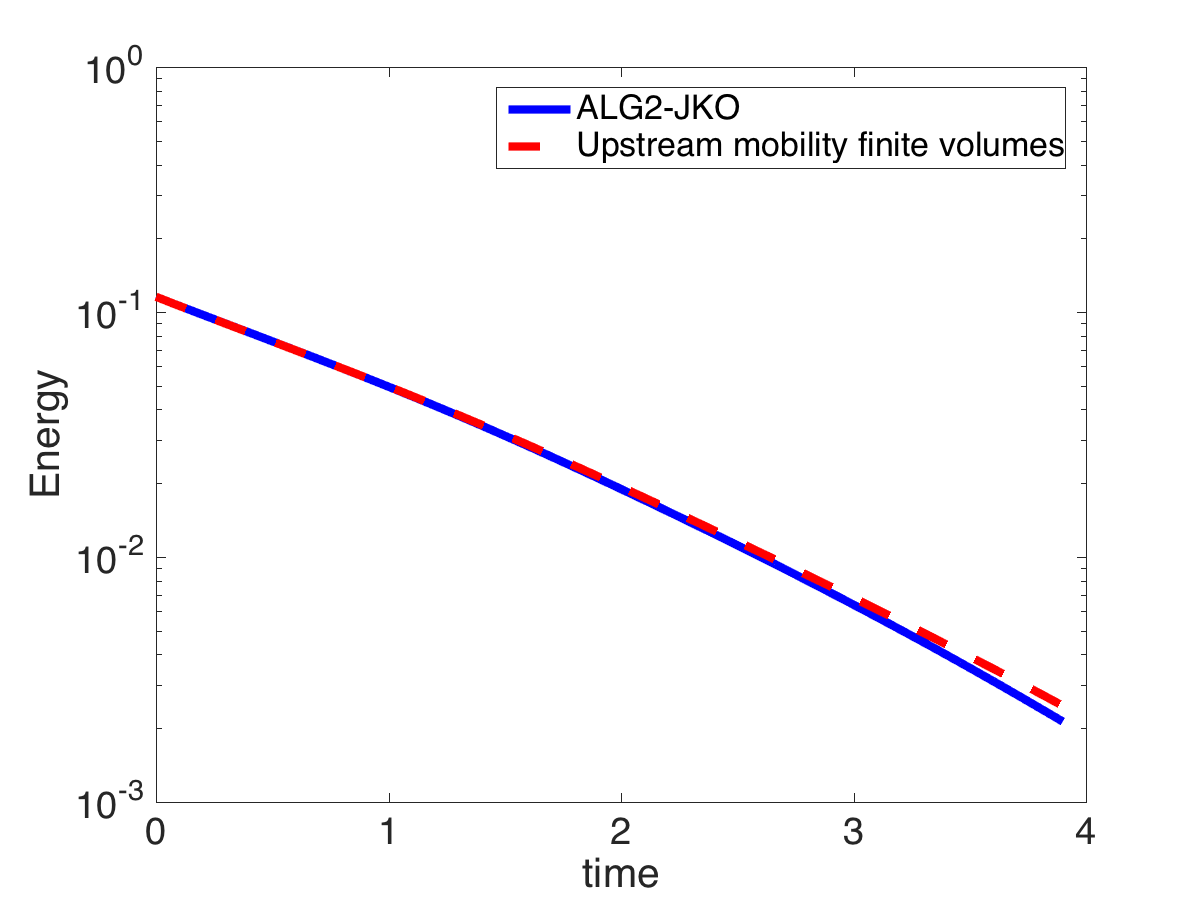}
\caption{\clem{Left: The steady state~\eqref{eq:steady}. Right: The relative energies computed with the ALG2-JKO scheme (blue) and the finite volume scheme (red).}}
\label{fig:entrel}
\end{figure}
Then one computes for both scheme the relative energy $\Ee(s_1) - \Ee(s_1^\infty)\geq 0$, 
that we plot as a function of time on Figure~\ref{fig:entrel}. The convergence towards the equilibrium appears to be exponential 
in both cases. 
}


\section{Conclusion}\label{sec:conclu}

We proposed to apply the ALG2-JKO scheme of~\cite{BCL15} to simulate multiphase porous media 
flows. The results have been compared to the widely used upstream mobility finite volume scheme. 
The ALG2-JKO scheme appears to be robust w.r.t. the capillary pressure function and overall w.r.t. 
the viscosity ratios. The method is parameter free (the only parameter $r$ has a rather low influence 
and is chosen equal to 1 in the computations) and is unconditionally converging whatever the time step. 
This is a great advantage when compared to the Newton method that may require very small time steps 
in presence of large viscosity ratios. Moreover, the ALG2-JKO scheme preserves the positivity 
of the saturations, the constraint on the sum of the saturations, and it is locally conservative. 
Its main drawback concerns the restriction to linear mobility function so that formulas~\eqref{eq: dual Phi}--\eqref{eq:def_K_alpha} 
hold (this can probably be extended to the non-physical case of concave mobilities~\cite{DNS09} 
but we did not push into this direction). Finally, let us stress that the code depends only at 
stage~\eqref{JKO step2 second subproblem} of the energy. Therefore, the extension of the ALG2-JKO 
approach to multiphase models with different energies 
(like for instance degenerate Cahn-Hilliard models~\cite{OE97, CMN_HAL}) is not
demanding once the code is written. 
A natural extension to this work would be to add source terms corresponding for instance to production wells.
This would for instance require to adapt the material of~\cite{GLM_ArXiv} to our context. 

\subsection*{Acknowledgements}
CC was supported by the French National Research Agency (ANR) through grant 
ANR-13-JS01-0007-01 (project GEOPOR) and ANR-11-LABX0007-01 (Labex CEMPI).
LM was partially supported by the Portuguese Science Fundation through FCT grant PTDC/MAT-STA/0975/2014. 
TOG was partially supported by the Fonds de la Recherche Scientifique - FNRS under Grant MIS F.4539.16.

\end{document}